\documentclass{amsart}
\usepackage{amssymb}

\numberwithin{equation}{section}
\input{tcilatex}

\begin{document}
\title[Strong maximum principle]{Strong maximum principle for mean curvature
operators on subriemannian manifolds}
\author{Jih-Hsin Cheng}
\address{Institute of Mathematics, Academia Sinica, Taipei and National
Center for Theoretical Sciences, Taipei Office, Taiwan, R.O.C.}
\email{cheng@math.sinica.edu.tw}
\thanks{}
\author{Hung-Lin Chiu}
\address{Department of Mathematics, National Central University, Chung Li,
32054, Taiwan, R.O.C.}
\email{hlchiu@math.ncu.edu.tw}
\urladdr{}
\author{Jenn-Fang Hwang}
\address{Institute of Mathematics, Academia Sinica, Taipei, Taiwan, R.O.C.}
\email{majfh@math.sinica.edu.tw}
\thanks{}
\author{Paul Yang}
\address{Department of Mathematics, Princeton University, Princeton, NJ
08544, U.S.A.}
\email{yang@Math.Princeton.EDU}
\urladdr{}
\thanks{}

\begin{abstract}
We study the strong maximum principle for horizontal ($p$-) mean curvature
operator and $p$-(sub)laplacian operator on subriemannian manifolds
including, in particular, Heisenberg groups and Heisenberg cylinders. Under
a certain Hormander type condition on vector fields, we show the strong
maximum principle holds in higher dimensions for two cases: (a) the touching
point is nonsingular; (b) the touching point is an isolated singular point
for one of comparison functions. For a background subriemannian manifold
with local symmetry of isometric translations, we have the strong maximum
principle for associated graphs which include, among others, intrinsic
graphs with constant horizontal ($p$-) mean curvature. As applications, we
show a rigidity result of horizontal ($p$-) minimal hypersurfaces in any
higher dimensional Heisenberg cylinder and a pseudo-halfspace theorem for
any Heisenberg group.
\end{abstract}

\maketitle




\bigskip



\section{\textbf{Introduction and statement of the results}}

E. Hopf probably is the first person who studied the strong maximum
principle (SMP in short) of elliptic operators in the generality. See his
paper \cite{Ho} of 1927 or Theorem 3.5 in \cite{GT}. For earlier results,
under more restrictive hypotheses, see references in \cite{PW}. This
principle has been extended to the case for certain quasilinear elliptic
operators of second order (\cite{GT}). In 1969 J.-M. Bony (\cite{Bony})
studied, among others, the SMP for linear operators of H\"{o}rmander type
including some known subelliptic operators. Bony's SMP has been applied to
study various geometric problems. See, for instance, B. Andrews' work on
noncollapsing in mean-convex mean curvature flow (\cite{An}) or S. Brendle's
solution to the Lawson conjecture (\cite{Br}). In Subsection A of the
Appendix, we give a brief review of Bony's SMP.

In this paper we first extend Bony's SMP to the quasilinear case and then
apply it to (generalized) mean curvature operators in subriemannian
geometry, including $p$-sublaplacian and usual horizontal ($p$-) mean
curvature. We consider quasilinear operators $Q$ of second order:%
\begin{equation}
Q\phi =a^{ij}(x,D\phi )D_{ij}\phi +b(x,\phi ,D\phi )  \label{I1}
\end{equation}

\noindent where $x=(x^{1},..,x^{m+1})$ is contained in a domain $\Omega $ of 
$R^{m+1},$ $m\geq 1$. The coefficients $a^{ij}(x,p)$ ($b(x,z,p),$ resp.) of $%
Q$ are assumed to be defined and $C^{\infty }$ smooth (for simplicity) for
all values of $(x,p)$ (($x,z,p),$ resp.) in the set $\Omega \times R^{m+1}$ (%
$\Omega \times R\times R^{m+1},$ resp.). Let $\pounds (X_{1},...,X_{r})$
denote the smallest $C^{\infty }$-module which contains $C^{\infty }$ vector
fields $X_{1},$ $...,$ $X_{r}$ on $\Omega $ and is closed under the Lie
bracket (see (\ref{A4}) in the Appendix for precise definition). The
following comparison principle is a straightforward application of Bony's
SMP (Theorem 3.1 in \cite{Bony} or see Theorem A7 in the Appendix).

\bigskip

\textbf{Theorem A.}\textit{\ Let }$\phi ,$\textit{\ }$\psi $\textit{\ }$\in $%
\textit{\ }$C^{\infty }(\Omega )$\textit{\ satisfy }$Q\phi \geq Q\psi $%
\textit{\ in }$\Omega .$\textit{\ Assume}

\textit{(1) }$(a^{ij})$\textit{\ is nonnegative and }$a^{ij}=a^{ji};$

\textit{(2) }$\frac{\partial b}{\partial z}\leq 0;$

\textit{(3) there exist vector fields }$X_{1},..,X_{r}$\textit{\ and }$Y$%
\textit{\ of class }$C^{\infty }$\textit{\ (depending on }$D\phi (x))$ 
\textit{such that}%
\begin{equation*}
a^{ij}(x,D\phi (x))D_{ij}=\sum_{k=1}^{r}X_{k}^{2}+Y.
\end{equation*}

\textit{\noindent Let }$\Gamma $\textit{\ be an integral curve of a vector
field }$Z$\textit{\ }$\in $\textit{\ }$\pounds (X_{1},...,X_{r}).$\textit{\
Suppose }$\phi -\psi $\textit{\ achieves a nonnegative maximum in }$\Omega $%
\textit{\ at a point of }$\Gamma $\textit{. Then the maximum is attained at
all points of }$\Gamma .$

\bigskip

Let $e_{1},$ ..., $e_{m+1}$ be independent $C^{\infty }$ vector fields on $%
\Omega .$ Consider second order quasilinear operators $Q^{\prime }$ of the
form:%
\begin{equation}
Q^{\prime }\phi =a^{ij}(x,e_{1}\phi ,..,e_{m+1}\phi )e_{i}e_{j}\phi
+b(x,\phi ,e_{1}\phi ,..,e_{m+1}\phi )  \label{I2}
\end{equation}

\noindent where $x=(x^{1},..,x^{m+1})$ is contained in a domain $\Omega $ of 
$R^{m+1},$ $m\geq 1.$ We have the following moving frame version of Theorem
A.

\bigskip

\textbf{Theorem A}$^{\prime }$\textbf{.}\textit{\ Let }$\phi ,$\textit{\ }$%
\psi $\textit{\ }$\in $\textit{\ }$C^{\infty }(\Omega )$\textit{\ satisfy }$%
Q^{\prime }\phi \geq Q^{\prime }\psi $\textit{\ in }$\Omega .$\textit{\
Assume}

\textit{(1) }$(a^{ij})$\textit{\ is nonnegative and }$a^{ij}=a^{ji};$

\textit{(2) }$\frac{\partial b}{\partial z}\leq 0;$

\textit{(3) there exist vector fields }$X_{1},..,X_{r}$\textit{\ and }$Y$%
\textit{\ of class }$C^{\infty }$\textit{\ (depending on }$e_{1}\phi
(x),..,e_{m+1}\phi (x))$ \textit{such that}%
\begin{equation*}
a^{ij}(x,e_{1}\phi (x),..,e_{m+1}\phi
(x))e_{i}e_{j}=\sum_{k=1}^{r}X_{k}^{2}+Y.
\end{equation*}

\textit{\noindent Let }$\Gamma $\textit{\ be an integral curve of a vector
field }$Z$\textit{\ }$\in $\textit{\ }$\pounds (X_{1},...,X_{r}).$\textit{\
Suppose }$\phi -\psi $\textit{\ achieves a nonnegative maximum in }$\Omega $%
\textit{\ at a point of }$\Gamma $\textit{. Then the maximum is attained at
all points of }$\Gamma .$

\bigskip

Let%
\begin{eqnarray*}
\tilde{a}^{ij}(x) &:&=a^{ij}(x,D\phi (x))\text{ in Theorem A} \\
( &:&=a^{ij}(x,e_{1}\phi (x),..,e_{m+1}\phi (x))\text{ in Theorem A}^{\prime
},\text{ resp.).}
\end{eqnarray*}%
\textit{\noindent }In practice, condition (3) in Theorem A (Theorem A$%
^{\prime },$ resp.) can be replaced by 
\begin{equation}
rank(\tilde{a}^{ij}(x))=\text{constant }\tilde{r}  \label{I2-0}
\end{equation}

\bigskip

\textbf{Theorem \~{A}. }\textit{Theorem A (Theorem A}$^{\prime },$\textit{\
resp.) holds if condition (3) is replaced by constant rank condition (\ref%
{I2-0}).}

\bigskip

In applications, we usually assume $\phi \leq \psi $ and $\phi =\psi $ at a
point $p_{0}$. Then we conclude $\phi \equiv \psi $ on a hypersurface $%
\Sigma $ containing $p_{0}$ if the Lie span $\pounds (X_{1},...,X_{r})$ = $%
C^{\infty }(\Sigma ,T\Sigma ).$ We are going to apply Theorem \~{A} to
generalized mean curvature or $p$-laplacian $H_{\phi ,p}$ with $p\geq 0.$

A subriemannian manifold is a ($C^{\infty })$ smooth manifold $M$ equipped
with a nonnegative inner product $<\cdot ,\cdot >^{\ast }$ on $T^{\ast }M,$
its cotangent bundle, i.e., $<v,v>^{\ast }$ $\geq $ $0$ for any cotangent
vector $v$. When $<\cdot ,\cdot >^{\ast }$\ is positive definite, the
definition is equivalent to the usual definition of Riemannian manifold
using positive definite inner product on tangent bundle $TM$. However for a
degenerate $<\cdot ,\cdot >^{\ast },$ it is difficult to define on the whole 
$TM$ instead of $T^{\ast }M$. So the above definition using $T^{\ast }M$ of
subriemannian manifold generalizes the notion of Riemannian manifold in a
unified way.

Let $\phi $ be a ($C^{\infty }$ smooth, say) defining function (i.e., $d\phi 
$ $\neq $ $0)$ on a subriemannian manifold $(M,<\cdot ,\cdot >^{\ast })$ of
dimension $m+1$ with $|d\phi |_{\ast }$ $:=$ ($<d\phi ,d\phi >^{\ast
})^{1/2} $ $\neq $ $0$ (note that at a point where $<\cdot ,\cdot >^{\ast }$
is degenerate, we may have $d\phi $ $\neq $ $0$ while $|d\phi |_{\ast }$ $=$ 
$0).$ Let $dv_{M}$ be a background volume form on $M,$ i.e., a given $(m+1)$%
-form which is nowhere vanishing$.$ For $M$ being Riemannian, we may take
the associated volume form as a background volume form. For degenerate $M,$
a background volume form is a choice independent of the (degenerate) metric $%
<\cdot ,\cdot >^{\ast }.$ With respect to $<\cdot ,\cdot >^{\ast }$ we can
then talk about the interior product of a $1$-form $\omega $ with $dv_{M}:$ $%
\omega $ $\lrcorner $ $dv_{M}$ being an $m$-form such that $\eta \wedge
(\omega \lrcorner $ $dv_{M})=<\eta ,\omega >^{\ast }dv_{M}$ for any $1$-form 
$\eta .$

We define a function $H_{\phi ,p}$ on $M$ with $p\geq 0$ by the following
formula:%
\begin{equation}
d(\frac{d\phi }{|d\phi |_{\ast }^{1-p}}\lrcorner dv_{M})=H_{\phi ,p}dv_{M}.
\label{I2-0-1}
\end{equation}%
\textit{\noindent }For $p$ $=$ $0,$ $H_{\phi ,p},$ denoted as $H_{\phi }$
often, is called (Riemannian, subriemannian, or horizontal) mean curvature
while, for $p$ $>$ $0,$ $H_{\phi ,p}$ is so called $p$-laplacian or $p$%
-sublaplacian. For the variational formulation, consider the following
energy functional:%
\begin{equation*}
\mathcal{F}_{p}(\phi ):=\int_{\Omega }(\frac{1}{p}|d\phi |_{\ast }^{p}+H\phi
)dv_{M}
\end{equation*}%
\textit{\noindent }($H$ being prescribed subriemannian mean curvature or $p$%
-sublaplacian) where $\Omega $ $\subset $ $M$ is a smooth bounded domain.
Let $\phi _{\varepsilon }$ $=$ $\phi +\varepsilon \rho $ where $\rho $ $\in $
$C_{0}^{\infty }(\Omega ).$ Compute the first variation of $\mathcal{F}_{p}:$
(omitting the volume form $dv_{M})$%
\begin{eqnarray}
&&\frac{d\mathcal{F}_{p}(\phi _{\varepsilon })}{d\varepsilon }|_{\varepsilon
=0\pm }  \label{I2-0-2} \\
&=&c_{p}\int_{S(\phi )}|d\rho |_{\ast }^{p}+\int_{\Omega \backslash S(\phi
)}|d\phi |_{\ast }^{p-2}<d\phi ,d\rho >^{\ast }+\int_{\Omega }H\rho  \notag
\end{eqnarray}%
\textit{\noindent }where $c_{p}$ $=$ $\pm 1$ for $p$ $=$ $1$, $c_{p}$ $=$ $0$
for $1$ $<$ $p$ $<$ $\infty ,$ and $S(\phi )$ is the set where $|d\phi
|_{\ast }$ $=$ $0,$ called the singular set of $\phi $ (cf. (1.4) in \cite%
{chy}). From (\ref{I2-0-2}) we learn that for $p$ $=$ $1,$ the first term
involving the singular set $S(\phi )$ is not negligible. So in the proof of
the maximum principle (or comparison theorem), we need to worry about the
size of $S(\phi )$ (see \cite{chmy}, \cite{chy} for more details). In Lemma
5.2 of this paper, we extend the maximum principle (comparison theorem) to
general subriemannian manifolds. This is necessary in order to show the SMP
near singular points (where $|d\phi |_{\ast }$ $=$ $0)$. In this paper we
mainly deal with the SMP near nonsingular points (where $|d\phi |_{\ast }$ $%
\neq $ $0,$ $|d\psi |_{\ast }$ $\neq $ $0).$ For the SMP near singular
points, we only discuss the situation that the reference singular point is
isolated for at least one comparison hypersurface. In general, the problem
of the SMP near singular points is still open.

Let $\phi $ and $\psi $ be defining functions for hypersurfaces $\Sigma _{1}$
and $\Sigma _{2}$ in a subriemannian manifold $(M,$ $<\cdot ,\cdot >^{\ast
}) $ of dimension $m+1,$ resp. ($m\geq 1$). I.e., $\Sigma _{1}$ ($\Sigma
_{2},$ resp.) is defined by $\phi $ $=$ $0$ ($\psi $ $=$ $0$, resp.).
Suppose $\Sigma _{1}$ and $\Sigma _{2}$ are tangent to each other at a point 
$p_{0}$ where $|d\phi |_{\ast }$ $\neq $ $0,$ $|d\psi |_{\ast }$ $\neq $ $0.$
Define $G$ $:$ $T^{\ast }M$ $\rightarrow $ $TM$ by $\omega (G(\eta ))$ $=$ $%
<\omega ,\eta >^{\ast }$ for $\omega ,$ $\eta $ $\in $ $T^{\ast }M.$ Let $%
\xi $ $:=$ $Range(G).$ Throughout this paper we assume 
\begin{equation}
\dim \xi =\text{constant }m+1-l  \label{I2-1}
\end{equation}%
\textit{\noindent }near $p_{0}$ with $l$ being a nonnegative integer unless
stated otherwise. We call $l$ the degree of degeneracy of $M.$ Note that $%
\dim \ker G$ $=$ $l.$ The following rank condition:%
\begin{equation}
rank(\pounds (X_{1},...,X_{m-l}))=m  \label{I2-2}
\end{equation}%
\textit{\noindent }for any local sections $X_{1},...,X_{m-l}$ of $\xi ,$
which are independent wherever defined, is important. It means that any ($%
m-l)$-dimensional subspace of local sections of $\xi $ can generate $m$%
-dimensional spaces. Let 
\begin{equation*}
\pounds (\xi \cap T\Sigma _{1})=\pounds (X_{1},...,X_{m-l})
\end{equation*}%
\noindent where $X_{1},$..., $X_{m-l}$ form a basis of local sections of $%
\xi \cap T\Sigma _{1}$ near $p_{0}.$ Similarly we can define%
\begin{equation*}
\pounds (\xi )=\pounds (X_{1},...,X_{m-l},X_{m+1-l})
\end{equation*}

\noindent where $X_{1},...,X_{m-l},X_{m+1-l}$ form a basis of local sections
of $\xi $ near $p_{0}.$ Note that both $\pounds (\xi \cap T\Sigma _{1})$ and 
$\pounds (\xi )$ are independent of choice of a basis of local sections.

\bigskip

\textbf{Theorem B. }\textit{Suppose we are in the situation described above,
in particular,} $|d\phi |_{\ast }$ $\neq $ $0,$ $|d\psi |_{\ast }$ $\neq $ $%
0 $\textit{\ at }$p_{0}.$ \textit{For }$p\geq 0,$ \textit{Assume }%
\begin{equation*}
H_{\psi ,p}+b(x,\psi ,D\psi )\mathit{\ }\leq \mathit{\ }H_{\phi ,p}+b(x,\phi
,D\phi )
\end{equation*}%
\noindent \textit{where }$b$\textit{\ is a }$C^{\infty }$\textit{\ smooth
function for all values of }$(x,z,\cdot )$\textit{\ and satisfies }$\frac{%
\partial b}{\partial z}\leq 0.$\textit{\ Moreover, assume }$\psi $\textit{\ }%
$\geq $\textit{\ }$\phi $\textit{\ near }$p_{0},$\textit{\ }$\psi $\textit{\ 
}$=$\textit{\ }$\phi $\textit{\ }$=$\textit{\ }$0$\textit{\ at }$p_{0}.$%
\textit{\ We have}

\textit{(a) if we further assume the rank condition (}$\ref{I2-2}\mathit{)}$ 
\textit{holds near }$p_{0}$\textit{,\ then }$\psi $\textit{\ }$=$\textit{\ }$%
\phi $\textit{\ }$=$\textit{\ }$0$\textit{\ on }$\Sigma _{1}$\textit{\ near }%
$p_{0}.$\textit{\ I.e., }$\Sigma _{2}$\textit{\ coincides with }$\Sigma _{1}$%
\textit{\ near }$p_{0}.$

\textit{(b) in the case of }$p>0,$\textit{\ if we further assume }$rank(%
\pounds (\xi ))$\textit{\ }$=$\textit{\ }$m+1,$\textit{\ then }$\psi $%
\textit{\ }$=$\textit{\ }$\phi $\textit{\ near }$p_{0}.$

\bigskip

Next we want to show that in a certain situation the assumption in Theorem B
can be achieved. Suppose we have a one-parameter family of diffeomorphisms $%
\Psi _{a},$ $a$ $\in $ $(-\delta ,\delta )$ for small $\delta ,$ say, in a
small neighborhood $U$ of $p_{0}$ in $M$ (where $\phi $ and $\psi $ are
defined), i.e., $\Psi _{a}$ $=$ $Id$ for $a$ $=$ $0$ and $\Psi _{a+b}$ $=$ $%
\Psi _{a}\circ \Psi _{b}$ wherever defined$.$

\bigskip

\textbf{Definition 1.1.} \textit{Let }$(M,$\textit{\ }$<\cdot ,\cdot >^{\ast
},dv_{M})$\textit{\ be a subriemannian manifold with a background volume
form }$dv_{M}.$\textit{\ }$(M,$\textit{\ }$<\cdot ,\cdot >^{\ast },dv_{M})$%
\textit{\ is said to have isometric translations (}$\Psi _{a}\mathit{)}$ 
\textit{near }$p_{0}$\textit{\ }$\in $\textit{\ }$M$\textit{\ if there
exists a one-parameter family of local diffeomorphisms }$\Psi _{a},$\textit{%
\ }$a$\textit{\ }$\in $\textit{\ }$(-\delta ,\delta )$\textit{\ for small }$%
\delta >0,$\textit{\ say, in a small neighborhood }$U$\textit{\ of }$p_{0}$%
\textit{, such that}

\textit{(a) (preserving }$<\cdot ,\cdot >^{\ast })$\textit{\ }$<\Psi
_{a}^{\ast }(\omega ),\Psi _{a}^{\ast }(\eta )>^{\ast }=<\omega ,\eta
>^{\ast }$\textit{\ for }$a$\textit{\ }$\in $\textit{\ }$(-\delta ,\delta )$%
\textit{\ and }$\omega ,\eta $\textit{\ }$\in $\textit{\ }$T^{\ast }U.$

\textit{(b) (preserving }$dv_{M})$\textit{\ }$\Psi _{a}^{\ast
}(dv_{M})=dv_{M}$\textit{\ for }$a$\textit{\ }$\in $\textit{\ }$(-\delta
,\delta ).$

\bigskip

We say the defining function $\phi $ of a (local) hypersurface $\Sigma $
passing through $p_{0}$ is compatible with \{$\Psi _{a}\}$ or \{$\Psi _{a}\}$
is compatible with the defining function $\phi $ if 
\begin{equation}
\phi (\Psi _{a}(x))=\phi (x)-a  \label{I3}
\end{equation}%
\noindent for $x$ $\in $ $\Sigma $ $\cap $ $U$ (and hence $x$ $\in $ $U)$
and $a\in (-\delta ,\delta ).$ \{$\Psi _{a}\}$ is said to be transversal to
a hypersurface $\Sigma $ of $U$ if 
\begin{equation*}
\frac{d\Psi _{a}(x)}{da}|_{a=0}\notin T_{x}\Sigma
\end{equation*}

\noindent for all $x$ $\in $ $\Sigma .$ We have a notion of (generalized)
"hypersurface" mean curvature ($p$-sublaplacian) $H_{\Sigma ,p}$ defined on $%
\Sigma \backslash S(\phi )$ as follows. At a point where $|d\phi |_{\ast }$ $%
\neq $ $0,$ for $p\geq 0,$ we define $p$-subriemannian area (or volume)
element $dv_{\phi ,p}$ for the hypersurface \{$\phi $ $=$ $c,$ a constant$\}$
by%
\begin{equation*}
dv_{\phi ,p}:=\frac{d\phi }{|d\phi |_{\ast }^{1-p}}\lrcorner dv_{M}.
\end{equation*}

\noindent Define unit $p$-normal $\nu _{p}$ to a hypersurface $\Sigma $ $:=$ 
$\{\phi =c\}$ by the formula%
\begin{equation*}
\nu _{p}\text{ }\lrcorner \text{ }dv_{M}=dv_{\phi ,p}.
\end{equation*}

\noindent We can now define "hypersurface" mean curvature $H_{\Sigma ,p}$ on 
$\Sigma $ through a variational formula:%
\begin{equation*}
\delta _{f\nu _{p}}\int_{\Sigma }dv_{\phi ,p}=\int_{\Sigma }fH_{\Sigma
,p}dv_{\phi ,p}
\end{equation*}

\noindent for $f\in C_{0}^{\infty }(\Sigma \backslash S(\phi )).$ In
Subsection B of the Appendix, we show that $H_{\phi ,p}$ $=$ $H_{\Sigma ,p}$
on $\Sigma \backslash S(\phi )$ (see Proposition B.1).

\bigskip

\textbf{Theorem C}. \textit{Suppose} $(M,$\textit{\ }$<\cdot ,\cdot >^{\ast
},dv_{M})$\textit{\ of dimension }$m+1$ \textit{has isometric translations }$%
\Psi _{a}$ \textit{near }$p_{0}$\textit{\ }$\in $\textit{\ }$M.$ \textit{Let 
}$U$\textit{\ be the neighborhood of }$p_{0}$\textit{\ in Definition 1.1.} 
\textit{\ Suppose }$\{\Psi _{a}\}$\textit{\ is transversal to }$\Sigma _{1}$%
\textit{\ and }$\Sigma _{2}$ \textit{in }$U.$ \textit{Choose }$\phi $\textit{%
\ and }$\psi $\textit{\ to be defining functions for hypersurfaces }$\Sigma
_{1}$\textit{\ and }$\Sigma _{2}$\textit{\ in }$U,$\textit{\ resp. (I.e., }$%
\Sigma _{1}$\textit{\ (}$\Sigma _{2},$\textit{\ resp.) is defined by }$\phi $%
\textit{\ }$=$\textit{\ }$0$\textit{\ (}$\psi $\textit{\ }$=$\textit{\ }$0$%
\textit{, resp.)), compatible with }$\{\Psi _{a}\}$\textit{. Suppose }$%
\Sigma _{1}$\textit{\ and }$\Sigma _{2}$\textit{\ are tangent to each other
at }$p_{0}$\textit{\ where }$|d\phi |_{\ast }$\textit{\ }$\neq $\textit{\ }$%
0,$\textit{\ }$|d\psi |_{\ast }$\textit{\ }$\neq $\textit{\ }$0$\textit{.\
Assume}

\textit{(1) For any }$q$\textit{\ }$\in $\textit{\ }$\Sigma _{1}$\textit{\ }$%
\cap $\textit{\ }$U,$\textit{\ there exists }$\delta $\textit{\ }$>$\textit{%
\ }$a(q)$\textit{\ }$\geq $\textit{\ }$0$\textit{\ such that }$\Psi
_{a(q)}(q)$\textit{\ }$\in $\textit{\ }$\Sigma _{2},$\textit{\ and}

\textit{(2) For some }$p\geq 0,$ $H_{\Sigma _{2},p}(\Psi _{a(q)}(q))$\textit{%
\ }$\leq $\textit{\ }$H_{\Sigma _{1},p}(q)$\textit{\ for any }$q$\textit{\ }$%
\in $\textit{\ }$\Sigma _{1}$\textit{\ }$\cap $\textit{\ }$U.$

\noindent \textit{Moreover, assume the rank conition (\ref{I2-2}) holds near 
}$p_{0}.$\textit{Then }$\Sigma _{2}$\textit{\ coincides with }$\Sigma _{1}$%
\textit{\ near }$p_{0}.$

\bigskip

Next we want to give a more analytic description in terms of graphs. We will
show the existence of some special coordinates for a subriemannian manifold
having local isometric translations as shown below.

\bigskip

\textbf{Theorem \^{C}}. \textit{Let }$(M,$\textit{\ }$<\cdot ,\cdot >^{\ast
},dv_{M})$\textit{\ be an }$(m+1)$-\textit{dimensional subriemannian
manifold with a background volume form }$dv_{M}.$ \textit{Suppose} $(M,$%
\textit{\ }$<\cdot ,\cdot >^{\ast },dv_{M})$ \textit{has isometric
translations }$\Psi _{a},$ $a$\textit{\ }$\in $\textit{\ }$(-\delta ,\delta
) $ for $\delta >0,$ \textit{near }$p_{0}$\textit{\ }$\in $\textit{\ }$M,$ 
\textit{transversal to a hypersurface }$\Sigma $\textit{\ passing through }$%
p_{0}.$\textit{\ Then we can find local coordinates }$x^{1},$\textit{\ }$%
x^{2},$\textit{\ ..., }$x^{m+1}$\textit{\ in a neighborhood }$V$\textit{\ of 
}$p_{0}$\textit{\ such that}

\textit{(1) }$\Sigma $\textit{\ is described by }$x^{m+1}$\textit{\ }$=$%
\textit{\ }$0$\textit{\ in }$V,$\textit{\ }$p_{0}$\textit{\ is the origin},%
\textit{\ }$x^{j}\circ \Psi _{a}$ $=$ $x^{j}$\textit{\ at }$q$ $\in $ $V$ 
\textit{for any }$a$ \textit{such that }$\Psi _{a}(q)$\textit{\ }$\in $%
\textit{\ }$V,$\textit{\ }$1\leq j\leq m,$ \textit{and }$x^{m+1}\circ \Psi
_{a}$\textit{\ }$=$\textit{\ }$x^{m+1}+a$\textit{\ }$;$

\textit{(2) }$A(\Psi _{a}(q))$\textit{\ }$=$\textit{\ }$A(q)$\textit{\ for
any }$q\in V$\textit{\ and any }$a$ \textit{such that }$\Psi _{a}(q)$\textit{%
\ }$\in $\textit{\ }$V$\textit{\ if we write }$dv_{M}(q)$\textit{\ }$=$%
\textit{\ }$A(q)dx^{1}\wedge ...\wedge dx^{m+1}.$

\bigskip

Call this system of special coordinates above in Theorem \^{C} (a system of)
translation-isometric coordinates.

\bigskip

\textbf{Definition 1.2}. \textit{Take a system of translation-isometric
coordinates }$x^{1},$\textit{\ }$x^{2},$\textit{\ ..., }$x^{m+1}$\textit{\
as in Theorem \^{C}}$.$\textit{\ A }$\Psi _{a}$\textit{\ graph is a graph
described by} 
\begin{equation*}
(x^{1},\mathit{\ }x^{2},\mathit{\ ..,}\text{ }x^{m},\text{ }u(x^{1},\mathit{%
\ }x^{2},\mathit{\ ..,}\text{ }x^{m})).
\end{equation*}

\bigskip

For a $\Psi _{a}$\textit{\ }graph, we take the defining function $\phi $ $=$ 
$u(x^{1},$\textit{\ }$x^{2},$\textit{\ ...}$x^{m})$ $-$ $x^{m+1}$ which is
compatible with \{$\Psi _{a}\}.$ For $p$ $\geq $ $0,$ we define 
\begin{equation}
H_{p}(u)(x^{1},..,x^{m}):=H_{\phi ,p}(x^{1},..,x^{m},u(x^{1},\mathit{\ }%
x^{2},\mathit{\ ...}x^{m}))  \label{I4}
\end{equation}%
\textit{\noindent }at $(x^{1},..,x^{m})$ where $|d\phi |_{\ast }$ $\neq $ $%
0. $ Making use of translation-isometric coordinates, we can reformulate
Theorem C as follows.

\bigskip

\textbf{Theorem C}$^{\prime }.$ \textit{Suppose} $(M,$\textit{\ }$<\cdot
,\cdot >^{\ast },dv_{M})$\textit{\ of dimension }$m+1$ \textit{has isometric
translations }$\Psi _{a}$ \textit{near }$p_{0}$\textit{\ }$\in $\textit{\ }$%
M,$ \textit{transversal to a hypersurface }$\Sigma $\textit{\ passing
through }$p_{0}.$ \textit{Take a system of translation-isometric coordinates}
$x^{1},$\textit{\ }$x^{2},$\textit{\ ..., }$x^{m+1}$\textit{\ in a
neighborhood }$V$\textit{\ of }$p_{0}$ \textit{such that }$x^{m+1}$\textit{\ 
}$=$\textit{\ }$0$\textit{\ on }$\Sigma .$\textit{\ Suppose }$u$\textit{\ }$%
(v,$\textit{\ resp.) }$:$\textit{\ }$\Sigma \cap V\rightarrow R$\textit{\
defines a graph }$\{(x^{1},$\textit{\ }$x^{2},$\textit{\ ...}$x^{m},$\textit{%
\ }$u(x^{1},$\textit{\ }$x^{2},$\textit{\ ...}$x^{m}))\}$\textit{\ (}$%
\{(x^{1},$\textit{\ }$x^{2},$\textit{\ ...}$x^{m},$\textit{\ }$v(x^{1},$%
\textit{\ }$x^{2},$\textit{\ ...}$x^{m}))\}$\textit{, resp.) }$\subset $%
\textit{\ }$V$ \textit{such that} $|d(u-x^{m+1})|_{\ast }$ $\neq $ $0$ 
\textit{(}$|d(v-x^{m+1})|_{\ast }$ $\neq $ $0,$\textit{\ resp.). Assume }

\textit{(1) }$v\geq u$\textit{\ on }$\Sigma \cap V$\textit{\ and }$v(0,..,0)$%
\textit{\ }$=$\textit{\ }$u(0,..,0)$ $=$ $0;$

\textit{(2) For some }$p\geq 0,$ $H_{p}(v)\leq H_{p}(u)$\textit{\ on }$%
\Sigma \cap V.$

\textit{Moreover, assume the rank condition (\ref{I2-2}) holds near }$p_{0}.$%
\textit{Then }$v$ $\equiv $ $u$\textit{\ in a neighborhood of }$p_{0}$ $\in $
$\Sigma .$

\bigskip

Let $M$ be the Heisenberg group $H_{n}$ considered as a pseudohermitian
manifold and hence a subriemannian manifold (see Appendix B for detailed
explanation). Suppose two hypersurfaces $\Sigma _{1}$ and $\Sigma _{2}$ in
Theorem C are (horizontal) graphs over the $x^{1}x^{2}...x^{2n}$ hyperplane,
defined by $v$ and $u,$ resp.. We can take $\Psi _{a}$ to be the translation
in the last coordinate by the amount $a$ : $\Psi _{a}(x^{1},$ $x^{2},$ $...,$
$x^{2n-1},$ $x^{2n},$ $z)$ $=$ $(x^{1},$ $x^{2},$ $...,$ $x^{2n-1},$ $%
x^{2n}, $ $z+a).$ The defining functions $\psi $ and $\phi $ for $\Sigma
_{2} $ and $\Sigma _{1},$ resp. are given by $%
v(x^{1},x^{2},...,x^{2n-1},x^{2n})$ $-$ $z$ and $%
u(x^{1},x^{2},...,x^{2n-1},x^{2n})$ $-$ $z$. We can then verify the
assumption of Theorem C or Theorem C$^{\prime }$ ($m$ $=$ $2n)$ and that
condition (2) for the case $p$ $=$ $0$ is equivalent to $H_{\vec{F}}(v)$%
\textit{\ }$\leq $\textit{\ }$H_{\vec{F}}(u)$ (see (\ref{0.1}) below) by
identifying $H_{\vec{F}}(v)$ and $H_{\vec{F}}(u)$ with the mean curvature of 
$\psi $ and $\phi ,$ resp., with respect to a certain subriemannian manifold 
$(M,$\textit{\ }$<\cdot ,\cdot >^{\ast },dv_{M})$ while condition (1) is the
same$.$ So Theorem C or Theorem C$^{\prime }$ includes the Heisenberg group
case (see Theorem F below and its proof in Section 4 for more details).

Another situation is that two hypersurfaces $\Sigma _{1}$ and $\Sigma _{2}$
are tangent at $p_{0}$ vertically in $H_{n},$ i.e., the common tangent space
at $p_{0}$ is a hyperplane $E$ perpendicular to the $x^{1}x^{2}...x^{2n}$
hyperplane. We can then find a one-parameter family of Heisenberg
translations in a direction of vector normal to this tangent space at $%
p_{0}. $ Let $l_{a}$ denote the Heisenberg translation in the direction $%
\frac{\partial }{\partial x^{1}}:$ 
\begin{equation}
l_{a}(x^{1},x^{2},...,x^{2n-1},x^{2n},z):=(x^{1}+a,x^{2},...,x^{2n-1},x^{2n},z-ax^{n+1}).
\label{I5}
\end{equation}

\bigskip

\textbf{Corollary D}. \textit{Suppose }$\Sigma _{1}$\textit{\ and }$\Sigma
_{2}$\textit{\ are tangent at }$p_{0}$ \textit{vertically in }$H_{n}$\textit{%
\ with }$n$\textit{\ }$\geq $\textit{\ }$2.$\textit{\ Suppose the
(Euclidean) unit normal to the tangent space at }$p_{0}$\textit{\ is }$-%
\frac{\partial }{\partial x^{1}}$\textit{\ without loss of generality}$.$%
\textit{\ For }$p\geq 0,$\textit{\ we} \textit{assume }%
\begin{equation*}
H_{\Sigma _{2},p}(l_{a(q)}(q))\leq H_{\Sigma _{1},p}(q)
\end{equation*}%
\textit{\noindent for }$q$\textit{\ }$\in $\textit{\ }$\Sigma _{1}$\textit{\
near }$p_{0}$ \textit{and }$l_{a(q)}(q)$\textit{\ }$\in $\textit{\ }$\Sigma
_{2}$\textit{\ with }$a(q)$\textit{\ }$\geq $\textit{\ }$0.$\textit{\ Then }$%
\Sigma _{2}$\textit{\ coincides with }$\Sigma _{1}$\textit{\ near }$p_{0}.$

\bigskip

As an example of Theorem C$^{\prime },$ we have the SMP for (horizontal)
mean curvature of $l_{a}$ graphs. See Corollary D$^{\prime }$ in Section 3.
Observe that translation-isometric coordinates with respect to the
Heisenberg translation in the direction $\frac{\partial }{\partial x_{1}}$
are closedly related to coordinates for an intrinsic graph. Recall (\cite%
{ASCV}) that an intrinsic graph $u$ in $H_{n}$ is a hypersurface of the form 
$(0,$ $\eta ^{2},$ $\eta ^{3},$ $..,$ $\eta ^{2n},$ $\tau )\circ (u(\eta
^{2},$ $\eta ^{3},$ $..,$ $\eta ^{2n},$ $\tau ),$ $0,$ $..,$ $0).$ Namely,
it is parametrized by $\eta ^{2},$ $\eta ^{3},$ $..,$ $\eta ^{2n},$ $\tau $
so that%
\begin{eqnarray*}
x^{1} &=&u(\eta ^{2},\eta ^{3},..,\eta ^{2n},\tau ), \\
x^{2} &=&\eta ^{2},...,x^{2n}=\eta ^{2n}, \\
z &=&\tau +\eta ^{n+1}u(\eta ^{2},\eta ^{3},..,\eta ^{2n},\tau )
\end{eqnarray*}%
\textit{\noindent }(see the Appendix for the definition of multiplication $%
\circ $ in $H_{n}).$ Let 
\begin{eqnarray*}
\mathring{e}_{2} &:&=\frac{\partial }{\partial \eta ^{2}}+\eta ^{n+2}\frac{%
\partial }{\partial \tau },..,\mathring{e}_{n}:=\frac{\partial }{\partial
\eta ^{n}}+\eta ^{2n}\frac{\partial }{\partial \tau }, \\
\mathring{e}_{n+1}^{u} &:&=\frac{\partial }{\partial \eta ^{n+1}}-2u\frac{%
\partial }{\partial \tau }, \\
\mathring{e}_{n+2} &:&=\frac{\partial }{\partial \eta ^{n+2}}-\eta ^{2}\frac{%
\partial }{\partial \tau },..,\mathring{e}_{2n}:=\frac{\partial }{\partial
\eta ^{2n}}-\eta ^{n}\frac{\partial }{\partial \tau }.
\end{eqnarray*}%
\textit{\noindent }Define a vector-valued operator $W^{u}$ by%
\begin{equation*}
W^{u}:=(\mathring{e}_{2},..,\mathring{e}_{n},\mathring{e}_{n+1}^{u},%
\mathring{e}_{n+2},..,\mathring{e}_{2n}).
\end{equation*}%
\textit{\noindent }The horizontal (or $p$-)mean curvature of an intrinsic
graph $u$ is given by%
\begin{equation}
H_{u(\eta ^{2},\eta ^{3},..,\eta ^{2n},\tau )-x^{1}}=W^{u}\cdot \left( \frac{%
W^{u}(u)}{\sqrt{1+|W^{u}(u)|^{2}}}\right)  \label{I6}
\end{equation}%
\textit{\noindent }at ($\eta ^{2},\eta ^{3},..,\eta ^{2n},\tau ,$ $u(\eta
^{2},\eta ^{3},..,\eta ^{2n},\tau ))$ (see (\ref{3.39}))$.$ In Section 3 we
provide more details and observe that an intrinsic graph is congruent with
an $l_{a}$ graph by a rotation in a certain situation. We therefore have the
SMP for intrinsic graphs with constant horizontal (or $p$-)mean curvature as
a special case of Theorem C$^{\prime }$.

\bigskip

\textbf{Theorem E.} \textit{Suppose }$n\geq 2.$ \textit{Let }$v$\textit{\ }$%
= $\textit{\ }$v(\eta ^{2},\eta ^{3},..,\eta ^{2n},\tau ),$\textit{\ }$u$%
\textit{\ }$=$\textit{\ }$u(\eta ^{2},\eta ^{3},..,\eta ^{2n},\tau )$\textit{%
\ be two (}$C^{\infty }$ \textit{smooth)} \textit{intrinsic graphs defined
on a neighborhood }$U$\textit{\ of }$p_{0}$\textit{\ }$=$\textit{\ }$(\eta
_{0}^{2},\eta _{0}^{3},..,\eta _{0}^{2n},\tau _{0}).$\textit{\ Assume }$v$%
\textit{\ }$=$\textit{\ }$u$\textit{\ at }$p_{0}$ \textit{where }$\nu _{\eta
^{n+1}}$\textit{\ }$\neq $\textit{\ }$0$\textit{\ and }$u_{\eta ^{n+1}}$%
\textit{\ }$\neq $\textit{\ }$0.$\textit{\ Suppose }$v$\textit{\ }$\geq $%
\textit{\ }$u$\textit{,}%
\begin{equation*}
W^{v}\cdot \left( \frac{W^{v}(v)}{\sqrt{1+|W^{v}(v)|^{2}}}\right) \leq
W^{u}\cdot \left( \frac{W^{u}(u)}{\sqrt{1+|W^{u}(u)|^{2}}}\right)
\end{equation*}%
\textit{\noindent in }$U,$\textit{\ and either }$v$\textit{\ or }$u$\textit{%
\ has constant}$.$\textit{horizontal (or }$p$\textit{-)mean curvature.\ Then 
}$v\equiv u$\textit{\ near }$p_{0}.$

\bigskip

We remark that the horizontal mean curvature operator of an intrinsic graph $%
u$ does not belong to the type (\ref{I1}) since the second order
coefficients contain $u$ itself. So Theorem E does not follow directly from
previous general theorems.

Let $\Omega $ be a (connected and open) domain of $R^{m}.$ Let $u,$ $v$ be
two $C^{2}$ smooth, real valued functions on $\Omega .$ Let $\vec{F}$ be a $%
C^{1}$ smooth vector field on $\Omega .$ Define the Legendrian (or
horizontal) normal $N_{\vec{F}}(u)$ ($N_{\vec{F}}(v),$ resp.) of $u$ ($v,$
resp.) by%
\begin{equation*}
N_{\vec{F}}(u):=\frac{\nabla u+\vec{F}}{|\nabla u+\vec{F}|}
\end{equation*}

\noindent ($N_{\vec{F}}(v):=\frac{\nabla v+\vec{F}}{|\nabla v+\vec{F}|},$
resp.) at points where $\nabla u+\vec{F}$ $\neq $ $0$ ($\nabla v+\vec{F}$ $%
\neq $ $0,$ resp.)$.$ Define the (generalized) horizontal (or $p$-) mean
curvature $H_{\vec{F}}(u)$ ($H_{\vec{F}}(v),$ resp.) by%
\begin{equation}
H_{\vec{F}}(u):=\func{div}N_{\vec{F}}(u)  \label{0.1}
\end{equation}

\noindent ($H_{\vec{F}}(v):=\func{div}N_{\vec{F}}(v),$ resp.). We call a
point $p_{0}$ singular with respect to $v$ if $\nabla v+\vec{F}$ $=$ $0$ at $%
p_{0}.$ Denote the set of all singular points with respect to $v$ by $S_{%
\vec{F}}(v).$

In a neighborhood $U$ of a nonsingular point $q_{0}$ $\in $ $\Omega
\backslash S_{\vec{F}}(v),$ let $N_{1}^{\perp }(v),$ $N_{2}^{\perp }(v),$
..., $N_{m-1}^{\perp }(v)$ be an orthonormal basis of the space
perpendicular to $N_{\vec{F}}(v)$. Let $\pounds (N_{1}^{\perp }(v),$ $%
N_{2}^{\perp }(v),$ ..., $N_{m-1}^{\perp }(v))$ denote the smallest $%
C^{\infty }$-module which contains $N_{1}^{\perp }(v),$ $N_{2}^{\perp }(v),$
..., $N_{m-1}^{\perp }(v),$ and is closed under the Lie bracket (see (\ref%
{A4}) in the Appendix for precise definition). The rank of $\pounds %
(N_{1}^{\perp }(v),$ $N_{2}^{\perp }(v),$ ..., $N_{m-1}^{\perp }(v))$ at a
point $q$ $\in $ $U$ is the dimension of the vector space spanned by the
vectors $Z(q)$ for all $Z$ $\in $ $\pounds (N_{1}^{\perp }(v),$ $%
N_{2}^{\perp }(v),$ ..., $N_{m-1}^{\perp }(v)).$ The following result is a
special, but important case of Theorem C or Theorem C$^{\prime }$ for degree
of degeneracy $l$ $=$ $1.$

\bigskip

\textbf{Theorem F.} \textit{Suppose }$m\geq 3,$ $H_{\vec{F}}(v)$\textit{\ }$%
\leq $\textit{\ }$H_{\vec{F}}(u),$\textit{\ }$v$ $\geq $ $u$\textit{\ in }$U$
$\subset $ $R^{m},$ \textit{which is a nonsingular domain for both }$v$%
\textit{\ and }$u,$\textit{\ and }$v=u$\textit{\ at }$p_{0}$\textit{\ }$\in $%
\textit{\ }$U.$\textit{\ Assume in }$U$ \textit{an orthonormal basis }$%
N_{1}^{\perp }(u),$ $N_{2}^{\perp }(u),$ ..., $N_{m-1}^{\perp }(u)$ \textit{%
of the space perpendicular to }$N_{\vec{F}}(u)$\textit{\ exists and} \textit{%
the rank of }$\pounds (N_{1}^{\perp }(u),$ $N_{2}^{\perp }(u),$ ..., $%
N_{m-1}^{\perp }(u))$ \textit{is constant} $m$ \textit{(similar condition
for }$N_{\vec{F}}(v),$\textit{\ resp.)}$.$ \textit{Then we have }$v\equiv u$ 
\textit{in} $U.$

\bigskip

Write 
\begin{equation}
N_{\alpha }^{\perp }(v)=\sum_{k=1}^{m}b_{\alpha }^{k}\partial _{k}.
\label{1.1}
\end{equation}

\textbf{Corollary G}. \textit{Suppose }$m\geq 3,$ $H_{\vec{F}}(u)$\textit{\ }%
$\leq $\textit{\ }$H_{\vec{F}}(v),$\textit{\ }$u$ $\geq $ $v$\textit{\ in }$%
U $ $\subset $ $R^{m},$ \textit{which is a nonsingular domain for both }$v$%
\textit{\ and }$u,$\textit{\ and }$u=v$\textit{\ at }$p_{0}$\textit{\ }$\in $%
\textit{\ }$U.$ \textit{Assume there exists a pair of (}$\alpha ,\beta ),$%
\textit{\ }$\alpha \neq \beta ,$\textit{\ such that}%
\begin{equation}
\sum_{k<j}(\partial _{k}F_{j}-\partial _{j}F_{k})(b_{\alpha }^{k}b_{\beta
}^{j}-b_{\alpha }^{j}b_{\beta }^{k})\neq 0  \label{1.2}
\end{equation}%
\noindent (\textit{similar condition for} $N_{\alpha }^{\perp }(u),$ $resp.)$
\textit{in }$U.$ \textit{Then we have }$u\equiv v$ \textit{in} $U.$

\bigskip

\textbf{Corollary H}. \textit{Suppose }$m\geq 3,$ $H_{\vec{F}}(u)$\textit{\ }%
$\leq $\textit{\ }$H_{\vec{F}}(v),$\textit{\ }$u$ $\geq $ $v$\textit{\ in }$%
U $ $\subset $ $R^{m},$ \textit{which is a nonsingular domain for both }$u$%
\textit{\ and }$v,$\textit{\ and }$u=v$\textit{\ at }$p_{0}$\textit{\ }$\in $%
\textit{\ }$U.$ \textit{Assume}%
\begin{equation}
rank(\partial _{k}F_{j}-\partial _{j}F_{k})\geq 3  \label{1.2.1}
\end{equation}

\noindent \textit{in }$U.$ \textit{Then we have }$u\equiv v$ \textit{in} $U.$

\bigskip

\textbf{Corollary I}. \textit{Suppose} $m$ $=$ $2n$ $\geq $ $4,$ $H_{\vec{F}%
}(u)$\textit{\ }$\leq $\textit{\ }$H_{\vec{F}}(v),$\textit{\ }$u$ $\geq $ $v$%
\textit{\ in }$U$ $\subset $ $R^{2n},$ \textit{which is a nonsingular domain
for both }$v$\textit{\ and }$u,$\textit{\ and }$u=v$\textit{\ at }$p_{0}$%
\textit{\ }$\in $\textit{\ }$U.$ Assume $\vec{F}$ $=$ $(-x^{2},$ $x^{1},$
..., $-x^{2n},$ $x^{2n-1}).$ \textit{Then we have }$u\equiv v$ \textit{in} $%
U.$

\bigskip

Corollary I provides the SMP of so called horizontal (or $p$-) mean
curvature for hypersurfaces given by graphs over a domain of the $%
x^{1}x^{2}...x^{2n}$ hyperplane in the Heisenberg group $H_{n}$ (identified
with $R^{2n+1}$ as a set).

We remark that when $m$ $=$ $2,$ Corollary I does not hold. That is, the SMP
of horizontal ($p$-)mean curvature for surfaces in $H_{1}$ does not hold
(although the maximum principle holds; see \cite{chmy} or Theorem $C^{\prime
\prime }$ below) as shown by the following example. Let $u$ $=$ $x^{1}x^{2}$ 
$+$ $(x^{2})^{2}$ and $v$ $=$ $x^{1}x^{2}.$ It follows from (\ref{0.1}) that 
$H_{\vec{F}}(u)$ $=$ $H_{\vec{F}}(v)$ $=$ $0$ in a nonsingular domain for $%
\vec{F}$ $=$ $(-x^{2},$ $x^{1}).$ Observe that $u$ $\geq $ $v$\textit{\ }in $%
U$ $=$ $\{x^{1}$ $>$ $0,$ $x^{1}+x^{2}$ $>$ $0\},$ a nonsingular domain for
both $u$ and $v,$ and $u$ $=$ $v$ at $(x^{1},$ $0)$ $\in $ $U.$ But
apparently $u$ $\neq $ $v$ in $U\backslash \{(x^{1},$ $0)$ $:$ $x^{1}$ $>$ $%
0\}.$

\bigskip

Define $\vec{G}^{b}$ for $\vec{G}$ $=$ $(G_{1},$ $...,$ $G_{m})$ by

\begin{equation*}
\vec{G}^{b}:=(\sum_{k=1}^{m}a^{1k}G_{k},\sum_{k=1}^{m}a^{2k}G_{k},...,%
\sum_{k=1}^{m}a^{mk}G_{k})
\end{equation*}

\noindent where $a^{jk\prime }s$ are real constants such that $a^{jk}+a^{kj}$
$=$ $0$ for $1$ $\leq $ $j,k$ $\leq $ $m.$ Note that $\vec{G}^{b}$ $=$ $\vec{%
G}^{\ast }$ for $m$ $=$ $2n,$ $a^{2j-1,2j}$ $=$ $-a^{2j,2j-1}$ $=$ $1,$ $%
1\leq j\leq n,$ $a^{jk}$ $=$ $0$ otherwise. When $p_{0}$ is an isolated
singular point of $v$, we still have the SMP.

\bigskip

\textbf{Theorem J.} \textit{Suppose }$m\geq 3,$\textit{\ }$v$ $\geq $ $u$%
\textit{\ in }$\Omega $ $\subset $ $R^{m},$\textit{\ such that }$\Omega \cap
S_{\vec{F}}(u)$\textit{\ }$=$\textit{\ }$\{p_{0}\}$\textit{\ (}$\Omega \cap
S_{\vec{F}}(v)$\textit{\ }$=$\textit{\ }$\{p_{0}\},$\textit{\ resp.) and }$%
v=u$\textit{\ at }$p_{0}.$ \textit{Suppose }$\mathcal{H}_{m-1}(\overline{S_{%
\vec{F}}(v)})$\textit{\ }$=$\textit{\ }$0$\textit{\ (}$\mathcal{H}_{m-1}(%
\overline{S_{\vec{F}}(u)})$\textit{\ }$=$\textit{\ }$0,$\textit{\ resp.) and 
}$\func{div}\vec{F}^{b}$\textit{\ }$>$\textit{\ }$0$\textit{\ (or }$\func{div%
}\vec{F}^{b}$\textit{\ }$<$\textit{\ }$0).$ \textit{Assume }$H_{\vec{F}}(v)$%
\textit{\ }$\leq $\textit{\ }$H_{\vec{F}}(u)$\textit{\ in }$\Omega
\backslash \{\{p_{0}\}\cup S_{\vec{F}}(v)\}$\textit{\ (}$\Omega \backslash
\{\{p_{0}\}\cup S_{\vec{F}}(u)\},$\textit{\ resp.) and for each point }$p$ $%
\in $ $\Omega \backslash \{p_{0}\},$ \textit{there is a neighborhood }$U$%
\textit{\ of }$p$\textit{\ in which an orthonormal basis }$N_{1}^{\perp
}(u), $\textit{\ }$N_{2}^{\perp }(u),$\textit{\ ..., }$N_{m-1}^{\perp }(u)$%
\textit{\ of the space perpendicular to }$N_{\vec{F}}(u)$\textit{.exists and
the rank of }$\pounds (N_{1}^{\perp }(u),$ $N_{2}^{\perp }(u),$ ..., $%
N_{m-1}^{\perp }(u))$ \textit{is constant} $m$ \textit{(similar condition
for }$N_{\vec{F}}(v),$\textit{\ resp.)}$.$ \textit{Then we have }$v\equiv u$ 
\textit{in} $\Omega .$

\bigskip

In the proof of Theorem J, we need to apply the following version of the
usual maximum principle (in the case of removable singularity) for $H_{\vec{F%
}}.$

\bigskip

\textbf{Theorem }$\mathbf{C}^{\prime \prime }$ (an extension of Theorem $%
C^{\prime }$ in \cite{chmy}).\textit{\ For a bounded domain }$\Omega $%
\textit{\ in }$R^{m}$\textit{\ with }$m$\textit{\ }$\geq $\textit{\ }$2,$%
\textit{\ let }$v,$\textit{\ }$u$\textit{\ }$\in $\textit{\ }$C^{2}(\Omega )$%
\textit{\ }$\cap $\textit{\ }$C^{0}(\bar{\Omega})$\textit{\ satisfy }%
\begin{eqnarray*}
H_{\vec{F}}(v) &\leq &H_{\vec{F}}(u)\text{ in }\Omega \backslash \{S_{\vec{F}%
}(u)\cup S_{\vec{F}}(v)\} \\
v &\geq &u\text{ on }\partial \Omega
\end{eqnarray*}%
\textit{\noindent Suppose }$\mathcal{H}_{m-1}(\overline{S_{\vec{F}}(u)\cup
S_{\vec{F}}(v)})$\textit{\ }$=$\textit{\ }$0$\textit{\ and }$\vec{F}$ $\in $ 
$C^{1}(\Omega )$ $\cap $ $C^{0}(\bar{\Omega})$\textit{\ satisfies }$\func{div%
}\vec{F}^{b}$\textit{\ }$>$\textit{\ }$0$\textit{\ (}$\func{div}\vec{F}^{b}$%
\textit{\ }$<$\textit{\ }$0,$\textit{\ resp.) in }$\Omega $\textit{. Then }$%
v\geq u$\textit{\ in }$\Omega .$

\bigskip

We remark that the condition $\func{div}\vec{F}^{b}$\textit{\ }$>$\textit{\ }%
$0$\textit{\ }or\textit{\ }$<$\textit{\ }$0$ was first used to extend
uniqueness results from even dimension to arbitrary dimension in \cite{ch2}.
In view of Corollary G (Corollary I, resp.), the condition on $N_{\vec{F}%
}(v) $ or $N_{\vec{F}}(u)$ in Theorem $F^{\prime }$ and Theorem J can be
replaced by (\ref{1.2}) ($\vec{F}$ $=$ $(-x^{2},$ $x^{1},$ ..., $-x^{2n},$ $%
x^{2n-1}), $ resp. for $m$ $=$ $2n$). Let $H_{n}$ denote the Heisenberg
group of dimension $2n$ $+$ $1.$ As a set, $H_{n}$ is $C^{n}\times R$ or $%
R^{2n}\times R.$ For a hypersurface $\Sigma $ in $H_{n}$ (which may not be a
graph over $R^{2n}$)$,$ the horizontal (or $p$-) mean curvature $H_{\Sigma }$
of $\Sigma $ for a defining function $\psi $ is given by%
\begin{equation}
H_{\Sigma }:=\func{div}_{b}\frac{\nabla _{b}\psi }{|\nabla _{b}\psi |}
\label{1.3}
\end{equation}%
\noindent where $\nabla _{b}$ and $\func{div}_{b}$ denote subgradient and
subdivergence in $H_{n}$, resp.. See Subsection B of the Appendix for
equivalent definitions of mean curvature in subriemannian geometry. Note
that for a graph $\Sigma $ over $R^{2n}$ defined by $u,$ $H_{\Sigma }$ may
be different from $H_{\vec{F}}(u)$ ($\vec{F}$ $=$ $(-x^{2},$ $x^{1},$ ..., $%
-x^{2n},$ $x^{2n-1})$ by sign. In fact, if we replace $\psi $ by $-\psi $ in
(\ref{1.3})$,$ $H_{\Sigma }$ becomes $-H_{\Sigma }.$

For the boundary $\Sigma $ of a ($C^{2}$ smooth, say) bounded domain $\Omega 
$ in $H_{n},$ we choose a defining function $\psi $ for $\Sigma ,$ such that 
$\psi $ $<$ $0$ in $\Omega .$ In this way $H_{\Sigma }$ is a positive
constant for a Pansu sphere given by the union of all the geodesics of
positive constant curvature joining the two poles (see, e.g., \cite{RR} for
the $m$ $=$ $2$ case and \cite{R} for the higher dimensional case).

\bigskip

\textbf{Theorem K.}\textit{\ Let }$\Sigma _{1}$\textit{\ and }$\Sigma _{2}$%
\textit{\ be two }$C^{2}$\textit{\ smooth, connected, orientable, closed
hypersurfaces of constant horizontal (}$p$\textit{-) mean curvature }$%
H_{\Sigma _{1}}$\textit{\ and }$H_{\Sigma _{2}},$\textit{\ resp. in }$H_{n},$%
\textit{\ }$n$\textit{\ }$\geq $\textit{\ }$2.$\textit{\ Suppose }$\Sigma
_{2}$\textit{\ is inscribed in }$\Sigma _{1},$\textit{\ i.e., }$\Sigma _{2}$%
\textit{\ is contained in the closure of the inside of }$\Sigma _{1}$\textit{%
\ and }$\Sigma _{1}$\textit{\ }$\cap $\textit{\ }$\Sigma _{2}$\textit{\ is
not empty. Assume }$H_{\Sigma _{2}}$\textit{\ }$\leq $\textit{\ }$H_{\Sigma
_{1}}$\textit{\ and }$\Sigma _{1}$\textit{\ }$\cap $\textit{\ }$\Sigma _{2}$%
\textit{\ contains a nonsingular (with respect to both }$\Sigma _{1}$\textit{%
\ and }$\Sigma _{2})$\textit{\ point or an isolated singular point of }$%
\Sigma _{1}$ \textit{(}$\Sigma _{2},$\textit{\ resp.)}$.$\textit{\ Moreover,
assume either }$\Sigma _{1}$\textit{\ or }$\Sigma _{2}$\textit{\ has only
isolated singular points.} \textit{Then }$\Sigma _{1}$\textit{\ }$\equiv $%
\textit{\ }$\Sigma _{2}.$

\textit{\bigskip }

For further applications we need to extend Theorem J to hypersurfaces of a
subriemannian manifold having isometric translations, with an isolated
singular point in touch. A point $\tilde{q}$ $\in $ $\tilde{\Sigma},$ a
hypersurface of a subriemannian manifold, is called singular if $\xi $ $%
\subset $ $T\tilde{\Sigma}$ at $\tilde{q}$ (this is equivalent to $|d\phi
|_{\ast }$ $=$ $0$ at $\tilde{q}$ for a defining function $\phi $ of $\tilde{%
\Sigma}$ mentioned previously$).$ For $\tilde{\Sigma}$ being a graph
described by $(x^{1},$\textit{\ }$x^{2},$\textit{\ ...}$x^{m},$\textit{\ }$%
w(x^{1},$\textit{\ }$x^{2},$\textit{\ ...}$x^{m})),$ ($x^{1},$\textit{\ }$%
x^{2},$\textit{\ ...}$x^{m})$ $\in $ $D,$ in local coordinates, we call a
point $q$ in $D$ singular for $w$ if $(q,w(q))$ is a singular point of $%
\tilde{\Sigma}$. Denote the set of all singular points in $D$ for $w$ by $%
S_{D}(w)$ or $S(w)$ if the domain of $w$ is clear in the context$.$

\bigskip

\textbf{Theorem J}$^{\prime }.$ \textit{Suppose} $(M,$\textit{\ }$<\cdot
,\cdot >^{\ast },dv_{M})$\textit{\ of dimension }$m+1$ \textit{has isometric
translations }$\Psi _{a}$ \textit{near }$p_{0}$\textit{\ }$\in $\textit{\ }$%
M,$ \textit{transversal to a hypersurface }$\Sigma $\textit{\ passing
through }$p_{0}.$ \textit{Take a system of translation-isometric coordinates}
$x^{1},$\textit{\ }$x^{2},$\textit{\ ..., }$x^{m+1}$\textit{\ in a
neighborhood }$\Omega $\textit{\ of }$p_{0}$ \textit{such that }$x^{m+1}$%
\textit{\ }$=$\textit{\ }$0$\textit{\ on (connected) }$\Sigma \cap \Omega .$%
\textit{\ Suppose }$u$\textit{\ }$(v,$\textit{\ resp.) }$:$\textit{\ }$%
\Sigma \cap \Omega \rightarrow R$\textit{\ defines a graph }$\{(x^{1},$%
\textit{\ }$x^{2},$\textit{\ ...}$x^{m},$\textit{\ }$u(x^{1},$\textit{\ }$%
x^{2},$\textit{\ ...}$x^{m}))\}$\textit{\ (}$\{(x^{1},$\textit{\ }$x^{2},$%
\textit{\ ...}$x^{m},$\textit{\ }$v(x^{1},$\textit{\ }$x^{2},$\textit{\ ...}$%
x^{m}))\}$\textit{, resp.) }$\subset $\textit{\ }$\Omega $\textit{. Assume }

\textit{(1) }$v\geq u$\textit{\ on }$\Sigma \cap \Omega $\textit{\ such that 
}$S_{\Sigma \cap \Omega }(u)$\textit{\ }$=$\textit{\ }$\{p_{0}\}$ \textit{(}$%
S_{\Sigma \cap \Omega }(v)$\textit{\ }$=$\textit{\ }$\{p_{0}\},$\textit{\
resp.)}$,$ $p_{0}$ $=$ $(0,..,0),$ \textit{and }$v(p_{0})$\textit{\ }$=$%
\textit{\ }$u(p_{0})$ $=$ $0;$

\textit{(2) }$H(v)\leq H(u)$\textit{\ in (}$\Sigma \cap \Omega )\backslash
\{\{p_{0}\}\cup S_{\Sigma \cap \Omega }(v)\}\mathit{\ ((\Sigma \cap }\Omega
)\backslash \{\{p_{0}\}\cup S_{\Sigma \cap \Omega }(u)\},\mathit{\ resp.)}.$

\textit{Suppose }$\mathcal{H}_{m-1}(\overline{S_{\Sigma \cap \Omega }(v)})$%
\textit{\ }$=$\textit{\ }$0$\textit{\ (}$\mathcal{H}_{m-1}(\overline{%
S_{\Sigma \cap \Omega }(u)})$\textit{\ }$=$\textit{\ }$0,$\textit{\ resp.)}$%
. $ \textit{Moreover, we assume the rank condition (\ref{I2-2}) holds}$.$ 
\textit{Then we have }$v\equiv u$ \textit{in} $\Sigma \cap \Omega .$

\bigskip

\textbf{Theorem \~{C}}$^{\prime \prime }.$ \textit{Suppose} $(M,$\textit{\ }$%
<\cdot ,\cdot >^{\ast },dv_{M})$\textit{\ of dimension }$m+1$\textit{\ has
isometric translations }$\Psi _{a}$ \textit{near }$p_{0}$\textit{\ }$\in $%
\textit{\ }$M,$ \textit{transversal to a hypersurface }$\Sigma $\textit{\
passing through }$p_{0}.$ \textit{Take a system of translation-isometric
coordinates} $x^{1},$\textit{\ }$x^{2},$\textit{\ ..., }$x^{m+1}$\textit{\
in an open neighborhood }$\Omega $\textit{\ of }$p_{0}$ \textit{such that }$%
x^{m+1}$\textit{\ }$=$\textit{\ }$0$\textit{\ on }$\Sigma \cap \Omega .$%
\textit{\ Let }$V$ $\subset $ $\bar{V}$ $\subset $ $\Omega $ \textit{be a
smaller open neighborhood of }$p_{0}.$ \textit{Let }$v,$\textit{\ }$u$%
\textit{\ }$\in $\textit{\ }$C^{2}(\Sigma \cap V)$\textit{\ }$\cap $\textit{%
\ }$C^{0}(\overline{\Sigma \cap V})$\textit{\ define graphs in }$\Omega $ 
\textit{and} \textit{satisfy}%
\begin{equation*}
H(v)\leq H(u)\text{ in (}\Sigma \cap V)\backslash \{S_{\Sigma \cap V}(u)\cup
S_{\Sigma \cap V}(v)\},
\end{equation*}%
\begin{equation*}
v\geq u\text{ on }\partial (\Sigma \cap V).
\end{equation*}%
\textit{\ \noindent Assume the rank condition (\ref{I2-2}) and }$\mathcal{H}%
_{m-1}(\overline{S_{\Sigma \cap V}(u)\cup S_{\Sigma \cap V}(v)})$\textit{\ }$%
=$\textit{\ }$0$\textit{.} \textit{Then }$v$ $\geq $ $u$ \textit{in }$\Sigma
\cap V.$

\bigskip

We can now apply Theorem C (or C$^{\prime })$ and Theorem J$^{\prime }$ to
prove a rigidity result for minimal hypersurfaces in a Heisenberg cylinder $%
(H_{n}\backslash \{0\},$\textit{\ }$\rho ^{-2}\Theta )$ with $n$\textit{\ }$%
\geq $\textit{\ }$2.$ Here $\Theta $ denotes the standard Heisenberg contact
form:%
\begin{equation*}
\Theta :=dz+\sum_{j=1}^{n}(x^{j}dx^{n+j}-x^{n+j}dx^{j}).
\end{equation*}

\noindent The associated Heisenberg distance function $\rho $ reads%
\begin{equation*}
\rho :=[(\sum_{K=1}^{2n}(x^{K})^{2})^{2}+4z^{2}]^{1/4}.
\end{equation*}

\noindent In Section 6 we discuss some basic geometry associated to the
contact form $\rho ^{-2}\Theta $ (using $x_{j},$ $y_{j}$ instead of $x^{j},$ 
$x^{n+j}$ and both interchangeably) before proving the following result.

\bigskip

\textbf{Theorem L.} \textit{Let }$\Sigma $\textit{\ be a closed (compact
with no boundary) immersed hypersurface in a Heisenberg cylinder }$%
(H_{n}\backslash \{0\},$\textit{\ }$\rho ^{-2}\Theta )$ \textit{with }$n$%
\textit{\ }$\geq $\textit{\ }$2.$\textit{\ Suppose either }

\textit{(a) }$H_{\Sigma }\leq 0$ \textit{or }

\textit{(b) }$H_{\Sigma }\geq 0$\textit{\ and the interior region of }$%
\Sigma $\textit{\ contains the origin of }$H_{n}$\textit{\ }

\noindent \textit{holds.} \textit{Then }$\Sigma $\textit{\ must be a
Heisenberg sphere defined by }$\rho ^{4}$\textit{\ }$=$\textit{\ }$c$\textit{%
\ for some constant }$c>0.$\textit{\ In particular, }$H_{\Sigma }$\textit{\ }%
$\equiv $\textit{\ }$0.$

\textit{\bigskip }

\textbf{Corollary M}. \textit{There does not exist a closed} \textit{%
immersed hypersurface of positive constant horizontal (}$p$\textit{-)mean
curvature} \textit{in a Heisenberg cylinder }$(H_{n}\backslash \{0\},$%
\textit{\ }$\rho ^{-2}\Theta )$\textit{\ with }$n\mathit{\ }\geq \mathit{\ }%
2,$ \textit{whose interior region contains the origin.}

\bigskip

Let $\varphi $ be a continuous function of $\tau $ $\in $ $[0,\infty ).$ We
have the following nonexistence result (pseudo-halfspace theorem).

\bigskip\ 

\textbf{Theorem N}. \textit{Let }$\Omega $\textit{\ be a domain of }$H_{n},$%
\textit{\ }$n$\textit{\ }$\geq $\textit{\ }$2,$\textit{\ defined by either }$%
z$\textit{\ }$>$\textit{\ }$\varphi (\sqrt{x_{1}^{2}+..+x_{2n}^{2}})$\textit{%
\ or }$x_{1}$\textit{\ }$>$\textit{\ }$\varphi (\sqrt{%
x_{2}^{2}+..+x_{2n}^{2}+z^{2}})$\textit{\ where }$\lim_{\tau \rightarrow
\infty }\varphi (\tau )$\textit{\ }$=$\textit{\ }$\infty $\textit{. Then
there does not exist any horizontal (}$p$\textit{-) minimal hypersurface
properly immersed in }$\Omega .$

\bigskip

The simplest example for Theorem N is $\varphi (\tau )$ $=$ $a\tau $ with $a$
$>$ $0.$ Call associated domains wedge-shaped. Theorem N tells us
nonexistence of horizontal ($p$-) minimal hypersurfaces in wedge-shaped
domains. But Theorem N does not hold for the case $a$ $=$ $0.$ That is,
halfspace theorem does not hold since there are catenoid type horizontal ($p$%
-)minimal hypersurfaces with finite height (\cite{RR0}) in $H_{n}$ for $n$ $%
\geq $ $2.$ On the other hand, we do have halfspace theorem for $H_{1}$ (see 
\cite{ch}). Hoffman and Meeks (\cite{HM}) first\ proved such a halfspace
theorem for $R^{3}.$ It fails for $R^{n}$ with $n\geq 4.$ But above type of
pseudo-halfspace theorem still holds for $R^{n}$ with $n\geq 4$ by a similar
reasoning.

There is another notion of mean curvature, called Levi-mean curvature, in
the study of real hypersurfaces in $C^{n}.$ We would like to remark that the
SMP for such mean curvature operators (generalized to pseudoconvex fully
nonlinear Levi-type curvature operators) has been proved by Montanari and
Lanconelli (\cite{ML}).

\bigskip

\textbf{Acknowledgements}. J.-H. Cheng, H.-L. Chiu, and J.-F. Hwang would
like to thank the Ministry of Science and Technology of Taiwan for the
support of the projects: MOST 104-2115 -M-001 -011 -MY2, 104-2115 -M-008
-003 -MY2, and 104-2115 -M-001 -009 -MY2, resp.. J.-H. Cheng is also
grateful to the National Center for Theoretical Sciences of Taiwan for the
constant support. P. Yang would like to thank the NSF of the U.S. for the
support of the project: DMS-1509505.

\bigskip

\section{Proofs of Theorems A, A$^{\prime }$, \~{A}, and B}

\proof
(\textbf{of Theorem A}) From the definition of $Q$ (see \ref{I1}), we
compute the difference of $Q\phi $ and $Q\psi $ as follows:%
\begin{eqnarray}
&&Q\phi -Q\psi  \label{B1} \\
&=&a^{ij}(x,D\phi )D_{ij}(\phi -\psi )+(a^{ij}(x,D\phi )-a^{ij}(x,D\psi
))D_{ij}\psi  \notag \\
&&+(b(x,\phi ,D\phi )-b(x,\phi ,D\psi ))+(b(x,\phi ,D\psi )-b(x,\psi ,D\psi
))\geq 0  \notag
\end{eqnarray}

\noindent by assumption. Writing%
\begin{eqnarray*}
w &=&\phi -\psi , \\
\tilde{a}^{ij}(x) &=&a^{ij}(x,D\phi ), \\
(a^{ij}(x,D\phi )-a^{ij}(x,D\psi ))D_{ij}\psi +(b(x,\phi ,D\phi )-b(x,\phi
,D\psi )) &=&\tilde{b}^{i}(x)D_{i}w, \\
b(x,\phi ,D\psi )-b(x,\psi ,D\psi ) &=&\tilde{a}(x)w,
\end{eqnarray*}

\noindent we get%
\begin{equation*}
Lw:=\tilde{a}^{ij}(x)D_{ij}w+\tilde{b}^{i}(x)D_{i}w+\tilde{a}(x)w\geq 0.
\end{equation*}

\noindent Noting that the quadratic form $(\tilde{a}^{ij}(x))$ is
nonnegative by condition (1) and $\tilde{a}(x)\leq 0$ by condition (2), we
can then apply Theorem A7 in the Appendix (Theorem 3.1 in \cite{Bony}) to
complete the proof.

\endproof%

\bigskip

\proof
(\textbf{of Theorem A}$^{\prime }$) Write $e_{i}=\alpha _{i}^{l}\partial
_{l}.$ Compute%
\begin{equation*}
a^{ij}e_{i}e_{j}\phi =a^{ij}\alpha _{i}^{l}\alpha _{j}^{k}\partial
_{l}\partial _{k}\phi +a^{ij}\alpha _{i}^{l}(\partial _{l}\alpha
_{j}^{k})\partial _{k}\phi .
\end{equation*}

\noindent Observe that $\tilde{a}^{lk}$ $:=$ $a^{ij}\alpha _{i}^{l}\alpha
_{j}^{k}$ satisfies $\tilde{a}^{lk}=\tilde{a}^{kl}$ since $a^{ij}=a^{ji},$
and $(\tilde{a}^{lk})$ is nonnegative since $(a^{ij})$ is nonnegative. Note
also that coefficients of first derivatives $\partial _{k}\phi $ do not rely
on the variable $\phi .$ So $Q^{\prime }\phi $ is of the form (\ref{I1}) for
a certain $Q\phi $ which satisfies the conditions (1), (2), (3) in Theorem
A. Thus the conclusion follows from Theorem A.

\endproof%

\bigskip

\proof
(\textbf{of Theorem \~{A}}) It suffices to show (\ref{I2-0}) implies
condition (3).

\bigskip

\endproof%

\bigskip

Given a subriemannian manifold $(M,<\cdot ,\cdot >^{\ast })$ of dimension $%
m+1,$ we recall that $G$ $:T^{\ast }M\rightarrow TM$ is defined by $\omega
(G(\eta ))$\textit{\ }$=$\textit{\ }$<\omega ,\eta >^{\ast }$\textit{\ }for%
\textit{\ }$\omega ,$\textit{\ }$\eta $\textit{\ }$\in $\textit{\ }$T^{\ast
}M.$ Define%
\begin{equation*}
<v,w>_{\ast }:=<\eta ,\zeta >^{\ast }
\end{equation*}%
\textit{\noindent }for $v,w$ $\in $ $\xi $ $:=$ $Range(G)$ and any choice $%
\eta $ $\in $ $G^{-1}(v),$ $\zeta $ $\in $ $G^{-1}(w).$ It is easy to see
that $<\cdot ,\cdot >_{\ast }$ is well defined. Assume 
\begin{equation*}
\dim \xi =\text{constant }m+1-l
\end{equation*}%
\textit{\noindent }for an integer $l,$ $0$ $\leq $ $l$ $\leq $ $m+1.$ So $%
<\cdot ,\cdot >_{\ast }$ is positive definite (for $l$ $\leq $ $m)$ on the
vector bundle $\xi $. There exist ($C^{\infty })$ smooth local sections $%
v_{1},$ $..,$ $v_{m+1-l}$ of $\xi ,$ orthonormal with respect to $<\cdot
,\cdot >_{\ast }($for $).$ We choose any smooth element $\eta ^{j}$ $\in $ $%
G^{-1}(v_{j}).$ It follows that 
\begin{equation*}
<\eta ^{i},\eta ^{j}>^{\ast }=\delta _{ij}
\end{equation*}%
\textit{\noindent }for $1$ $\leq $ $i,j$ $\leq $ $m+1-l.$ Now given a
background volume form $dv_{M},$ for $l$ $\geq $ $1,$ we then choose smooth
independent sections $\eta ^{m+2-l},$ $...,$ $\eta ^{m+1}$ of $KerG$ $%
\subset $ $T^{\ast }M$ such that $<\eta ^{j},\cdot >^{\ast }$ $=$ $<\cdot
,\eta ^{j}>^{\ast }$ $=$ $0$ for $m+2-l$ $\leq $ $j$ $\leq $ $m+1$ and 
\begin{equation*}
dv_{M}=\eta ^{1}\wedge \mathit{\ }...\wedge \mathit{\ }\eta ^{m+1}
\end{equation*}

\noindent (note that we have freedom to choose a scalar multiple of $\eta
^{m+2-l},$ $...,$ $\eta ^{m+1}).$ For $l$ $=$ $0,$ $<\cdot ,\cdot >^{\ast }$
is a Riemannian metric on $T^{\ast }M$ and $\eta ^{1}\wedge \mathit{\ }%
...\wedge \mathit{\ }\eta ^{m+1}$ is the Riemannian volume form (up to a
sign). So a given volume form $dv_{M}$ is a nonzero scalar multiple of $\eta
^{1}\wedge \mathit{\ }...\wedge \mathit{\ }\eta ^{m+1}.$ We have shown

\bigskip

\textbf{Lemma 2.1}. \textit{Let }$(M,<\cdot ,\cdot >^{\ast },dv_{M})$\textit{%
\ be a subriemannian manifold of dimension }$m+1$ \textit{with a background
volume form }$dv_{M}.$ \textit{Assume}%
\begin{equation*}
\dim \xi =\text{constant }m+1-l.
\end{equation*}%
\textit{\noindent (cf. (\ref{I2-1})) Then locally we can choose a suitable (}%
$C^{\infty }$\textit{\ smooth) coframe }$\eta ^{1},$\textit{\ }$..,$\textit{%
\ }$\eta ^{m+1}$\textit{\ such that }$<\eta ^{i},\eta ^{j}>^{\ast }$\textit{%
\ }$=$\textit{\ }$\delta _{ij}$\textit{\ for }$1$\textit{\ }$\leq $\textit{\ 
}$i,j$\textit{\ }$\leq $\textit{\ }$m+1-l$\textit{, }$<\eta ^{i},\eta
^{j}>^{\ast }$\textit{\ }$=$\textit{\ }$0$\textit{\ otherwise, and }$dv_{M}$%
\textit{\ = }$\eta ^{1}\wedge $\textit{\ }$...\wedge $\textit{\ }$\eta
^{m+1} $ \textit{for }$1$\textit{\ }$\leq $\textit{\ }$l$\textit{\ }$\leq $%
\textit{\ }$m+1$\textit{\ while }$dv_{M}$\textit{\ is a nonzero scalar
multiple of Riemannian volume form for }$l$\textit{\ }$=$\textit{\ }$0.$

\textit{\bigskip }

\proof
\textbf{(of Theorem B) }Take a coframe $\omega ^{1},$ $\omega ^{2},$ $...,$ $%
\omega ^{m+1}$ in $T^{\ast }M$ (near $p_{0}),$ such that $dv_{M}$ $=$ $%
\omega ^{1}\wedge \omega ^{2}$ $\wedge $ $...$ $\wedge $ $\omega ^{m+1}.$
Let $e_{1},$ $e_{2},$ $...,$ $e_{m+1}$ be the dual frame in $TM.$ Compute%
\begin{eqnarray}
d\phi \lrcorner \text{ }dv_{M} &=&(e_{i}\phi )\omega ^{i}\lrcorner \text{ }%
\omega ^{1}\wedge \omega ^{2}\wedge ...\wedge \omega ^{m+1}  \label{3.3} \\
&=&(e_{i}\phi )g^{ij}(-1)^{j-1}\omega ^{1}\wedge ..\wedge \hat{\omega}%
^{j}\wedge ..\wedge \omega ^{m+1}  \notag
\end{eqnarray}

\noindent where $g^{ij}$ $=$ $<\omega ^{i},\omega ^{j}>^{\ast }$ and $\hat{%
\omega}^{j}$ means deleting $\omega ^{j}.$ From (\ref{3.3}) and $|d\phi
|_{\ast }$ $\neq $ $0$ at $p_{0},$ we then compute (near $p_{0})$%
\begin{eqnarray}
&&d(\frac{d\phi }{|d\phi |_{\ast }^{1-p}}\lrcorner \text{ }dv_{M})
\label{3.4} \\
&=&d(\frac{e_{i}\phi }{|d\phi |_{\ast }^{1-p}})g^{ij}(-1)^{j-1}\wedge \omega
^{1}\wedge ..\wedge \hat{\omega}^{j}\wedge ..\wedge \omega ^{m+1}  \notag \\
&&+\frac{e_{i}\phi }{|d\phi |_{\ast }^{1-p}}(-1)^{j-1}d(g^{ij}\omega
^{1}\wedge ..\wedge \hat{\omega}^{j}\wedge ..\wedge \omega ^{m+1}).  \notag
\end{eqnarray}

\noindent The first term of the right-hand side in (\ref{3.4}) is the term
of second order in $\phi .$ We compute it as follows:%
\begin{eqnarray}
&&d(\frac{e_{i}\phi }{|d\phi |_{\ast }^{1-p}})g^{ij}(-1)^{j-1}\wedge \omega
^{1}\wedge ..\wedge \hat{\omega}^{j}\wedge ..\wedge \omega ^{m+1}
\label{3.5} \\
&=&e_{k}(\frac{e_{i}\phi }{|d\phi |_{\ast }^{1-p}})\omega
^{k}g^{ij}(-1)^{j-1}\wedge \omega ^{1}\wedge ..\wedge \hat{\omega}^{j}\wedge
..\wedge \omega ^{m+1}  \notag \\
&=&e_{j}(\frac{e_{i}\phi }{|d\phi |_{\ast }^{1-p}})g^{ij}dv_{M}.  \notag
\end{eqnarray}

\noindent So by (\ref{I2-0-1}), (\ref{3.4}), and (\ref{3.5}), we have%
\begin{equation}
H_{\phi ,p}=e_{j}(\frac{e_{i}\phi }{|d\phi |_{\ast }^{1-p}})g^{ij}+\text{%
first order terms in }\phi .  \label{3.6}
\end{equation}

\noindent Note that the first order terms in (\ref{3.6}) do not depend on
the variable $\phi $ itself.

Write 
\begin{eqnarray}
e_{j}(\frac{e_{i}\phi }{|d\phi |_{\ast }^{1-p}})g^{ij} &=&g^{ij}\frac{%
e_{i}e_{j}\phi }{|d\phi |_{\ast }^{1-p}}  \label{3.7} \\
&&+(p-1)\frac{g^{ij}(e_{i}\phi )(e_{k}\phi )(e_{j}e_{l}\phi )g^{kl}}{|d\phi
|_{\ast }^{3-p}}  \notag \\
&&+(p-1)\frac{g^{ij}(e_{i}\phi )(e_{k}\phi )(e_{l}\phi )(e_{j}g^{kl})}{%
2|d\phi |_{\ast }^{3-p}}.  \notag
\end{eqnarray}

\noindent Define the second order operator $a^{ij}(x,e_{1}\phi
,..,e_{m+1}\phi )e_{i}e_{j}$ by%
\begin{equation}
a^{ij}(x,e_{1}\phi ,..,e_{m+1}\phi ):=\frac{g^{ij}(x)}{|d\phi |_{\ast }^{1-p}%
}+(p-1)\frac{(e^{i}\phi )(e^{j}\phi )}{|d\phi |_{\ast }^{3-p}}  \label{3.12}
\end{equation}

\noindent where $e^{i}:=g^{ij}e_{j}.$ Observe that $H_{\phi ,p}$ is
independent of the choice of (co)frames. By Lemma 2.1, we can choose a
suitable coframe (field), denoted as $\tilde{\omega}^{1},$ $\tilde{\omega}%
^{2},$ $...,$ $\tilde{\omega}^{m+1},$ such that $\tilde{g}^{ij}$ $=$ $\delta
_{ij}$ for $1$ $\leq $ $i,j$ $\leq $ $m+1-l,$ $\tilde{g}^{ij}$ $=$ $0$
otherwise, and $dv_{M}$ $=$ $\tilde{\omega}^{1}\wedge \tilde{\omega}^{2}$ $%
\wedge $ $...$ $\wedge $ $\tilde{\omega}^{m+1}$ for $l$ $\geq $ $1$. If $l$ $%
=$ $0,$ $dv_{M}$ is a nonzero scalar multiple of the Riemannian volume form $%
\tilde{\omega}^{1}\wedge \tilde{\omega}^{2}$ $\wedge $ $...$ $\wedge $ $%
\tilde{\omega}^{m+1}$. So we have the same form of second order term and the
later argument still works. Let \{$\tilde{e}_{1},$ $..,$ $\tilde{e}_{m+1}\}$
be dual to $\{\tilde{\omega}^{1},$ $\tilde{\omega}^{2},$ $...,$ $\tilde{%
\omega}^{m+1}\}.$ Let%
\begin{equation}
\check{e}_{m+1-l}:=\frac{1}{|d\phi (x)|_{\ast }}\sum_{j=1}^{m+1-l}\tilde{e}%
_{j}\phi (x)\tilde{e}_{j}\in \xi  \label{3.13}
\end{equation}%
\noindent where $\xi $ := $Range(G)$ is spanned by the orthonormal basis $%
\tilde{e}_{1},$ $..,$ $\tilde{e}_{m+1-l}.$ Choose another system of
orthonormal vectors $\check{e}_{1},$ $..,$ $\check{e}_{m-l}$ perpendicular
to $\check{e}_{m+1-l}$ in $\xi .$ Also let $\check{e}_{m+2-l}$ $=$ $\tilde{e}%
_{m+2-l},$ $..,$ $\check{e}_{m+1}$ $=$ $\tilde{e}_{m+1}$. Consider 
\begin{eqnarray*}
\tilde{a}^{ij}(x) &:&=a^{ij}(x,\tilde{e}_{1}\phi (x),..,\tilde{e}_{m+1}\phi
(x)) \\
&=&\frac{\delta _{ij}}{|d\phi (x)|_{\ast }^{1-p}}+(p-1)\frac{\tilde{e}%
_{i}\phi (x)\tilde{e}_{j}\phi (x)}{|d\phi (x)|_{\ast }^{3-p}}
\end{eqnarray*}%
\noindent (viewed as a function of $x)$ for $1\leq i,j\leq m+1-l;$ $=$ $(p-1)%
\frac{(\tilde{e}_{i}\phi )(\tilde{e}_{j}\phi )}{|d\phi |_{\ast }^{3-p}}$
otherwise. Compute%
\begin{eqnarray}
&&\tilde{a}^{ij}(x)\tilde{e}_{i}\tilde{e}_{j}  \label{3.14} \\
&=&\frac{1}{|d\phi (x)|_{\ast }^{1-p}}\{\sum_{j=1}^{m+1-l}\tilde{e}%
_{j}^{2}+(p-1)\sum_{i,j=1}^{m+1}\frac{\tilde{e}_{i}\phi (x)\tilde{e}_{j}\phi
(x)}{|d\phi (x)|_{\ast }^{2}}\tilde{e}_{i}\tilde{e}_{j}\}  \notag \\
&=&\frac{1}{|d\phi (x)|_{\ast }^{1-p}}\{\sum_{j=1}^{m+1-l}\check{e}_{j}^{2}+%
\text{first order}+(p-1)\check{e}_{m+1-l}^{2}+\text{first order\}}  \notag \\
&=&\frac{1}{|d\phi (x)|_{\ast }^{1-p}}\{\sum_{j=1}^{m-l}\check{e}_{j}^{2}+p%
\check{e}_{m+1-l}^{2}\}+\text{first order}  \notag \\
&=&\sum_{j=1}^{m-l}X_{j}^{2}+(\frac{\sqrt{p}}{|d\phi (x)|_{\ast }^{(1-p)/2}}%
\check{e}_{m+1-l})^{2}+\text{first order}  \notag
\end{eqnarray}%
\noindent where 
\begin{equation*}
X_{j}:=\frac{\check{e}_{j}}{|d\phi (x)|_{\ast }^{(1-p)/2}}.
\end{equation*}%
\noindent Observe that $\check{e}_{j}\phi $ $=$ $0,$ and hence we have%
\begin{equation}
X_{j}\in \xi \cap T\{\phi =0\}.  \label{3.23}
\end{equation}

\noindent i.e., $X_{j}$ lies in the tangent space of hypersurface defined by 
$\phi $ $=$ $0.$ It is not hard to see that $(a^{ij})$ is symmetric and
nonnegative by Cauchy-Schwarz inequality. In view of (\ref{3.5}), (\ref{3.6}%
), and (\ref{3.7}), we learn that $H_{\phi ,p}$ $+b(x,\phi ,D\phi )$ is an
operator of type $Q^{\prime }\phi $ in (\ref{I2}). Observe that condition
(3) in Theorem A$^{\prime }$ holds by (\ref{3.14}). By assumption we have%
\begin{equation}
rank(\pounds (X_{1,}...,X_{m-l}))=m.  \label{3.24}
\end{equation}

\noindent Since $\dim \{\phi =0\}$ $=$ $m,$ the integral curves of all $Z$ $%
\in $ $\pounds (X_{1,}...,X_{m-l})$ will cover a neighborhood of $\{\phi
=0\} $ by (\ref{3.23}) and (\ref{3.24}). (a) follows from Theorem A$^{\prime
}.$

In case $p$ $>$ $0$, $X_{1},$ $X_{2},$ $..,$ $X_{m-l},$ $\frac{\sqrt{p}}{%
|d\phi (x)|_{\ast }^{(1-p)/2}}\check{e}_{m+1-l}$ form a basis of $\xi .$ By
the assumption $rank(\pounds (\xi ))$\textit{\ }$=$\textit{\ }$m+1,$ the
integral curves of all $Z$ $\in $ $\pounds (\xi )$ will cover a neighborhood
of $p_{0}$ in $M.$ (b) follows from (\ref{3.14}) and Theorem A$^{\prime }.$

\endproof%

\bigskip

\section{Graphs under symmetry and proofs of Theorems C, \^{C}, C$^{\prime
}, $ E and Corollary D}

Next we want to study when the conditions in Theorem B are satisfied. Let us
start with a general subriemannian manifold $(M,$ $<\cdot ,\cdot >^{\ast })$
where $<\cdot ,\cdot >^{\ast }$ is a nonnegative definite inner product on $%
T^{\ast }M.$ Take a background volume form $dv_{M}$. Let $\Sigma _{0}$ be a
(local) hypersurface in $M.$ Consider a one-parameter family of
diffeomorphisms $\Psi _{a}$ $:$ $M\rightarrow M$, i.e., $\Psi _{0}$ $=$ $%
Identity,$ $\Psi _{a+b}=\Psi _{a}\circ \Psi _{b}.$ We ask when the
hypersurface $\Sigma _{a}$ $:=$ $\Psi _{a}(\Sigma _{0})$ has the same mean
curvature (function) as $\Sigma _{0}.$

\bigskip

\textbf{Proposition 3.1}. \textit{Suppose} $(M,$\textit{\ }$<\cdot ,\cdot
>^{\ast },dv_{M})$ \textit{has isometric translations }$\Psi _{a}$ \textit{%
near }$p_{0}$\textit{\ }$\in $\textit{\ }$M,$ \textit{compatible with a
defining function }$\phi $ \textit{(see Section 1 for the definition)}$.$ 
\textit{Then for any }$p\geq 0,$ \textit{we have} $H_{\phi ,p}(\Psi _{a}(x))$%
\textit{\ }$=$\textit{\ }$H_{\phi ,p}(x).$

\bigskip

\proof
Recall that $H_{\phi ,p}$ is defined by%
\begin{equation*}
d(\frac{d\phi }{|d\phi |_{\ast }^{1-p}}\lrcorner \text{ }dv_{M}):=H_{\phi
,p}dv_{M}.
\end{equation*}

\noindent Pulling back the above identity by $\Psi _{a},$ we get%
\begin{equation}
d(\frac{d(\phi \circ \Psi _{a})}{|d(\phi \circ \Psi _{a})|_{\ast }^{1-p}}%
\lrcorner \text{ }dv_{M})=(H_{\phi ,p}\circ \Psi _{a})dv_{M}  \label{3.25}
\end{equation}

\noindent by (a) and (b) in Definition 1.1. By the compatibility of $\Psi
_{a}$ with $\phi $ (see (\ref{I3}))$,$ we compute the left hand side of (\ref%
{3.25}): 
\begin{eqnarray}
&&d(\frac{d(\phi \circ \Psi _{a})}{|d(\phi \circ \Psi _{a})|_{\ast }^{1-p}}%
\lrcorner \text{ }dv_{M})  \label{3.26} \\
&=&d(\frac{d\phi }{|d\phi |_{\ast }^{1-p}}\lrcorner \text{ }dv_{M})=H_{\phi
,p}dv_{M}.  \notag
\end{eqnarray}

\noindent Comparing (\ref{3.25}) with (\ref{3.26}), we obtain $H_{\phi
,p}\circ \Psi _{a}=H_{\phi ,p}.$

\endproof%

\bigskip

\proof
\textbf{(of Theorem C)} For any $x$ near $p_{0}$ there exist a number $b$
and $q$ $\in $ $\Sigma _{1}(\cap U)$ such that 
\begin{equation*}
x=\Psi _{b}(q).
\end{equation*}

\noindent Let $\tilde{q}$ $:=$ $\Psi _{a(q)}(q)$\textit{\ }$\in $\textit{\ }$%
\Sigma _{2}$ with $a(q)$ $\geq $ $0$ by condition (1). Choose defining
functions $\psi $ and $\phi $ compatible with $\Psi _{a}$ : first define $%
\psi (\Psi _{a}(x))$ $=$ $\psi (x)$ $-$ $a$ ($\phi (\Psi _{a}(x))$ $=$ $\phi
(x)$ $-$ $a,$ resp.) for $x$ $\in $ $\Sigma _{2}$ ($\Sigma _{1},$ resp.)$.$
The same formula then holds for any $x$ near $p_{0}$. We compute%
\begin{eqnarray}
\psi (x) &=&\psi (\Psi _{b}(q))=\psi (\Psi _{b-a(q)}(\tilde{q}))
\label{3.27} \\
&=&\psi (\tilde{q})-(b-a(q))  \notag \\
&=&0-b+a(q)\text{ \ (since }\tilde{q}\in \Sigma _{2})  \notag \\
&\geq &\phi (q)-b\text{ \ (since }q\in \Sigma _{1}\text{ and }a(q)\geq 0) 
\notag \\
&=&\phi (\Psi _{-b}(x))-b  \notag \\
&=&\phi (x)-(-b)-b=\phi (x).  \notag
\end{eqnarray}

On the other hand, we compute%
\begin{eqnarray}
H_{\psi ,p}(x) &=&H_{\psi ,p}(\Psi _{b}(q))=H_{\psi ,p}(\Psi _{b-a(q)}(%
\tilde{q}))  \label{3.28} \\
&=&H_{\psi ,p}(\tilde{q})\text{ \ (by Proposition 3.1)}  \notag \\
&=&H_{\Sigma _{2},p}(\tilde{q})\text{ \ (by Proposition B.1 in the Appendix)}
\notag \\
&=&H_{\Sigma _{2},p}(\Psi _{a(q)}(q))  \notag \\
&\leq &H_{\Sigma _{1},p}(q)\text{ \ (by condition (2))}  \notag \\
&=&H_{\phi ,p}(q)\text{ \ (by Proposition B.1})  \notag \\
&=&H_{\phi ,p}(\Psi _{b}(q))\text{ \ (by Proposition 3.1)}  \notag \\
&=&H_{\phi ,p}(x).  \notag
\end{eqnarray}

\noindent The result follows from (\ref{3.27}), (\ref{3.28}), and Theorem B
with $b$ $\equiv $ $0$.

\endproof%

\bigskip

\proof
\textbf{(of Corollary D) }Observe that the Heisenberg group $H_{n}$ can be
viewed as a subriemannian manifold\textbf{. }See Subsection B in the
Appendix. \{$l_{a}\},$ the Heisenberg translations, is a one-parameter
family of symmetries transversal to $\Sigma _{1}$ and $\Sigma _{2}$ near $%
p_{0}.$ Also we can choose $\phi $\ and $\psi $\ to be defining functions
for hypersurfaces $\Sigma _{1}$\ and $\Sigma _{2},$\ resp. (I.e., $\Sigma
_{1}$\ ($\Sigma _{2},$\ resp.) is defined by $\phi $\ $=$\ $0$\ ($\psi $\ $=$%
\ $0$, resp.)), compatible with $\{l_{a}\}$\textit{. }

We claim that $rank(\pounds (\xi \cap T\Sigma _{1}))$\textit{\ }$=$\textit{\ 
}$m$\textit{\ }near\textit{\ }$p_{0}.$ Observe that $\dim \xi \cap T\Sigma
_{1}$ $=$ $2n-1$ $\geq $ $3$ for $n$ $\geq $ $2.$ We can find a $J$%
-invariant pair of nonzero vectors $X,$ $JX$ in $\xi \cap T\Sigma _{1}.$ The
Lie bracket $[X,JX]$ generates the direction out of contact distribution $%
\xi .$ Therefore $rank(\pounds (\xi \cap T\Sigma _{1}))$\textit{\ }$=$%
\textit{\ }$2n$ $=$ $m$ near\textit{\ }$p_{0}.$ We then conclude the result
by Theorem C.

\endproof%

\bigskip

\proof
\textbf{(of Theorem \^{C})} Take any system of local coordinates $\hat{x}%
^{1},$\textit{\ }$\hat{x}^{2},$\textit{\ ..., }$\hat{x}^{m}$ on $\Sigma \cap
U$ where $U$ is a neighborhood of $p_{0}$, such that $\hat{x}^{j}(p_{0})$ $=$
$0.$ Observe that $V$ $:=$ $\cup _{a\in (-\delta ^{\prime },\delta ^{\prime
})}\Psi _{a}(\Sigma \cap U)$ for smaller $\delta ^{\prime }$ $<$ $\delta $
is a neighborhood of $p_{0}$ due to transversality of $\Psi _{a}$ to $\Sigma
.$Define the last coordinate $x^{m+1}$ and $x^{1},$\textit{\ }$x^{2},$%
\textit{\ ..., }$x^{m}$ in $V$ by 
\begin{eqnarray}
(1)\text{ }x^{m+1}(\Psi _{a}(\hat{q})) &=&a,  \label{3.31} \\
(2)\text{ }x^{j}(\Psi _{a}(\hat{q})) &=&\hat{x}^{j}(\hat{q}),\text{ }1\leq
j\leq m  \notag
\end{eqnarray}

\noindent for any $a\in (-\delta ^{\prime },\delta ^{\prime })$ and any $%
\hat{q}\in \Sigma \cap U$. Write an arbitrary point $q$ $\in $ $V$ as $q$ $=$
$\Psi _{b}(\hat{q})$ for $b$ $\in $ $(-\delta ^{\prime },\delta ^{\prime }).$
Choose $a$ such that $a+b$ $\in $ $(-\delta ^{\prime },\delta ^{\prime }).$
Compute%
\begin{eqnarray}
x^{j}(\Psi _{a}(q)) &=&x^{j}(\Psi _{a}(\Psi _{b}(\hat{q})))  \label{3.32} \\
&=&x^{j}(\Psi _{a+b}(\hat{q}))  \notag \\
&=&\hat{x}^{j}(\hat{q})\text{ \ (by (2) of (\ref{3.31}))}  \notag \\
&=&x^{j}(\Psi _{b}(\hat{q}))\text{ \ (by (2) of (\ref{3.31}))}  \notag \\
&=&x^{j}(q).  \notag
\end{eqnarray}

\noindent Similarly we have%
\begin{eqnarray}
x^{m+1}(\Psi _{a}(q)) &=&x^{m+1}(\Psi _{a}(\Psi _{b}(\hat{q})))  \label{3.33}
\\
&=&x^{m+1}(\Psi _{a+b}(\hat{q}))  \notag \\
&=&a+b\text{ \ (by (1) of (\ref{3.31}))}  \notag \\
&=&a+x^{m+1}(\Psi _{b}(\hat{q}))\text{ \ (by (1) of (\ref{3.31}))}  \notag \\
&=&a+x^{m+1}(q).  \notag
\end{eqnarray}

\noindent We have proved (1). Moreover, we conclude (2) by (b) of Definition
1.1 and (1).

\endproof%

\bigskip

\proof
(\textbf{of Theorem C}$^{\prime }$) Let $\phi $ $:=$ $u(x^{1},$\textit{\ }$%
x^{2},$\textit{\ ...}$x^{m})$ $-$ $x^{m+1}$ and $\psi $ $:=$ $v(x^{1},$%
\textit{\ }$x^{2},$\textit{\ ...}$x^{m})$ $-$ $x^{m+1}.$ Observe that $v\geq
u$ implies 
\begin{equation}
\psi \geq \phi .  \label{3.33-1}
\end{equation}%
\noindent For any $q$ near $p_{0},$ there are $a$ $=$ $x^{m+1}(q)$ $%
-v(x_{0}^{1},\mathit{\ ..,}x_{0}^{m}),$ $b$ $=$ $x^{m+1}(q)$ $-$ $%
u(x_{0}^{1},\mathit{\ ..,}x_{0}^{m})\ $($\in (-\delta ,\delta ))$ in which $%
x_{0}^{1}$ $=$ $x^{1}(q),$\textit{\ .., }$x_{0}^{m}$ $=$ $x^{m}(q),$ such
that 
\begin{eqnarray}
q &=&\Psi _{a}(x_{0}^{1},\mathit{\ ..,\ }x_{0}^{m},v(x_{0}^{1},\mathit{\
..,\ }x_{0}^{m}))  \label{3.34} \\
&=&\Psi _{b}(x_{0}^{1},\mathit{\ ..,\ }x_{0}^{m},u(x_{0}^{1},\mathit{\ ..,\ }%
x_{0}^{m})).  \notag
\end{eqnarray}%
\noindent by transversality of isometric translations$.$ Note that $\Psi
_{a} $ is compatible with $\phi $ and $\psi .$ From (\ref{3.34}) we compute%
\begin{eqnarray}
H_{\psi ,p}(q) &=&H_{\psi ,p}(\Psi _{a}(x_{0}^{1},\mathit{\ ..,\ }%
x_{0}^{m},v(x_{0}^{1},\mathit{\ ..,\ }x_{0}^{m}))  \label{3.35} \\
&=&H_{\psi ,p}(x_{0}^{1},\mathit{\ ..,\ }x_{0}^{m},v(x_{0}^{1},\mathit{\
..,\ }x_{0}^{m}))\text{ \ (by Proposition 3.1)}  \notag \\
&=&H_{p}(v)(x_{0}^{1},\mathit{\ ..,\ }x_{0}^{m})\text{ \ (by (\ref{I4}))} 
\notag \\
&\leq &H_{p}(u)(x_{0}^{1},\mathit{\ ..,\ }x_{0}^{m})\text{ \ (by assumption
(2))}  \notag \\
&=&H_{\phi ,p}(x_{0}^{1},\mathit{\ ..,\ }x_{0}^{m},u(x_{0}^{1},\mathit{\
..,\ }x_{0}^{m}))\text{ \ (by (\ref{I4}))}  \notag \\
&=&H_{\phi ,p}(\Psi _{b}(x_{0}^{1},\mathit{\ ..,\ }x_{0}^{m},u(x_{0}^{1},%
\mathit{\ ..,\ }x_{0}^{m}))\text{..(by Proposition 3.1)}  \notag \\
&=&H_{\phi ,p}(q)\text{ \ (by }(\ref{3.34})).  \notag
\end{eqnarray}

\noindent In view of (\ref{3.33-1}) and (\ref{3.35}), we have $\psi $ $=$ $%
\phi $ on $\Sigma _{1}$ $:=$ $\{\phi =0\}$ by Theorem B (a)$.$ It follows
that $v$ $\equiv $ $u$\textit{\ }in a neighborhood of $p_{0}$ $\in $ $\Sigma
.$

\endproof%

\bigskip

The argument in (\ref{3.35}) shows the following fact.

\bigskip

\textbf{Corollary 3.2. }\textit{Suppose} $(M,$\textit{\ }$<\cdot ,\cdot
>^{\ast },dv_{M})$\textit{\ of dimension }$m+1$\textit{\ with }$m$\textit{\ }%
$\geq $\textit{\ }$3$ \textit{has isometric translations }$\Psi _{a}$ 
\textit{near }$p$\textit{\ }$\in $\textit{\ }$M,$ \textit{transversal to a
hypersurface }$\Sigma $\textit{\ passing through }$p.$ \textit{Take a system
of translation-isometric coordinates} $x^{1},$\textit{\ }$x^{2},$\textit{\
..., }$x^{m+1}$\textit{\ in a neighborhood }$V$\textit{\ of }$p$ \textit{%
such that }$x^{m+1}$\textit{\ }$=$\textit{\ }$0$\textit{\ on }$\Sigma .$%
\textit{\ Suppose }$u$\textit{\ }$:$\textit{\ }$\Sigma \cap V\rightarrow R$%
\textit{\ defines a graph }$\{(x^{1},$\textit{\ }$x^{2},$\textit{\ ...}$%
x^{m},$\textit{\ }$u(x^{1},$\textit{\ }$x^{2},$\textit{\ ...}$x^{m}))\}$%
\textit{\ }$\subset $\textit{\ }$V$\textit{. Let }$\phi (x^{1},\ x^{2},%
\mathit{\ ..,.}x^{m},$ $x^{m+1}):=$ $u(x^{1},$\textit{\ }$x^{2},$\textit{\
...}$x^{m})$ $-$ $x^{m+1}.$ \textit{Then }$H_{\phi }(q)$ $=$ $H(u)(\pi
_{\Sigma }(q))$ \textit{wherever defined, in which }$\pi _{\Sigma }$\textit{%
\ is the projection to }$\Sigma $\textit{\ (i.e.}, ($x^{1},\ x^{2},\mathit{\
..,.}x^{m},$ $x^{m+1})$ $\rightarrow $ $(x^{1},$\textit{\ }$x^{2},$\textit{\
...}$x^{m})$ \textit{in translation-isometric coordinates}$)$

\bigskip

We now want to apply our theory of translation-isometric coordinates to the
situation of $l_{a}$ graphs. Observe that $\eta ^{2}=x^{2},$ $\eta
^{3}=x^{3},..,$ $\eta ^{2n}=x^{2n},$ and $\tau $ $:=$ $z+x^{1}x^{n+1}$ are
invariant under $l_{a},$ Heisenberg translations in the direction $\frac{%
\partial }{\partial x^{1}}.$ So an $l_{a}$ graph $u$ is a hypersurface $%
(u(\eta ^{2},\eta ^{3},..,\eta ^{2n},\tau ),$ $0,..,$ $0)\circ (0,$ $\eta
^{2},$ $..,$ $\eta ^{2n},$ $\tau )$ of $H_{n},$ parametrized by $\eta ^{2},$ 
$\eta ^{3},$ $..,$ $\eta ^{2n},$ $\tau $ so that 
\begin{eqnarray}
x^{1} &=&u(\eta ^{2},\eta ^{3},..,\eta ^{2n},\tau ),  \label{3.a} \\
x^{2} &=&\eta ^{2},...,x^{2n}=\eta ^{2n},  \notag \\
z &=&\tau -\eta ^{n+1}u(\eta ^{2},\eta ^{3},..,\eta ^{2n},\tau ).  \notag
\end{eqnarray}

\noindent In coordinates $\eta ^{2},\eta ^{3},..,\eta ^{2n},\tau ,$ and $%
x^{1},$ we write the standard contact form $\Theta $ for $H_{n}$ as follows:%
\begin{eqnarray}
\Theta
&:&=dz+x^{1}dx^{n+1}-x^{n+1}dx^{1}+%
\sum_{j=2}^{n}(x^{j}dx^{n+j}-x^{n+j}dx^{j})  \label{3.b} \\
&=&d\tau -2\eta ^{n+1}dx^{1}+\sum_{j=2}^{n}(\eta ^{j}d\eta ^{n+j}-\eta
^{n+j}d\eta ^{j}).  \notag
\end{eqnarray}

\noindent Observe that%
\begin{eqnarray*}
\mathring{e}_{2} &:&=\frac{\partial }{\partial \eta ^{2}}+\eta ^{n+2}\frac{%
\partial }{\partial \tau },..,\mathring{e}_{n}:=\frac{\partial }{\partial
\eta ^{n}}+\eta ^{2n}\frac{\partial }{\partial \tau }, \\
\mathring{e}_{1} &:&=\frac{\partial }{\partial x^{1}}+2\eta ^{n+1}\frac{%
\partial }{\partial \tau },\mathring{e}_{n+1}:=\frac{\partial }{\partial
\eta ^{n+1}}, \\
\mathring{e}_{n+2} &:&=\frac{\partial }{\partial \eta ^{n+2}}-\eta ^{2}\frac{%
\partial }{\partial \tau },..,\mathring{e}_{2n}:=\frac{\partial }{\partial
\eta ^{2n}}-\eta ^{n}\frac{\partial }{\partial \tau }.
\end{eqnarray*}

\noindent form an orthonormal basis of $\ker \Theta $ with respect to the
Levi metric $\frac{1}{2}d\Theta (\cdot ,J\cdot )$ (see Subsection B in the
Appendix for more explanation). We remark that the above $\mathring{e}_{j},$ 
$\mathring{e}_{n+j}$ are the same as the vectors $\hat{e}_{j},$ $\hat{e}%
_{j^{\prime }},$ resp. in Subsection B of the Appendix, but expressed in
different coordinates. Let $\phi $ $:=$ $u(\eta ^{2},\eta ^{3},..,\eta
^{2n},\tau )$ $-$ $x^{1}.$ We compute%
\begin{eqnarray}
\mathring{e}_{1}\phi &=&-1+2\eta ^{n+1}\partial _{\tau }u,\text{ }\mathring{e%
}_{n+1}\phi =\partial _{\eta ^{n+1}}u,  \label{3.c} \\
\mathring{e}_{2}\phi &=&\partial _{\eta ^{2}}u+\eta ^{n+2}\partial _{\tau
}u,..,\mathring{e}_{2n}\phi =\partial _{\eta ^{2n}}u-\eta ^{n}\partial
_{\tau }u.  \notag
\end{eqnarray}

\noindent Hence $d\phi $ $=$ ($\mathring{e}_{1}\phi )dx^{1}+\sum_{j=2}^{2n}(%
\mathring{e}_{j}\phi )d\eta ^{j}+(\partial _{\tau }\phi )\Theta $ has the
length with respect to the subriemannian metric (\ref{H2}) as follows:%
\begin{eqnarray}
|d\phi |_{\ast }^{2} &=&\sum_{j=1}^{2n}(\mathring{e}_{j}\phi
)^{2}=|W(u-x^{1})|^{2}  \label{3.d} \\
&=&(1-2\eta ^{n+1}\partial _{\tau }u)^{2}+(\partial _{\eta ^{n+1}}u)^{2} 
\notag \\
&&+\sum_{j=2}^{n}[(\partial _{\eta ^{j}}u+\eta ^{n+j}\partial _{\tau
}u)^{2}+(\partial _{\eta ^{n+j}}u-\eta ^{j}\partial _{\tau }u)^{2}]  \notag
\\
&=&1-4\eta ^{n+1}\partial _{\tau }u+|Wu|^{2}  \notag
\end{eqnarray}

\noindent where $W$ is the vector-valued operator $(\mathring{e}_{1},%
\mathring{e}_{2},..,\mathring{e}_{2n}).$ Note that the standard volume form $%
dV_{H_{n}}$ of $H_{n}$ equals $dx^{1}\wedge d\eta ^{2}\wedge ..\wedge d\eta
^{2n}\wedge \Theta .$

From (\ref{3.1}) we compute $H_{\phi }$ as follows:%
\begin{eqnarray*}
&&d(\frac{d\phi }{|d\phi |_{\ast }}\lrcorner \text{ }dV_{H_{n}}) \\
&=&d\{[\frac{\mathring{e}_{1}\phi }{|d\phi |_{\ast }}dx^{1}+\sum_{j=2}^{2n}%
\frac{\mathring{e}_{j}\phi }{|d\phi |_{\ast }}d\eta ^{j}+\frac{\partial
_{\tau }\phi }{|d\phi |_{\ast }}\Theta ]\lrcorner \text{ }dV_{H_{n}}\} \\
&=&[\mathring{e}_{1}(\frac{\mathring{e}_{1}\phi }{|d\phi |_{\ast }}%
)+\sum_{j=2}^{2n}\mathring{e}_{j}(\frac{\mathring{e}_{j}\phi }{|d\phi
|_{\ast }})]dV_{H_{n}}
\end{eqnarray*}

\noindent So we have%
\begin{eqnarray*}
H_{\phi } &=&\mathring{e}_{1}(\frac{\mathring{e}_{1}\phi }{|d\phi |_{\ast }}%
)+\sum_{j=2}^{2n}\mathring{e}_{j}(\frac{\mathring{e}_{j}\phi }{|d\phi
|_{\ast }})=W\cdot \left( \frac{W(u-x^{1})}{|W(u-x^{1})|}\right) \\
&=&2\eta ^{n+1}\partial _{\tau }\left( \frac{-1+2\eta ^{n+1}\partial _{\tau
}u}{\sqrt{1-4\eta ^{n+1}\partial _{\tau }u+|Wu|^{2}}}\right) \\
&&+\partial _{\eta ^{n+1}}\left( \frac{\partial _{\eta ^{n+1}}u}{\sqrt{%
1-4\eta ^{n+1}\partial _{\tau }u+|Wu|^{2}}}\right) \\
&&+\sum_{j=2}^{n}[(\partial _{\eta ^{j}}+\eta ^{n+j}\partial _{\tau })\left( 
\frac{\partial _{\eta ^{j}}u+\eta ^{n+j}\partial _{\tau }u}{\sqrt{1-4\eta
^{n+1}\partial _{\tau }u+|Wu|^{2}}}\right) \\
&&+(\partial _{\eta ^{n+j}}-\eta ^{j}\partial _{\tau })\left( \frac{\partial
_{\eta ^{n+j}}u-\eta ^{j}\partial _{\tau }u}{\sqrt{1-4\eta ^{n+1}\partial
_{\tau }u+|Wu|^{2}}}\right) ].
\end{eqnarray*}

\noindent So $H(u)$ has the above expression (see (\ref{I4}) with $p$ $=$ $0$
for the definition). Let $\Sigma _{1}$ denote the hypersurface defined by $%
u(\eta ^{2},\eta ^{3},..,\eta ^{2n},\tau )$ $=$ $x^{1}.$ As in the proof of
Corollary D, we observe that $\dim \xi \cap T\Sigma _{1}$ $=$ $2n-1$ $\geq $ 
$3$ for $n$ $\geq $ $2.$ We can then find a $J$-invariant pair of nonzero
vectors $X,$ $JX$ in $\xi \cap T\Sigma _{1}.$ The Lie bracket $[X,JX]$
generates a direction out of contact distribution $\xi .$ Therefore $rank(%
\pounds (\xi \cap T\Sigma _{1}))$\textit{\ }$=$\textit{\ }$2n$ $=$ $m$ near%
\textit{\ }$p_{0}.$ By Theorem C$^{\prime }$ we have

\bigskip

\textbf{Corollary D}$^{\prime }.$ \textit{Suppose }$n$ $\geq $ $2.$ \textit{%
Let }$x^{1}=u(\eta ^{2},\eta ^{3},..,\eta ^{2n},\tau )$\textit{\ and }$%
x^{1}=v(\eta ^{2},\eta ^{3},..,\eta ^{2n},\tau )$\textit{\ be two }$l_{a}$%
\textit{\ graphs defined on a common domain }$\Omega .$\textit{\ Assume }$%
v=u $\textit{\ at }$p_{0}$\textit{\ }$\in $\textit{\ }$\Omega ,$\textit{\ }$%
v\geq u$\textit{\ in }$\Omega $\textit{, and }$H(v)$\textit{\ }$\leq $%
\textit{\ }$H(u)$\textit{\ in }$\Omega .$\textit{\ Then }$v\equiv u$\textit{%
\ near }$p_{0}.$

\bigskip

We now want to apply our theory of translation-isometric coordinates to the
situation of intrinsic graphs. Let us recall (\cite{ASCV}) that an intrinsic
graph $u$ is a hypersurface of $H_{n},$ parametrized by $\eta ^{2},$ $\eta
^{3},$ $..,$ $\eta ^{2n},$ $\tau $ so that $(x^{1},$ $..,$ $x^{2n},$ $z)$ $=$
$(0,$ $\eta ^{2}$ $,\eta ^{3},..,\eta ^{2n},$ $\tau )\cdot (u(\eta ^{2},\eta
^{3},..,\eta ^{2n},\tau ),$ $0,..,$ $0)$ or%
\begin{eqnarray}
x^{1} &=&u(\eta ^{2},\eta ^{3},..,\eta ^{2n},\tau ),  \label{3.35-0} \\
x^{2} &=&\eta ^{2},...,x^{2n}=\eta ^{2n},  \notag \\
z &=&\tau +\eta ^{n+1}u(\eta ^{2},\eta ^{3},..,\eta ^{2n},\tau ).  \notag
\end{eqnarray}

\noindent In coordinates $\eta ^{2},\eta ^{3},..,\eta ^{2n},\tau ,$ and $%
x^{1},$ we write the standard contact form $\Theta $ for $H_{n}$ as follows:%
\begin{eqnarray}
\Theta
&:&=dz+x^{1}dx^{n+1}-x^{n+1}dx^{1}+%
\sum_{j=2}^{n}(x^{j}dx^{n+j}-x^{n+j}dx^{j})  \label{3.35-1} \\
&=&d\tau +2x^{1}d\eta ^{n+1}+\sum_{j=2}^{n}(\eta ^{j}d\eta ^{n+j}-\eta
^{n+j}d\eta ^{j}).  \notag
\end{eqnarray}

\noindent Observe that 
\begin{eqnarray*}
\mathring{e}_{1} &:&=\frac{\partial }{\partial x^{1}},\text{ }\mathring{e}%
_{2}:=\frac{\partial }{\partial \eta ^{2}}+\eta ^{n+2}\frac{\partial }{%
\partial \tau },..,\mathring{e}_{n}:=\frac{\partial }{\partial \eta ^{n}}%
+\eta ^{2n}\frac{\partial }{\partial \tau }, \\
\mathring{e}_{n+1} &:&=\frac{\partial }{\partial \eta ^{n+1}}-2x^{1}\frac{%
\partial }{\partial \tau }, \\
\mathring{e}_{n+2} &:&=\frac{\partial }{\partial \eta ^{n+2}}-\eta ^{2}\frac{%
\partial }{\partial \tau },..,\mathring{e}_{2n}:=\frac{\partial }{\partial
\eta ^{2n}}-\eta ^{n}\frac{\partial }{\partial \tau }
\end{eqnarray*}

\noindent form an orthonormal basis of $\ker \Theta $ with respect to the
Levi metric $\frac{1}{2}d\Theta (\cdot ,J\cdot )$ (see Subsection B in the
Appendix for more explanation). Let $\phi $ $:=$ $u(\eta ^{2},\eta
^{3},..,\eta ^{2n},\tau )$ $-$ $x^{1}.$ We compute%
\begin{eqnarray}
\mathring{e}_{1}\phi &=&-1,\mathring{e}_{2}\phi =\partial _{\eta ^{2}}u+\eta
^{n+2}\partial _{\tau }u,...,  \label{3.36} \\
\mathring{e}_{n+1}\phi &=&\partial _{\eta ^{n+1}}u-2x^{1}\partial _{\tau
}u,..,\mathring{e}_{2n}\phi =\partial _{\eta ^{2n}}u-\eta ^{n}\partial
_{\tau }u.  \notag
\end{eqnarray}

\noindent Hence $d\phi $ $=$ ($\mathring{e}_{1}\phi )dx^{1}+\sum_{j=2}^{2n}(%
\mathring{e}_{j}\phi )d\eta ^{j}+(\partial _{\tau }\phi )\Theta $ has the
length with respect to the subriemannian metric (\ref{H2}) as follows:%
\begin{equation}
|d\phi |_{\ast }=\sqrt{\sum_{j=1}^{2n}(\mathring{e}_{j}\phi )^{2}}=\sqrt{%
1+|W^{u}(u)|^{2}}  \label{3.36-1}
\end{equation}

\noindent on the graph described by $x^{1}$ $=$ $u(\eta ^{2},\eta
^{3},..,\eta ^{2n},\tau ).$ Here $W^{u}:=(\mathring{e}_{2},..,\mathring{e}%
_{n},\mathring{e}_{n+1}^{u},\mathring{e}_{n+2},..,\mathring{e}_{2n})$ and $%
\mathring{e}_{n+1}^{u}$ $:=$ $\frac{\partial }{\partial \eta ^{n+1}}-2u\frac{%
\partial }{\partial \tau }.$ Note that the standard volume form $dV_{H_{n}}$
of $H_{n}$ equals $dx^{1}\wedge d\eta ^{2}\wedge ..\wedge d\eta ^{2n}\wedge
\Theta .$ On $\Sigma $ $:=$ $\{\phi =0\}$ we have $0$ $=$ $d\phi $ $=$ ($%
\mathring{e}_{1}\phi )dx^{1}+\sum_{j=2}^{2n}(\mathring{e}_{j}\phi )d\eta
^{j}+(\partial _{\tau }\phi )\Theta .$ So by (\ref{3.36}) we get%
\begin{equation}
dx^{1}=\sum_{j=2}^{2n}(\mathring{e}_{j}\phi )d\eta ^{j}+(\partial _{\tau
}\phi )\Theta .  \label{3.37}
\end{equation}%
\noindent Now we compute the area (or volume) element $dv_{\phi }$ for the
hypersurface $\{\phi =\phi (p_{0})\}:$%
\begin{eqnarray}
dv_{\phi } &:&=\frac{d\phi }{|d\phi |_{\ast }}\lrcorner \text{ }dV_{H_{n}}
\label{3.38} \\
&=&\frac{1}{|d\phi |_{\ast }}[(\mathring{e}_{1}\phi )d\eta ^{2}\wedge
..\wedge d\eta ^{2n}\wedge \Theta  \notag \\
&&+(-1)^{j-1}\sum_{j=2}^{2n}(\mathring{e}_{j}\phi )dx^{1}\wedge d\eta
^{2}\wedge .d\hat{\eta}^{j}.\wedge d\eta ^{2n}\wedge \Theta ]  \notag \\
&=&\frac{1}{|d\phi |_{\ast }}[-1-\sum_{j=2}^{2n}(\mathring{e}_{j}\phi
)^{2}]d\eta ^{2}\wedge ..\wedge d\eta ^{2n}\wedge \Theta \text{ (by (\ref%
{3.37}))}  \notag \\
&=&-|d\phi |_{\ast }d\eta ^{2}\wedge ..\wedge d\eta ^{2n}\wedge d\tau  \notag
\\
&=&-\sqrt{1+|W^{u}(u)|^{2}}d\eta ^{2}\wedge ..\wedge d\eta ^{2n}\wedge d\tau
\notag
\end{eqnarray}

\noindent by (\ref{3.35-1}) and (\ref{3.36-1}). From (\ref{3.1}) we compute $%
H_{\phi }$ as follows:%
\begin{eqnarray*}
&&d(\frac{d\phi }{|d\phi |_{\ast }}\lrcorner \text{ }dV_{H_{n}}) \\
&=&d\{[\frac{\mathring{e}_{1}\phi }{|d\phi |_{\ast }}dx^{1}+\sum_{j=2}^{2n}%
\frac{\mathring{e}_{j}\phi }{|d\phi |_{\ast }}d\eta ^{j}+\frac{\partial
_{\tau }\phi }{|d\phi |_{\ast }}\Theta ]\lrcorner \text{ }dV_{H_{n}}\} \\
&=&[\mathring{e}_{1}(\frac{\mathring{e}_{1}\phi }{|d\phi |_{\ast }}%
)+\sum_{j=2}^{2n}\mathring{e}_{j}(\frac{\mathring{e}_{j}\phi }{|d\phi
|_{\ast }})]dV_{H_{n}} \\
&=&\sum_{j=2}^{2n}\mathring{e}_{j}(\frac{\mathring{e}_{j}\phi }{|d\phi
|_{\ast }})dV_{H_{n}}
\end{eqnarray*}

\noindent since $|d\phi |_{\ast }$ is independent of $x^{1}$ by (\ref{3.36-1}%
) and hence $\mathring{e}_{1}(\frac{\mathring{e}_{1}\phi }{|d\phi |_{\ast }}%
) $ $=$ $\partial _{x^{1}}(\frac{-1}{|d\phi |_{\ast }})$ $=$ $0.$ So we have 
\begin{eqnarray}
H_{\phi } &=&\sum_{j=2}^{2n}\mathring{e}_{j}(\frac{\mathring{e}_{j}\phi }{%
|d\phi |_{\ast }})  \label{3.39} \\
&=&W^{u}\cdot \left( \frac{W^{u}(u)}{\sqrt{1+|W^{u}(u)|^{2}}}\right)  \notag
\end{eqnarray}

\noindent at ($\eta ^{2},\eta ^{3},..,\eta ^{2n},\tau ,$ $u(\eta ^{2},\eta
^{3},..,\eta ^{2n},\tau ))$ (cf. (\ref{I6})).

\bigskip

\proof
(\textbf{of Theorem E) }We observe that an intrinsic graph is congruent with
an $l_{a}$\textbf{\ }graph by a rotation. Suppose we have an intrinsic graph
described by (\ref{3.35-0}). Define a rotation $\Psi $ $:$ $(x^{1},$ $..,$ $%
x^{2n},$ $z)$ $\in $ $H_{n}$ $\rightarrow $ $(\tilde{x}^{1},$ $..,$ $\tilde{x%
}^{2n},$ $\tilde{z})$ $\in $ $H_{n}$ by%
\begin{eqnarray}
\tilde{x}^{1} &=&x^{n+1},\tilde{x}^{n+1}=-x^{1},  \label{3.40} \\
\tilde{x}^{j} &=&x^{j},1\leq j\leq 2n,j\neq 1,j\neq n+1,  \notag \\
\tilde{z} &=&z.  \notag
\end{eqnarray}

\noindent From $u(\eta ^{2},\eta ^{3},..,\eta ^{2n},\tau )$ $-$ $x^{1}$ $=$ $%
0$ and $u_{\eta ^{n+1}}$\textit{\ }$\neq $\textit{\ }$0$ at $p_{0},$ we can
express $\eta ^{n+1}$ $=$ $\eta ^{n+1}(x^{1},$ $\eta ^{2},$ $..,$ $\hat{\eta}%
^{n+1},$ $..,$ $\eta ^{2n},$ $\tau )$ near $(u(p_{0}),$ $\eta _{0}^{2},$ $%
.., $ $\hat{\eta}_{0}^{n+1},$ $..,$ $\eta _{0}^{2n},$ $\tau _{0})$ by
implicit function theorem. Here $\hat{\eta}^{n+1}$ means $\eta ^{n+1}$
deleted. In view of (\ref{3.40}) and (\ref{3.35-0}), we get%
\begin{eqnarray}
\tilde{x}^{1} &=&\eta ^{n+1}(x^{1},\eta ^{2},..,\hat{\eta}^{n+1},..,\eta
^{2n},\tau ),  \label{3.41} \\
\tilde{x}^{j} &=&\eta ^{j},\text{ }2\leq j\leq 2n,\text{ }j\neq n+1,  \notag
\\
\tilde{x}^{n+1} &=&-x^{1},  \notag \\
\tilde{z} &=&z=\tau +x^{1}\eta ^{n+1}(x^{1},\eta ^{2},..,\hat{\eta}%
^{n+1},..,\eta ^{2n},\tau ).  \notag
\end{eqnarray}

\noindent Let $\zeta =-x^{1}$ and%
\begin{eqnarray}
&&\tilde{\eta}^{n+1}(\eta ^{2},..,\zeta ,\text{ }..,\eta ^{2n},\tau )\text{ (%
}\hat{\eta}^{n+1}\text{ is replaced by }\zeta )  \label{3.42} \\
&:&=\eta ^{n+1}(x^{1},\eta ^{2},..,\hat{\eta}^{n+1},..,\eta ^{2n},\tau ). 
\notag
\end{eqnarray}

\noindent It follows from (\ref{3.41}) and (\ref{3.42}) that the image of an
intrinsic graph under the rotation $\Psi $ can be depicted as $(\tilde{\eta}%
^{n+1}(\eta ^{2},..,\zeta ,$ $..,\eta ^{2n},\tau ),$ $0,$ $..,$ $0)\circ (0,$
$\eta ^{2},..,\zeta ,$ $..,\eta ^{2n},\tau ),$ an $l_{a}$ graph. Let $\tilde{%
\xi}^{n+1}$ denote the $l_{a}$ graph corresponding to the intrinsic graph $v$
under the rotation $\Psi .$ Near $p_{0}$ the condition $v$\textit{\ }$\geq $%
\textit{\ }$u$ implies $\tilde{\xi}^{n+1}$ $\geq $ $\tilde{\eta}^{n+1}$ ($%
\tilde{\xi}^{n+1}$ $\leq $ $\tilde{\eta}^{n+1},$ resp.) if at $p_{0},$ $\nu
_{\eta ^{n+1}}$\textit{\ }$=$\textit{\ }$u_{\eta ^{n+1}}$\textit{\ }$<$%
\textit{\ }$0$ ($>0,$ resp.). On the other hand, $H_{\tilde{\xi}^{n+1}}$ ($%
H_{\tilde{\eta}^{n+1}},$ resp.) is the same as $H_{v}$ ($H_{u},$ resp.) at
the same point in the graph for the case of $\nu _{\eta ^{n+1}}$\textit{\ }$%
= $\textit{\ }$u_{\eta ^{n+1}}$\textit{\ }$<$\textit{\ }$0$ at $p_{0}.$
Therefore the condition $H_{v}$ $\leq $ $H_{u}$ is reduced to $H_{\tilde{\xi}%
^{n+1}}$ $\leq $ $H_{\tilde{\eta}^{n+1}}$ in some corresponding small
neighborhood when either one is constant$.$\ For the case of $\nu _{\eta
^{n+1}}$\textit{\ }$=$\textit{\ }$u_{\eta ^{n+1}}$\textit{\ }$>$\textit{\ }$%
0 $ at $p_{0},$ $H_{\tilde{\xi}^{n+1}}$ ($H_{\tilde{\eta}^{n+1}},$ resp.) is
the same as $-H_{v}$ (-$H_{u},$ resp.) at the same point in the graph. So $%
H_{v}$ $\leq $ $H_{u}$ is equivalent to $H_{\tilde{\xi}^{n+1}}$ $\geq $ $H_{%
\tilde{\eta}^{n+1}}$ in some corresponding small neighborhood when either
one is constant$.$ Now the theorem follows from Corollary D$^{\prime }.$

\endproof%

\bigskip

We remark that if both $v$\ and $u$\ do not have constant$.$horizontal (or $%
p $-)mean curvature, then we won't be able to reduce $H_{v}$ $\leq $ $H_{u}$
to $H_{\tilde{\xi}^{n+1}}$ $\leq $ (or $\geq )$ $H_{\tilde{\eta}^{n+1}}$ in
general. The reason is that we are comparing horizontal (or $p$-)mean
curvature at different pairs of points on two hypersurfaces.

\bigskip

\section{\textbf{Proofs of Theorem F, Corollaries G, H, and I}}

\proof
(\textbf{of Theorem F}) We want to show that $H_{\vec{F}}(u)$ at $(x^{1},$
.., $x^{m})$ is the same as $H_{\phi }$ for $\phi $ $=$ $u(x^{1},$ .., $%
x^{m})$ $-$ $x^{m+1}$ at $(x^{1},$ .., $x^{m},$ $u(x^{1},$ .., $x^{m}))$
with respect to a certain subriemannian structure on $H_{n}.$ Write $\vec{F}$
$=$ $(F_{1},$ $...,$ $F_{m}).$ Let%
\begin{equation}
\omega ^{m+1}=dx^{m+1}+\sum_{j=1}^{m}F_{j}dx^{j}\text{.}  \label{F0}
\end{equation}

\noindent Define subriemannian metric $<\cdot ,\cdot >_{F}$ on $T^{\ast
}H_{n}$ by%
\begin{eqnarray*}
&<&dx^{i},dx^{j}>_{F}=\delta _{ij},\text{ }1\leq i,j\leq m, \\
&<&\eta ,\omega ^{m+1}>_{F}=0
\end{eqnarray*}

\noindent for any one form $\eta $. The vectors dual to $dx^{1},$ .., $%
dx^{m},$ $\omega ^{m+1}$ read%
\begin{eqnarray*}
e_{j}^{F} &=&\frac{\partial }{\partial x^{j}}-F_{j}\frac{\partial }{\partial
x^{m+1}},\text{ }1\leq j\leq m \\
e_{m+1} &=&\frac{\partial }{\partial x^{m+1}}.
\end{eqnarray*}

\noindent Take the standard volume form $dV_{H_{n}}$ $:=$ $dx^{1}\wedge
...\wedge dx^{m}\wedge \omega ^{m+1}$ $=$ $dx^{1}\wedge ...\wedge
dx^{m}\wedge dx^{m+1}.$ We compute $H_{\phi }$ as follows:%
\begin{eqnarray}
&&d(\frac{d\phi }{|d\phi |_{<\cdot ,\cdot >_{F}}}\lrcorner _{<\cdot ,\cdot
>_{F}}dV_{H_{n}})  \label{F1} \\
&=&\sum_{j=1}^{m}e_{j}^{F}(\frac{e_{j}^{F}(\phi )}{|d\phi |_{<\cdot ,\cdot
>_{F}}})dV_{H_{n}}  \notag
\end{eqnarray}

\noindent where 
\begin{equation*}
|d\phi |_{<\cdot ,\cdot >_{F}}=(\sum_{j=1}^{m}(e_{j}^{F}(\phi ))^{2})^{1/2}.
\end{equation*}

\noindent From (\ref{F1}) we have%
\begin{equation}
H_{\phi }=\sum_{j=1}^{m}e_{j}^{F}(\frac{e_{j}^{F}(\phi )}{|d\phi |_{<\cdot
,\cdot >_{F}}}).  \label{F2}
\end{equation}

\noindent On the other hand, we compute, for $\phi $ $=$ $u(x^{1},$ .., $%
x^{m})$ $-$ $x^{m+1},$%
\begin{equation*}
e_{j}^{F}(\phi )=\frac{\partial u}{\partial x^{j}}+F_{j}.
\end{equation*}

\noindent Hence we have%
\begin{equation*}
|d\phi |_{<\cdot ,\cdot >_{F}}=(\sum_{j=1}^{m}(\frac{\partial u}{\partial
x^{j}}+F_{j})^{2})^{1/2}=|\nabla u+\vec{F}|_{R^{m}}.
\end{equation*}

\noindent Observe that $\frac{e_{j}^{F}(\phi )}{|d\phi |_{<\cdot ,\cdot
>_{F}}}$ is independent of $x^{m+1}.$ It follows that 
\begin{equation}
\sum_{j=1}^{m}e_{j}^{F}(\frac{e_{j}^{F}(\phi )}{|d\phi |_{<\cdot ,\cdot
>_{F}}})=\func{div}(\frac{\nabla u+\vec{F}}{|\nabla u+\vec{F}|_{R^{m}}}).
\label{F3}
\end{equation}

\noindent From (\ref{F2}), (\ref{F3}), we have shown $H_{\vec{F}}(u)$ at $%
(x^{1},$ .., $x^{m})$ equals $H_{\phi }$ at $(x^{1},$ .., $x^{m},$ $u(x^{1},$
.., $x^{m}))$. It is obvious that the translation along $x^{m+1}$ preserves $%
<\cdot ,\cdot >_{F}$ and $dV_{H_{n}}.$ Moreover, $x^{1},$ .., $x^{m},$ $%
x^{m+1}$ are translation-isometric coordinates with respect to the $%
x^{1}..x^{m}$ hyperplane. Observe that $N_{1}^{\perp }(u),$ $N_{2}^{\perp
}(u),$ ..., $N_{m-1}^{\perp }(u)$ are the $x^{1}..x^{m}$ hyperplane
projection of (a choice of) $X_{1},$ .., $X_{m-1}$ in (\ref{3.23}), resp..
The conclusion follows from Theorem C$^{\prime }$ for $l$ $=$ $1.$

\endproof%

\bigskip

\proof
(\textbf{of Corollary G}) Let 
\begin{eqnarray*}
\Theta _{v} &:&=dv+F_{j}dx^{j} \\
&=&\sum_{j=1}^{m}(v_{j}+F_{j})dx^{j}.
\end{eqnarray*}

\noindent Observe that $\Theta _{v}(N_{\gamma }^{\perp }(v))$ $=$ $0$ since $%
N_{\alpha }^{\perp }(v)$ is perpendicular to $\nabla v+\vec{F}.$ It follows
that%
\begin{eqnarray}
-\Theta _{v}([N_{\alpha }^{\perp }(v),N_{\beta }^{\perp }(v)]) &=&d\Theta
_{v}(N_{\alpha }^{\perp }(v),N_{\beta }^{\perp }(v))  \label{2.20} \\
&=&\sum_{k,j=1}^{m}(\partial _{k}F_{j})dx^{k}\wedge dx^{j}(N_{\alpha
}^{\perp }(v),N_{\beta }^{\perp }(v))  \notag \\
&=&\sum_{k,j=1}^{m}(\partial _{k}F_{j})(b_{\alpha }^{k}b_{\beta
}^{j}-b_{\alpha }^{j}b_{\beta }^{k})  \notag \\
&=&\sum_{k<j}(\partial _{k}F_{j}-\partial _{j}F_{k})(b_{\alpha }^{k}b_{\beta
}^{j}-b_{\alpha }^{j}b_{\beta }^{k})  \notag
\end{eqnarray}

\noindent by (\ref{1.1}). By assumption we get $\Theta _{v}([N_{\alpha
}^{\perp }(v),N_{\beta }^{\perp }(v)])$ $\neq $ $0$ in view of (\ref{2.20}).
So $[N_{\alpha }^{\perp }(v),N_{\beta }^{\perp }(v)]$ is nonzero and
transversal to the kernel of $\Theta _{v},$ the space spanned by all $%
N_{\gamma }^{\perp }(v),$ $\gamma $ $=$ $1,$ $2,$ ..., $m-1.$ We then have
that the rank of $\pounds (N_{1}^{\perp }(v),N_{2}^{\perp
}(v),...,N_{m-1}^{\perp }(v)))$ is constant $m$ in $U$.\ Now the conclusion
follows from Theorem F.

\endproof%

\bigskip

We now want to find an intrinsic criterion for the Lie span condition $rank(%
\pounds (N_{1}^{\perp }(u),$ $N_{2}^{\perp }(u),$ ..., $N_{m-1}^{\perp
}(u))) $ $=$ $m$ in Theorem F or more generally $rank(\pounds (\xi \cap
T\Sigma _{1}))$\textit{\ }$=$\textit{\ }$m$ in Theorem C to hold. Suppose $%
\Sigma $ is a hypersurface in a subriemannian manifold $(M,$ $<\cdot ,\cdot
>^{\ast }) $ of dimension $m+1.$ Define $G$ $:$ $T^{\ast }M$ $\rightarrow $ $%
TM$ by $\omega (G(\eta ))$ $=$ $<\omega ,\eta >^{\ast }$ for $\omega ,$ $%
\eta $ $\in $ $T^{\ast }M.$ Let $\xi $ $:=$ $Range(G)$ $\subset $ $TM.$ We
assume (as always) 
\begin{equation*}
\dim \xi =\text{constant }m+1-l
\end{equation*}%
\noindent where $l$ = $\dim \ker G.$ We wonder when%
\begin{equation*}
\pounds (\xi \cap T\Sigma )=\pounds (X_{1},...,X_{m-l}),
\end{equation*}%
\noindent where $X_{1},$..., $X_{m-l}$ form a basis of local sections of $%
\xi \cap T\Sigma $ (near $p_{0}$ where $\xi $ $\neq $ $T\Sigma ),$ has rank $%
m.$ Let us start with the case $l$ $=$ $1.$ Take a nonzero $1$-form $\theta $
$\in $ $\ker G.$ It is easy to see that $\theta $ annihilates vectors in $%
\xi .$ So%
\begin{equation}
\theta (X_{j})=0  \label{F4}
\end{equation}

\noindent for $1$ $\leq $ $j$ $\leq $ $m-1.$ We want Lie bracket of a pair
of $X_{j}$ to generate a direction not in $\xi $ (which will imply $\pounds %
(X_{1},...,X_{m-1})$ has rank $m).$ Suppose the converse holds. Then we have%
\begin{eqnarray}
&&d\theta (X_{i},X_{j})  \label{F5} \\
&=&X_{i}(\theta (X_{j}))-X_{j}(\theta (X_{i}))-\theta ([X_{i},X_{j}])  \notag
\\
&=&-\theta ([X_{i},X_{j}])=0  \notag
\end{eqnarray}

\noindent by (\ref{F4}) for $1$ $\leq $ $j$ $\leq $ $m-1$. It follows from (%
\ref{F5}) that $d\theta |_{\xi \times \xi },$ the bilinear form $d\theta $
restricted to $\xi \times \xi ,$ has rank $\leq $ 2. So we have proved the
following result.

\bigskip

\textbf{Proposition 4.1}. \textit{Let }$\Sigma $\textit{\ be a hypersurface
in a subriemannian manifold }$(M,$\textit{\ }$<\cdot ,\cdot >^{\ast })$%
\textit{\ of dimension }$m+1.$\textit{\ Suppose }$\dim \xi $\textit{\ }$=$%
\textit{\ constant }$m$\textit{\ (i.e., }$l$\textit{\ }$=$\textit{\ }$1).$%
\textit{\ Assume}%
\begin{equation}
rank(d\theta |_{\xi \times \xi })\geq 3.  \label{F6}
\end{equation}%
\textit{\ \noindent Then there holds (\ref{I2-2}) and hence we have}%
\begin{equation*}
rank(\pounds (\xi \cap T\Sigma ))=m.
\end{equation*}

\bigskip

We remark that condition (\ref{F6}) is independent of the choice of $\theta
. $ Since $\dim \ker G$ $=$ $1,$ another nonzero choice $\tilde{\theta}$ $%
\in $ $\ker G$ will be a nonzero multiple of $\theta .$ That is, $\tilde{%
\theta}$ $=$ $\lambda \theta $ with $\lambda $ $\neq $ $0.$ It follows that 
\begin{eqnarray*}
d\tilde{\theta}|_{\xi \times \xi } &=&\lambda d\theta |_{\xi \times \xi
}+(d\lambda \wedge \theta )|_{\xi \times \xi } \\
&=&\lambda d\theta |_{\xi \times \xi }.
\end{eqnarray*}%
\textit{\noindent }So $d\tilde{\theta}|_{\xi \times \xi }$ has the same rank
as $d\theta |_{\xi \times \xi }.$

\bigskip

\proof
(\textbf{of Corollary H}) From the proof of Theorem F, we learn that $H_{%
\vec{F}}(u)$ at $(x^{1},$ .., $x^{m})$ is the same as $H_{\phi }$ for $\phi $
$=$ $u(x^{1},$ .., $x^{m})$ $-$ $x^{m+1}$ at $(x^{1},$ .., $x^{m},$ $%
u(x^{1}, $ .., $x^{m}))$ with respect to a certain subriemannian structure
on $H_{n}$ with $\ m$ $=$ $2n.$ This subriemannian structure has the
property that $\dim \xi $\textit{\ }$=$\textit{\ constant }$m$\textit{\
(i.e., }$l$\textit{\ }$=$\textit{\ }$1).$ We take 
\begin{equation*}
\theta =\omega ^{m+1}=dx^{m+1}+\sum_{j=1}^{m}F_{j}dx^{j}
\end{equation*}%
\textit{\noindent }(see (\ref{F0})). Compute%
\begin{equation*}
d\theta =\sum_{k<j}(\partial _{k}F_{j}-\partial _{j}F_{k})dx^{k}\wedge
dx^{j}.
\end{equation*}%
\textit{\noindent }Therefore\textit{\ }condition\textit{\ }(\ref{1.2.1})
means $rank(d\theta |_{\xi \times \xi })$ $\geq $ $3.$ Note that $%
N_{1}^{\perp }(v),$ $N_{2}^{\perp }(v),$ ..., $N_{m-1}^{\perp }(v)$ are the $%
x^{1}..x^{m}$ hyperplane projection of (a choice of) $X_{1},$ .., $X_{m-1},$
resp. in (\ref{3.23}) for $\phi $ $=$ $v-x^{m+1}$. The conclusion follows
from Proposition 4.1 and Theorem F.

\endproof%

\bigskip

\proof
(\textbf{of Corollary I}) For such an $\vec{F}$ $=$ $(-x^{2},$ $x^{1},$ $..,$
$-x^{2n},$ $x^{2n-1}),$ we compute%
\begin{eqnarray}
\partial _{2k-1}F_{2k}-\partial _{2k}F_{2k-1} &=&2\text{ for }k=1,...,n;
\label{2.21} \\
\partial _{k}F_{j}-\partial _{j}F_{k} &=&0\text{ }(k<j)\text{ otherwise}. 
\notag
\end{eqnarray}%
\noindent So the rank of the matrix ($\partial _{k}F_{j}-\partial _{j}F_{k})$
equals $2n$ by (\ref{2.21}). By assumption we have $m$ $=$ $2n$ $\geq $ $4$ $%
>$ $3.$ The conclusion follows from Corollary H.

\endproof%

\bigskip

We are going to give a sufficient condition for the rank estimate in the
case of general $l.$ Let $\theta ^{1},$ $..,$ $\theta ^{l}$ be a basis of $%
\ker G.$ Choose an (codimension 1, resp.) orthonormal basis $X_{1},$ $%
..,X_{m+1-l}$ ($X_{1},$ $..,X_{m-l},$ resp.) of $\xi $ and its dual forms $%
\omega ^{1},$ $..,$ $\omega ^{m+1-l}$ ($\omega ^{1},$ $..,$ $\omega ^{m-l},$
resp.) (not unique). Then we can find unique vector fields $T_{1},$ $..,$ $%
T_{l}$ such that \{$X_{1},$ $..,X_{m+1-l},$ $T_{1},$ $..,$ $T_{l}\}$ is dual
to $\{\omega ^{1},$ $..,$ $\omega ^{m+1-l},$ $\theta ^{1},$ $..,$ $\theta
^{l}\}.$ Write%
\begin{equation}
\lbrack X_{i},X_{j}]=\sum_{\alpha =1}^{l}W_{ij}^{\alpha }T_{\alpha }\text{ }%
\func{mod}\text{ }\xi  \label{F7}
\end{equation}

\noindent for $1$ $\leq $ $i.j$ $\leq $ $m+1-l$ ($1$ $\leq $ $i.j$ $\leq $ $%
m-l,$ resp.)$.$ For fixed $\alpha $ we view ($W_{ij}^{\alpha })$ as an $%
(m+1-l)\times (m+1-l)$ matrix. Since it is skew-symmetric, it has $\frac{%
(m+1-l)(m-l)}{2}$ degrees of freedom.

\bigskip

\textbf{Proposition 4.2. }\textit{Suppose }%
\begin{eqnarray}
\frac{(m+1-l)(m-l)}{2} &\geq &l  \label{F8} \\
(\frac{(m-l)(m-l-1)}{2} &\geq &l,\text{ resp.)}  \notag
\end{eqnarray}%
\textit{\noindent Then }$rank(\pounds (\xi ))$\textit{\ }$=$\textit{\ }$m+1$%
\textit{\ (}$rank(\pounds (X_{1},...,X_{m-l}))$ $\geq $ $m$\textit{, resp.)
if }$d\theta ^{\alpha }$\textit{\ }$\lrcorner $\textit{\ }$\omega ^{1}\wedge 
$\textit{\ }$..$\textit{\ }$\wedge \omega ^{m+1-l}$\textit{\ (}$\omega
^{1}\wedge $\textit{\ }$..$\textit{\ }$\wedge \omega ^{m-l},$\textit{\
resp.), }$\alpha $\textit{\ }$=$\textit{\ }$1,$\textit{\ }$..,$\textit{\ }$%
l, $\textit{\ are linearly independent, or equivalently }$[W_{ij}^{\alpha }]$%
\textit{\ as (}$m+1-l)$\textit{\ }$\times $\textit{\ (}$m+1-l)$\textit{\
matrix ((}$m-l)$\textit{\ }$\times $\textit{\ (}$m-l)$\textit{\ matrix,
resp.)}$,$\textit{\ }$\alpha $\textit{\ }$=$\textit{\ }$1,$\textit{\ }$..,$%
\textit{\ }$l,$\textit{\ are linearly independent.}

\bigskip

\proof
In view of (\ref{F7}) and (\ref{F8}), the linear independence of $%
[W_{ij}^{\alpha }],$ $\alpha $\textit{\ }$=$\textit{\ }$1,$\textit{\ }$..,$%
\textit{\ }$l,$ implies $T_{1},$ $..,$ $T_{l}$ $\in $ $\pounds (\xi )$ ($%
\pounds (X_{1},...,X_{m-l}),$ resp.) by basic linear algebra. It follows
that $rank(\pounds (\xi ))$\textit{\ }$\geq $\textit{\ }($m+1-l)$ $+$ $l$ $=$
$m+1$\textit{\ }($rank(\pounds (X_{1},...,X_{m-l}))$ $\geq $ ($m-l)$ $+$ $l$ 
$=$ $m,$ resp.). On the other hand, clearly we have $rank(\pounds (\xi ))$ $%
\leq $ $m+1.$ By the formula 
\begin{equation*}
d\theta ^{\alpha }(X_{i},X_{j})=X_{i}(\theta ^{\alpha }(X_{j}))-X_{j}(\theta
^{\alpha }(X_{i}))-\theta ^{\alpha }([X_{i},X_{j}]),
\end{equation*}

\noindent we obtain

\begin{equation*}
d\theta ^{\alpha }=-\sum_{i<j}W_{ij}^{\alpha }\omega ^{i}\wedge \omega ^{j}%
\text{ }\func{mod}\text{ }\theta ^{1},..,\theta ^{l}
\end{equation*}

\noindent by (\ref{F7}). We then learn that two conditions are equivalent.

\endproof%

\bigskip

\section{Singularities and proofs of Theorems $C^{\prime \prime }$, J, K, $%
\tilde{C}^{\prime \prime },$J$^{\prime }$}

\proof
\textbf{(of Theorem }$C^{\prime \prime }$\textbf{) }By Theorem $5.2^{\prime
} $ in \cite{chmy}, we have $N_{\vec{F}}(v)$ $=$ $N_{\vec{F}}(u)$ in $\Omega
^{+}\backslash \{S_{\vec{F}}(u)\cup S_{\vec{F}}(v)\}$ where $\Omega ^{+}$ $%
:= $ $\{p$ $\in $ $\Omega $ $|$ $u(p)-v(p)$ $>$ $0\}.$ Suppose $\Omega ^{+}$
is not empty. Then applying Lemma $5.3^{\prime }$ in \cite{chmy} with $%
\Omega $ $=$ $\Omega ^{+}$ and $\vec{F}^{\ast }$ replaced by $\vec{F}^{b}$
in the proof$,$ we obtain $v$ $=$ $u$ in $\Omega ^{+},$ a contradiction. We
have proved $v\geq u$\textit{\ in }$\Omega .$

\endproof%

\bigskip

\proof
\textbf{(of Theorem J)} Case 1: Suppose $v$ $=$ $u$ at $q$ $\neq $ $p_{0}.$
Since $p_{0}$ is the only singular point of $u$ ($v,$ resp.)$,$ $q$ is
nonsingular with respect to $u$ ($v,$ resp.)$.$ On the other hand, that $v$ $%
-$ $u$ $\geq $ 0 and $v$ $-$ $u$ $=$ $0$ at $q$ implies that $\nabla (v-u)$ $%
=$ $0$ at $q.$ It follows that at $q$, 
\begin{eqnarray*}
\nabla v+\vec{F} &=&\nabla u+\vec{F}\neq 0 \\
(\nabla u+\vec{F} &=&\nabla v+\vec{F}\neq 0,\text{ resp.)}
\end{eqnarray*}%
\noindent since $q$ is nonsingular with respect to $u$ ($v,$ resp.)$.$ So $q$
is also nonsingular with respect to $v$ ($u,$ resp.)$.$ Now by Theorem F we
obtain $v$ $\equiv $ $u$ in a connected component $W$ of nonsingular (with
respect to both $v$ and $u)$ set, containing $q.$ We claim $W$ $=$ $\Omega
\backslash \{p_{0}\}.$ Otherwise there is a point $q^{\prime }$ $\in $ $(S_{%
\vec{F}}(v)$ $\cap $ $\bar{W})\backslash \{p_{0}\}$ ($(S_{\vec{F}}(u)$ $\cap 
$ $\bar{W})\backslash \{p_{0}\},$ resp.) at which $\nabla u+\vec{F}$ $=$ $%
\nabla v+\vec{F}$ $=$ $0$ ($\nabla v+\vec{F}$ $=$ $\nabla u+\vec{F}$ $=$ $0,$
resp.)$.$ Therefore $q^{\prime }$ $\in $ $S_{\vec{F}}(u)$ ($q^{\prime }$ $%
\in $ $S_{\vec{F}}(v),$resp.)$,$ a contradiction to $p_{0}$ being the only
singular point of $u$ ($v,$ resp.).\ Hence $v$ $\equiv $ $u$ in $\Omega .$

Case 2: Suppose $v$ $>$ $u$ in $\Omega \backslash \{p_{0}\}.$ So there is a
small ball $B$ centered at $p_{0}$ such that $v$ $>$ $u$ in $B\backslash
\{p_{0}\}.$ It follows that 
\begin{equation*}
v\geq u+c
\end{equation*}

\noindent on $\partial B$ for some constant $c$ $>$ $0.$ By Theorem $%
C^{\prime \prime }$ (a version of the usual maximum principle), we conclude
that $v\geq u+c$ in $B.$ But $v(p_{0})$ $=$ $u(p_{0})$ implies $0$ $\geq $ $%
c,$ a contradiction. We have shown the impossibility of this case.

\endproof%

\bigskip

\proof
\textbf{(of Theorem K)} We first observe that the assumption $\mathcal{H}%
_{m-1}(\overline{S_{\vec{F}}(u)})$\textit{\ }$=$\textit{\ }$0$\textit{\ }or%
\textit{\ }$\mathcal{H}_{m-1}(\overline{S_{\vec{F}}(v)})$\textit{\ }$=$%
\textit{\ }$0$ in Theorem J is satisfied for $m$ $=$ $2n$ and $\vec{F}$ $=$ $%
(-x^{2},$ $x^{1},$ ..., $-x^{2n},$ $x^{2n-1})$ by Lemma 5.4 in \cite{chmy}.
Also $\func{div}\vec{F}^{b}$\textit{\ = }$2n$ $>$ $0.$ In view of Corollary
D, Corollary I, and Theorem J, we finally reach a situation that $\Sigma
_{1} $ and $\Sigma _{2}$ are tangent at a nonisolated singular point for
both $\Sigma _{1}$ and $\Sigma _{2}$ if they don't coincide completely$.$
This contradicts the assumption that \textit{either }$\Sigma _{1}$\textit{\
or }$\Sigma _{2}$\textit{\ has only isolated singular points.}

\endproof%

\bigskip

\textbf{Lemma 5.1}. \textit{Suppose }$|d\psi |_{\ast }$\textit{\ }$\neq $%
\textit{\ }$0,$\textit{\ }$|d\phi |_{\ast }$\textit{\ }$\neq $\textit{\ }$0.$%
\textit{\ Then the following formula}%
\begin{eqnarray}
&<&d\psi -d\phi ,\frac{d\psi }{|d\psi |_{\ast }}-\frac{d\phi }{|d\phi
|_{\ast }}>^{\ast }  \label{S0} \\
&=&\frac{1}{2}(|d\psi |_{\ast }+|d\phi |_{\ast })|\frac{d\psi }{|d\psi
|_{\ast }}-\frac{d\phi }{|d\phi |_{\ast }}|_{\ast }^{2}  \notag
\end{eqnarray}

\textit{holds.}

\bigskip

\proof
We compute 
\begin{eqnarray}
&<&d\psi -d\phi ,\frac{d\psi }{|d\psi |_{\ast }}-\frac{d\phi }{|d\phi
|_{\ast }}>^{\ast }  \label{S0-1} \\
&=&|d\psi |_{\ast }+|d\phi |_{\ast }-\frac{<d\phi ,d\psi >^{\ast }}{|d\psi
|_{\ast }}-\frac{<d\psi ,d\phi >^{\ast }}{|d\phi |_{\ast }}  \notag \\
&=&(|d\psi |_{\ast }+|d\phi |_{\ast })(1-\cos \vartheta )  \notag
\end{eqnarray}

\noindent where we write $<d\phi ,d\psi >^{\ast }$ $=$ $<d\psi ,d\phi
>^{\ast }$ $=$ $|d\psi |_{\ast }$ $|d\phi |_{\ast }$ $\cos \vartheta .$ On
the other hand, we have%
\begin{eqnarray}
&&|\frac{d\psi }{|d\psi |_{\ast }}-\frac{d\phi }{|d\phi |_{\ast }}|_{\ast
}^{2}  \label{S0-2} \\
&=&|\frac{d\psi }{|d\psi |_{\ast }}|_{\ast }^{2}+|\frac{d\phi }{|d\phi
|_{\ast }}|_{\ast }^{2}-2\frac{<d\psi ,d\phi >^{\ast }}{|d\psi |_{\ast
}|d\phi |_{\ast }}  \notag \\
&=&2(1-\cos \vartheta )  \notag
\end{eqnarray}

\noindent Now (\ref{S0}) follows from (\ref{S0-1}) and (\ref{S0-2}).

\endproof%

\bigskip

We remark that the formula (\ref{S0}) in vector form first appeared in \cite%
{Mi}, \cite{Hw}, and \cite{CK} independently. The version for Heisenberg
group appeared in Lemma 5.1$^{\prime }$ of \cite{chmy}). For the case of
bounded variation, the reader is referred to \cite{ch1}. For the definition
of $S(u)$ ($S(v),$ resp.), see the paragraph before Theorem J$^{\prime }$ in
Section 1.

\bigskip

\textbf{Lemma 5.2}. \textit{Suppose} $(M,$\textit{\ }$<\cdot ,\cdot >^{\ast
},dv_{M})$\textit{\ of dimension }$m+1$ \textit{has isometric translations }$%
\Psi _{a}$ \textit{near }$p_{0}$\textit{\ }$\in $\textit{\ }$M,$ \textit{%
transversal to a hypersurface }$\Sigma $\textit{\ passing through }$p_{0}.$ 
\textit{Take a system of translation-isometric coordinates} $x^{1},$\textit{%
\ }$x^{2},$\textit{\ ..., }$x^{m+1}$\textit{\ in an open neighborhood }$%
\Omega $\textit{\ of }$p_{0}$ \textit{such that }$x^{m+1}$\textit{\ }$=$%
\textit{\ }$0$\textit{\ on }$\Sigma \cap \Omega .$\textit{\ Let }$V$ $%
\subset $ $\bar{V}$ $\subset $ $\Omega $ \textit{be a smaller open
neighborhood of }$p_{0}.$ \textit{Let }$v,$\textit{\ }$u$\textit{\ }$\in $%
\textit{\ }$C^{2}(\Sigma \cap V)$\textit{\ }$\cap $\textit{\ }$C^{0}(%
\overline{\Sigma \cap V})$\textit{\ define graphs in }$\Omega $ and satisfy%
\begin{equation}
H(v)\leq H(u)\text{ in (}\Sigma \cap V)\backslash \{S(u)\cup S(v)\},
\label{1.4}
\end{equation}%
\begin{equation}
v\geq u\text{ on }\partial (\Sigma \cap V).  \label{1.5}
\end{equation}%
\textit{\noindent Assume }$\mathcal{H}_{m-1}(\overline{S(u)\cup S(v)})$%
\textit{\ }$=$\textit{\ }$0$\textit{. Let }$\psi $ :$=$ $v-x^{m+1},$ $\phi $ 
$:=$ $u-x^{m+1}.$ \textit{Then }$\frac{d\psi }{|d\psi |_{\ast }}$ $=$ $\frac{%
d\phi }{|d\phi |_{\ast }}$ \textit{mod }$\ker G$\textit{\ in (}$\Sigma \cap
V)^{+}\backslash \overline{S(u)\cup S(v)}$ \textit{where (}$\Sigma \cap
V)^{+}$\textit{\ }$:=$\textit{\ }$\{q\in \Sigma \cap V:$\textit{\ }$u(q)$%
\textit{\ }$-$\textit{\ }$v(q)$\textit{\ }$>$\textit{\ }$0\}.$

\bigskip

\proof
First $\mathcal{H}_{m-1}(\overline{S(u)\cup S(v)})$\textit{\ }$=$\textit{\ }$%
0$ means that given any $\varepsilon $ $>$ $0,$ we can find countably many
ball $B_{j,\varepsilon },$ $j$ = $1,$ $2,$ ... such that $\overline{S(u)\cup
S(v)}$ $\subset $ $\cup _{j=1}^{\infty }B_{j,\varepsilon }$ and $%
\sum_{j=1}^{\infty }\mathcal{H}_{m-1}(\partial B_{j,\varepsilon })$ $<$ $%
\varepsilon $ and we can arrange $\cup _{j=1}^{\infty }B_{j,\varepsilon
_{1}} $ $\subset $ $\cup _{j=1}^{\infty }B_{j,\varepsilon _{2}}$ for $%
\varepsilon _{1}$ $<$ $\varepsilon _{2}.$ Since $\overline{S(u)\cup S(v)}$
is compact, we can find finitely many $B_{j,\varepsilon }$ 's, say $j$ $=$ $%
1,$ $2,$ $...,$ $n(\varepsilon ),$ still covering $\overline{S(u)\cup S(v)}$%
. Suppose ($\Sigma \cap V)^{+}$ is not empty. Then by Sard's theorem there
exists a sequence of positive number $\delta _{i}$ converging to 0 as $i$
goes to infinity, such that ($\Sigma \cap V)_{i}^{+}$ $:=$ $\{p\in \Sigma
\cap V:$ $u(p)-v(p)$ $>$ $\delta _{i}\}$ is not empty and $\partial (\Sigma
\cap V)_{i}^{+}\backslash (S(u)\cup S(v))$ is $C^{2}$ smooth. Note that $%
\partial (\Sigma \cap V)_{i}^{+}$ $\cap $ $\partial (\Sigma \cap V)$ $%
\subset $ $S(u)\cup S(v)$ by (\ref{1.5}). Choose $a$ $>$ $0$ (independent of 
$\varepsilon $ and $\delta _{i})$ such that [$(\Sigma \cap
V)_{i}^{+}\backslash $ $\cup _{j=1}^{n(\varepsilon )}B_{j,\varepsilon }]$ $%
\times $ $[-a,a]$ $\subset $ $\Omega .$ Now we consider%
\begin{equation}
I_{\varepsilon }^{i}:=\doint\limits_{\partial ([(\Sigma \cap
V)_{i}^{+}\backslash \cup _{j=1}^{n(\varepsilon )}B_{j,\varepsilon }]\times
\lbrack -a,a])}\tan ^{-1}(\psi -\phi )(\frac{d\psi }{|d\psi |_{\ast }}-\frac{%
d\phi }{|d\phi |_{\ast }})\text{ }\lrcorner \text{ }dv_{M}.  \label{S0-3}
\end{equation}

\noindent By Stokes' theorem we have%
\begin{eqnarray}
I_{\varepsilon }^{i} &=&\int_{[(\Sigma \cap V)_{i}^{+}\backslash \cup
_{j=1}^{n(\varepsilon )}B_{j,\varepsilon }]\times \lbrack -a,a]}\{\frac{%
d\psi -d\phi }{1+(\psi -\phi )^{2}}\wedge \lbrack (\frac{d\psi }{|d\psi
|_{\ast }}-\frac{d\phi }{|d\phi |_{\ast }})\text{ }\lrcorner \text{ }dv_{M}]
\label{S1} \\
&&+\tan ^{-1}(\psi -\phi )(H_{\psi }-H_{\phi })dv_{M}\}  \notag
\end{eqnarray}

\noindent Observe that $\psi -\phi $ $=$ $v-u$ $<$ $-\delta _{i}$ $<$ $0$ in
($\Sigma \cap V)_{i}^{+}$ $\times $ $[-a,a]$ and $H_{\psi }(q)$ $-$ $H_{\phi
}(q)$ = $H(v)(\pi _{\Sigma }(q))$ $-$ $H(u)(\pi _{\Sigma }(q))$ $\leq $ $0$
by Corollary 3.2 and (\ref{1.4}), where $\pi _{\Sigma }$ is the projection
to $\Sigma $ ($(x^{1},$ $..,x^{m},$ $x^{m+1})$ $\rightarrow $ $(x^{1},$ $%
..,x^{m})).$ So the second term in the right hand side of (\ref{S1}) is
nonnegative. As to the first term, we observe that $\eta $ $\wedge $ $%
(\omega $ $\lrcorner $ $dv_{M})$ $=$ $<\eta ,\omega >^{\ast }$ $dv_{M}$ for
any $1$-forms $\eta $ and $\omega .$ It then follows from (\ref{S0}) in
Lemma 5.1 and (\ref{S1}) that 
\begin{equation}
I_{\varepsilon }^{i}\geq \int_{\lbrack (\Sigma \cap V)_{i}^{+}\backslash
\cup _{j=1}^{n(\varepsilon )}B_{j,\varepsilon }]\times \lbrack -a,a]}\frac{%
|d\psi |_{\ast }+|d\phi |_{\ast }}{2[1+(\psi -\phi )^{2}]}|\frac{d\psi }{%
|d\psi |_{\ast }}-\frac{d\phi }{|d\phi |_{\ast }}|_{\ast }^{2}dv_{M}
\label{S2}
\end{equation}

\noindent On the other hand, look at the boundary integral in (\ref{S0-3}): $%
\psi -\phi $ $=$ $v-u$ $=$ $\delta _{i}$ on $\partial $($\Sigma \cap
V)_{i}^{+}\times \lbrack -a_{i},a_{i}]$ $\rightarrow $ $0$ as $i\rightarrow
\infty ,$ the boundary integral on $[(\Sigma \cap V)_{i}^{+}\backslash \cup
_{j=1}^{n(\varepsilon )}B_{j,\varepsilon }]$ $\times $ $\{-a_{i}\}$ cancels
with the boundary integral on $[(\Sigma \cap V)_{i}^{+}\backslash \cup
_{j=1}^{n(\varepsilon )}B_{j,\varepsilon }]$ $\times $ $\{a_{i}\}$ due to
translation invariance and orientation, and $\frac{d\psi }{|d\psi |_{\ast }}$
$-$ $\frac{d\phi }{|d\phi |_{\ast }}$ is bounded while $\sum_{j=1}^{n(%
\varepsilon )}\mathcal{H}_{m-1}(\partial B_{j,\varepsilon })$ $<$ $%
\varepsilon .$ So we conclude that $I_{\varepsilon }^{i}\leq \varepsilon $
for $i$ = $i(\varepsilon )$ large enough from (\ref{S0-3}). Together with (%
\ref{S2}) we conclude that%
\begin{equation*}
|\frac{d\psi }{|d\psi |_{\ast }}-\frac{d\phi }{|d\phi |_{\ast }}|_{\ast }=0
\end{equation*}

\noindent It follows that $\frac{d\psi }{|d\psi |_{\ast }}$ $=$ $\frac{d\phi 
}{|d\phi |_{\ast }}$ mod\textit{\ }$\ker G$\textit{\ }in\textit{\ (}$\Sigma
\cap V)^{+}\backslash \overline{S(u)\cup S(v)}.$

\endproof%

\bigskip

\textbf{Lemma 5.3}. \textit{Suppose} $(M,$\textit{\ }$<\cdot ,\cdot >^{\ast
},dv_{M})$\textit{\ of dimension }$m+1$\textit{\ has isometric translations }%
$\Psi _{a}$ \textit{near }$p_{0}$\textit{\ }$\in $\textit{\ }$M,$ \textit{%
transversal to a hypersurface }$\Sigma $\textit{\ passing through }$p_{0}.$ 
\textit{Take a system of translation-isometric coordinates} $x^{1},$\textit{%
\ }$x^{2},$\textit{\ ..., }$x^{m+1}$\textit{\ in an open neighborhood }$%
\Omega $\textit{\ of }$p_{0}$ \textit{such that }$x^{m+1}$\textit{\ }$=$%
\textit{\ }$0$\textit{\ on }$\Sigma \cap \Omega $ \textit{and }$p_{0}$%
\textit{\ is the origin }$(0,$\textit{\ }$..,$\textit{\ }$0).$\textit{\ Let }%
$v,u$\textit{\ : }$\Sigma \cap \Omega $\textit{\ }$\rightarrow $\textit{\ }$%
R $\textit{\ be two graphs in }$\Omega $\textit{.} \textit{Let }$\psi $ :$=$ 
$v(x^{1},\mathit{\ }x^{2},\mathit{\ ..,}$ $x^{m})$ $-$ $x^{m+1},$ $\phi $ $%
:= $ $u(x^{1},\mathit{\ }x^{2},\mathit{\ ..,}x^{m})$ $-$ $x^{m+1}.$ \textit{%
Suppose }$\frac{d\psi }{|d\psi |_{\ast }}$ $=$ $\frac{d\phi }{|d\phi |_{\ast
}}$\textit{\ mod }$\ker G$\textit{. Moreover, we assume the rank condition (%
\ref{I2-2})}$.$\textit{\ Then }$\nabla v$ $=$ $\nabla u$\textit{\ in }$%
\Sigma \cap \Omega .$

\bigskip

\proof
Note that $\xi $ $=$ $\cap _{\theta \in \ker G}\ker \theta .$ So $\frac{%
d\psi }{|d\psi |_{\ast }}=\frac{d\phi }{|d\phi |_{\ast }}$\textit{\ mod }$%
\ker G$ implies $\xi \cap T_{p}\{\psi =c(p)\}$ $=$ $\xi \cap T_{p}\{\phi
=0\}.$ By a translation (depending on $p)$ in the $x^{m+1}$ direction, we
can translate the hypersurface $\{\psi =c(p)\}$ to the fixed hypersurface $%
\{\psi =0\}.$ Since translations are isometries, we then have%
\begin{equation}
\xi \cap T_{(\bar{p},v(\bar{p}))}\{\psi =0\}=\Phi (\bar{p})_{\ast }(\xi \cap
T_{(\bar{p},u(\bar{p}))}\{\phi =0\})  \label{S3}
\end{equation}

\noindent where $\bar{p}$ $\in $ $\Sigma \cap \Omega $ and $\Phi (\bar{p})$
is a translation in the $x^{m+1}$ direction, depending on $\bar{p}.$ By (\ref%
{I2-2}) we have $\pounds (\xi \cap T_{(\bar{p},v(\bar{p}))}\{\psi =0\})$ $=$ 
$\Gamma (T_{(\bar{p},v(\bar{p}))}\{\psi =0\})$ ( $\pounds (\xi \cap T_{(\bar{%
p},u(\bar{p}))}\{\phi =0\})$ $=$ $\Gamma (T_{(\bar{p},u(\bar{p}))}\{\phi
=0\}),$ resp.) since $rank(\pounds (\xi \cap T_{(\bar{p},v(\bar{p}))}\{\psi
=0\}))$ $=$ $m$ ($rank(\pounds (\xi \cap T_{(\bar{p},u(\bar{p}))}\{\phi
=0\}))$ $=$ $m,$ resp.) and $\xi \cap T_{(\bar{p},v(\bar{p}))}\{\psi =0\}$ $%
\subset $ $T_{(\bar{p},v(\bar{p}))}\{\psi =0\}$ ($\xi \cap T_{(\bar{p},u(%
\bar{p}))}\{\phi =0\}$ $\subset $ $T_{(\bar{p},u(\bar{p}))}\{\phi =0\},$
resp.). Here $\Gamma (E)$ denotes the space of all ($C^{\infty }$ smooth)
sections of the vector bundle $E.$ It follows from (\ref{S3}) that%
\begin{equation*}
T_{(\bar{p},v(\bar{p}))}\{\psi =0\}=\Phi (\bar{p})_{\ast }(T_{(\bar{p},u(%
\bar{p}))}\{\phi =0\}).
\end{equation*}

\noindent We then have $\nabla v=\nabla u$ in $\Sigma \cap \Omega .$

\endproof%

\bigskip

\proof
\textbf{(of Theorem \~{C}}$^{\prime \prime })$ By Lemma 5.2 we get $\frac{%
d\psi }{|d\psi |_{\ast }}$ $=$ $\frac{d\phi }{|d\phi |_{\ast }}$ mod $\ker G$%
\ in ($\Sigma \cap V)^{+}\backslash \overline{S(u)\cup S(v)}$ where ($\Sigma
\cap V)^{+}$\ $:=$\ $\{q\in \Sigma \cap V:$\ $u(q)$\ $-$\ $v(q)$\ $>$\ $0\}.$
By Lemma 5.3 we get $\nabla v$ $=$ $\nabla u$ in ($\Sigma \cap V)^{+}$ and
hence $u-v$ $=$ constant $>$ $0$ in ($\Sigma \cap V)^{+}.$ On the other
hand, we have $v\geq u$ on $\partial (\Sigma \cap V)$ by assumption. We get
contradiction by continuity of $v$ and $u.$ So ($\Sigma \cap V)^{+}$ is an
empty set. We then conclude that $v$ $\geq $ $u$ in\textit{\ }$\Sigma \cap
V. $

\endproof%

\bigskip

\proof
\textbf{(of Theorem J}$^{\prime })$ The idea is similar as in the proof of
Theorem J.

Case 1: Suppose $v=u$ at $q$ $\neq $ $p_{0}.$ Observe that $q$ is
nonsingular with respect to $u$ ($v$, resp.) since $p_{0}$ is the only
singular point of $u$ ($v,$ resp.). On the other hand, we have $\nabla (v-u)$
$=$ $0$ at $q$ since $v-u$ $\geq $ $0$ and $v-u$ $=$ $0$ at $q$. It follows
that 
\begin{eqnarray*}
d\psi &=&dv-dx^{m+1} \\
&=&du-dx^{m+1}=d\phi
\end{eqnarray*}

\noindent at $(q,$ $v(q)$ $=$ $(q,u(q)).$ So $\xi $ $\nsubseteq $ $\ker
d\phi $ at $(q,$ $u(q))$ implies $\xi $ $\nsubseteq $ $\ker d\psi $ at $%
(q,v(q)).$ That is, $q$ is also nonsingular with respect to $v$ ($u,$
resp.). By Theorem C$^{\prime }$ we obtain $v\equiv u$ in a connected
component $W$ of nonsingular (with respect to both $v$ and $u)$ set,
containing $q.$ We claim $W$ $=$ ($\Sigma \cap \Omega )\backslash \{p_{0}\}.$
Otherwise there is a point $q^{\prime }$ $\in $ ($S_{\Sigma \cap \Omega
}(v)\cap \bar{W})\backslash \{p_{0}\}$ (($S_{\Sigma \cap \Omega }(u)\cap 
\bar{W})\backslash \{p_{0}\},$ resp.) at which $\nabla v$ $=$ $\nabla u,$
and hence $\xi $ $\subset $ $\ker d\psi $ $=$ $\ker d\phi $ at $(q^{\prime
}, $ $v(q^{\prime })$ $=$ $(q^{\prime },u(q^{\prime })).$ So $q^{\prime }$
is also a singular point of $u$ ($v,$ resp.)$,$ a contradiction to $p_{0}$
being the only singular point of $u$ ($v,$ resp.). Hence $v\equiv u$ in $%
\Sigma \cap \Omega .$

Case 2: Suppose $v$ $>$ $u$ in ($\Sigma \cap \Omega )\backslash \{p_{0}\}.$
So there is a small ball $B$ $\subset $ $\Sigma \cap \Omega ,$ centered at $%
p_{0}$ such that $v$ $>$ $u$ in $B\backslash \{p_{0}\}.$ It follows that 
\begin{equation*}
v\geq u+c
\end{equation*}

\noindent on $\partial B$ for some constant $c$ $>$ $0.$ By Theorem $\tilde{C%
}^{\prime \prime }$, we conclude that $v\geq u+c$ in $B.$ But $v(p_{0})$ $=$ 
$u(p_{0})$ implies $0$ $\geq $ $c,$ a contradiction. We have shown the
impossibility of this case.

\endproof%

\bigskip

\section{Applications: uniqueness and nonexistence}

Consider the Heisenberg cylinder $H_{n}\backslash \{0\}$ with $CR$ structure
same as $H_{n}$ and contact form%
\begin{equation*}
\theta =\frac{1}{\rho ^{2}}\Theta ,
\end{equation*}

\noindent denoted as $(H_{n}\backslash \{0\},$ $\rho ^{-2}\Theta ),$ where $%
\Theta $ $:=$ $dz+\sum_{j=1}^{n}(x_{j}dy_{j}-y_{j}dx_{j})$ and $\rho $ $:=$ $%
[(\sum_{j=1}^{n}(x_{j}^{2}+y_{j}^{2}))^{2}+4z^{2}]^{1/4}.$ Here $x_{1},$ $%
y_{1},$ $..,x_{n},$ $y_{n},$ $z$ denote the coordinates of $H_{n}.$
Topologically $H_{n}\backslash \{0\}$ is homeomorphic to $S^{2n}\times R^{+}$
through the map%
\begin{equation*}
(x_{1},y_{1},..,x_{n},y_{n},z)\rightarrow ((\frac{x_{1}}{\rho },..,\frac{%
y_{n}}{\rho },\frac{z}{\rho ^{2}}),\rho )
\end{equation*}

\noindent where the Heisenberg sphere $S^{2n}$ $\subset $ $H_{n}$ is defined
by $\rho $ $=$ $1.$ Next we want to compute horizontal ($p$-) mean curvature
of a hypersurface of $(H_{n}\backslash \{0\},$ $\rho ^{-2}\Theta ),$
described by a defining function $\phi .$ Take an orthonormal basis $e_{I}$ $%
:=$ $\rho \hat{e}_{I},$ $1\leq I\leq 2n,$ with respect to the Levi metric
(see Subsection B in the Appendix) and $T$ $:=$ $\rho ^{2}\frac{\partial }{%
\partial z}$ where 
\begin{eqnarray*}
\hat{e}_{I} &=&\hat{e}_{j}=\frac{\partial }{\partial x_{j}}+y_{j}\frac{%
\partial }{\partial z},\text{ for }1\leq I=j\leq n \\
\hat{e}_{I} &=&\hat{e}_{j^{\prime }}=\frac{\partial }{\partial y_{j}}-x_{j}%
\frac{\partial }{\partial z},\text{ for }n+1\leq I=j^{\prime }=n+j\leq 2n.
\end{eqnarray*}

\noindent So the dual coframe is $\theta ^{I}$ $=$ $\rho ^{-1}dx^{I}$ ($%
x^{j}=x_{j},$ $x^{n+j}=y_{j}$ for $1\leq j\leq n)$ and $\theta $ $=$ $\rho
^{-2}\Theta $ and the associated subriemannian metric $<\cdot ,\cdot >^{\ast
}$ on $H_{n}\backslash \{0\}$ is given by%
\begin{equation*}
<\theta ^{I},\theta ^{K}>^{\ast }=\delta _{IK},\text{ }<\cdot ,\theta
>^{\ast }=<\theta ,\cdot >^{\ast }=0
\end{equation*}

\noindent (cf. (\ref{H2})). We have the volume form 
\begin{equation}
dV:=\theta ^{1}\wedge \theta ^{2}\wedge ...\wedge \theta ^{2n}\wedge \theta
=\rho ^{-(2n+2)}dx^{1}\wedge dx^{2}\wedge ...\wedge dx^{2n}\wedge \Theta .
\label{App-1}
\end{equation}%
\noindent Compute 
\begin{eqnarray}
d\phi &=&(e_{I}\phi )\theta ^{I}+(T\phi )\theta  \label{App0} \\
&=&(\hat{e}_{I}\phi )dx^{I}+\frac{\partial \phi }{\partial z}\Theta ,  \notag
\end{eqnarray}%
\noindent and hence%
\begin{eqnarray}
|d\phi |_{\ast }^{2} &:&=<d\phi ,d\phi >^{\ast }  \label{App1} \\
&=&\sum_{I=1}^{2n}(e_{I}\phi )^{2}=\rho ^{2}\sum_{I=1}^{2n}(\hat{e}_{I}\phi
)^{2}.  \notag
\end{eqnarray}

From (\ref{App-1}), \ref{App0}, we compute 
\begin{eqnarray}
&&\frac{d\phi }{|d\phi |_{\ast }}\text{ }\lrcorner \text{ }dV  \label{App2}
\\
&=&(-1)^{I-1}\frac{e_{I}\phi }{|d\phi |_{\ast }}\theta ^{1}\wedge ..\wedge 
\hat{\theta}^{I}\wedge ..\wedge \theta ^{2n}\wedge \theta \text{ }  \notag \\
&&\text{(}\hat{\theta}^{I}\text{ means }\theta ^{I}\text{ deleted)}  \notag
\\
&=&(-1)^{I-1}\frac{\hat{e}_{I}\phi }{|d\phi |_{H_{n}}}\rho
^{-(2n+1)}dx^{1}\wedge ..\wedge d\hat{x}^{I}\wedge ..\wedge dx^{2n}\wedge
\Theta  \notag \\
&=&\rho ^{-(2n+1)}\frac{d\phi }{|d\phi |_{H_{n}}}\text{ }\lrcorner _{H_{n}}%
\text{ }dV_{H_{n}}.  \notag
\end{eqnarray}

\noindent where $dV_{H_{n}}$ denotes the standard volume form of $H_{n}$,
which is $dx^{1}\wedge ..\wedge d\hat{x}^{I}\wedge ..\wedge dx^{2n}\wedge
\Theta $ and $|\cdot |_{H_{n}}$ denotes the length with respect to the
standard subriemannian metric in $H_{n}$ (see (\ref{H1}))$.$ It follows that 
$|d\phi |_{H_{n}}$ = ($\sum_{I=1}^{2n}(\hat{e}_{I}\phi )^{2})^{1/2}.$ Taking
exterior differentiation of (\ref{App2}) gives%
\begin{eqnarray}
H_{\phi }dV &:&=d(\frac{d\phi }{|d\phi |_{\ast }}\text{ }\lrcorner \text{ }%
dV)  \label{App3} \\
&=&\rho ^{-(2n+1)}d(\frac{d\phi }{|d\phi |_{H_{n}}}\text{ }\lrcorner _{H_{n}}%
\text{ }dV_{H_{n}})  \notag \\
&&+\frac{\hat{e}_{I}\phi }{|d\phi |_{H_{n}}}\hat{e}_{I}(\rho
^{-(2n+1)})dV_{H_{n}}  \notag \\
&=&[\rho ^{-(2n+1)}\hat{H}_{\phi }-(2n+1)\rho ^{-2n-2}\frac{\hat{e}_{I}\phi 
}{|d\phi |_{H_{n}}}\hat{e}_{I}\rho ]dV_{H_{n}}  \notag
\end{eqnarray}

\noindent where $\hat{H}_{\phi }$ denotes the horizontal ($p$-) mean
curvature with respect to the standard subriemannian metric in $H_{n}.$
Observe that $\frac{\hat{e}_{I}\phi }{|d\phi |_{H_{n}}}\hat{e}_{I}$ is the
horizontal normal to hypersurfaces defined by $\phi $ $=$ constant, denoted
as $e_{2n}^{\phi }.$ By (\ref{App-1}) and (\ref{App3}), we obtain%
\begin{equation}
H_{\phi }=\rho \hat{H}_{\phi }-(2n+1)e_{2n}^{\phi }\rho  \label{App4}
\end{equation}

\noindent (cf. Lemma 7.2 in \cite{chmy}). For $\phi $ $=$ $u(r)-z$ where $%
r=[\sum_{I=1}^{2n}(x^{I})^{2}]^{1/2},$ we want to get a formula for $H_{\phi
}$ in terms of $u$ and its derivatives. First we compute $|d\phi |_{H_{n}}$
as follows:%
\begin{eqnarray}
|d\phi |_{H_{n}}^{2} &=&\sum_{I=1}^{2n}(\hat{e}_{I}\phi )^{2}  \label{App5}
\\
&=&\sum_{I=1}^{2n}(u^{\prime }(r)\partial _{I}r-x^{I^{\prime }})^{2}  \notag
\\
&=&\sum_{I=1}^{2n}(u^{\prime }(r)\frac{x^{I}}{r}-x^{I^{\prime
}})^{2}=(u^{\prime }(r))^{2}+r^{2}  \notag
\end{eqnarray}

\noindent where $x^{I^{\prime }}$ $:=$ $x^{n+j}$ for $I$ $=$ $j$, $%
x^{I^{\prime }}$ $:=$ $-x^{j}$ for $I$ $=$ $n+j,$ $1\leq j\leq n$ and note
that $\sum_{I=1}^{2n}x^{I}x^{I^{\prime }}$ $=$ $0.$ Next noting that $%
\partial _{z}(\frac{\hat{e}_{I}\phi }{|d\phi |_{H_{n}}})$ $=$ $0,$ we compute%
\begin{eqnarray}
\hat{H}_{\phi } &=&\hat{e}_{I}(\frac{\hat{e}_{I}(\phi )}{|d\phi |_{H_{n}}})
\label{App6} \\
&=&\func{div}_{R^{2n}}(\frac{u^{\prime }(r)\nabla r-(x^{I^{\prime }})}{\sqrt{%
(u^{\prime }(r))^{2}+r^{2}}})  \notag
\end{eqnarray}

\noindent by (\ref{App5}). Observe that $\partial _{I}x^{I^{\prime }}=0$ and 
\begin{eqnarray*}
&&\sum_{I=1}^{2n}x^{I^{\prime }}\partial _{r}(\sqrt{(u^{\prime
}(r))^{2}+r^{2}})^{-1}\partial _{I}r \\
&=&r^{-1}\partial _{r}(\sqrt{(u^{\prime }(r))^{2}+r^{2}})^{-1}%
\sum_{I=1}^{2n}x^{I^{\prime }}x^{I}=0.
\end{eqnarray*}

\noindent So we can reduce (\ref{App6}) to%
\begin{eqnarray}
\hat{H}_{\phi } &=&\func{div}_{R^{2n}}(\frac{u^{\prime }(r)\nabla r}{\sqrt{%
(u^{\prime }(r))^{2}+r^{2}}})  \label{App7} \\
&=&\func{div}_{R^{2n}}(\frac{u^{\prime }(r)r^{2n-1}}{\sqrt{(u^{\prime
}(r))^{2}+r^{2}}}\frac{\nabla r}{r^{2n-1}})  \notag \\
&=&(\nabla \frac{u^{\prime }(r)r^{2n-1}}{\sqrt{(u^{\prime }(r))^{2}+r^{2}}}%
)\cdot \frac{\nabla r}{r^{2n-1}}\text{ (since }\func{div}_{R^{2n}}(\frac{%
\nabla r}{r^{2n-1}})=0)  \notag \\
&=&\frac{d}{dr}(\frac{u^{\prime }(r)r^{2n-1}}{\sqrt{(u^{\prime
}(r))^{2}+r^{2}}})\frac{1}{r^{2n-1}}\text{ (since }|\nabla r|^{2}=1).  \notag
\end{eqnarray}

\noindent On the other hand, we compute%
\begin{equation}
e_{2n}^{\phi }\rho =\frac{\hat{e}_{I}\phi }{|d\phi |_{H_{n}}}\hat{e}_{I}\rho
=\frac{r^{2}(ru^{\prime }(r)-2u(r))}{\rho ^{3}\sqrt{(u^{\prime
}(r))^{2}+r^{2}}}  \label{App8}
\end{equation}

\noindent for $\phi $ $=$ $u(r)-z.$ Substituting (\ref{App7}) and (\ref{App8}%
) into (\ref{App4}), we obtain%
\begin{eqnarray}
H_{\phi } &=&\frac{\rho }{r^{2n-1}}\frac{d}{dr}(\frac{u^{\prime }(r)r^{2n-1}%
}{\sqrt{(u^{\prime }(r))^{2}+r^{2}}})  \label{App9} \\
&&-(2n+1)\frac{r^{2}(ru^{\prime }(r)-2u(r))}{\rho ^{3}\sqrt{(u^{\prime
}(r))^{2}+r^{2}}}.  \notag
\end{eqnarray}

\noindent For $\phi $ $=$ $u(r)-z$ with $u(r)$ $=$ $cr^{2},$ $c$ being a
constant, we get 
\begin{eqnarray}
H_{\phi } &=&\frac{2(2n-1)c}{\sqrt{1+4c^{2}}}\frac{\rho }{r}  \label{App10}
\\
&=&\frac{2(2n-1)c}{(1+4c^{2})^{1/4}}  \notag
\end{eqnarray}

\noindent at points where $\phi $ $=$ $cr^{2}-z$ $=$ $0.$ That is to say,
the hypersurface in the Heisenberg cylinder $(H_{n}\backslash \{0\},$ $\rho
^{-2}\Theta ),$ defined by $z=cr^{2},$ has constant horizontal ($p$-) mean
curvature as shown in (\ref{App10}) (note that at $p_{0}$ $\in $ $\Sigma $ $%
:=$ $\{\phi $ $=$ $c\},$ we have $H_{\phi }(p_{0})$ $=$ $H_{\Sigma }(p_{0}),$
the horizontal ($p$-) mean curvature of $\Sigma $. See Proposition B.1 in
the Appendix). Next we want to compute the horizontal ($p$-) mean curvature
of Heisenberg spheres defined by $\rho ^{4}$ $=$ $c$ $>$ $0.$ Let $\phi $ $=$
$\rho ^{4}$ $-$ $c.$ We compute%
\begin{eqnarray*}
\hat{e}_{I}\phi &=&(\partial _{I}+x^{I^{\prime }}\partial
_{z})((r^{2})^{2}+4z^{2}) \\
&=&4(r^{2}x^{I}+2x^{I^{\prime }}z).
\end{eqnarray*}

\noindent It follows that 
\begin{eqnarray}
|d\phi |_{H_{n}}^{2} &=&\sum_{I=1}^{2n}(\hat{e}_{I}\phi )^{2}
\label{App10-1} \\
&=&16(r^{6}+4r^{2}z^{2})=16r^{2}\rho ^{4}.  \notag
\end{eqnarray}

\noindent Then a straightforward computation shows%
\begin{equation}
\hat{H}_{\phi }=(2n+1)\frac{r}{\rho ^{2}}.  \label{App11}
\end{equation}

\noindent On the other hand, since $4\rho ^{3}\hat{e}_{I}\rho $ $=$ $\hat{e}%
_{I}\phi $ for $\phi $ $=$ $\rho ^{4}$ $-$ $c,$ we have%
\begin{eqnarray}
e_{2n}^{\phi }\rho &=&\sum_{I=1}^{2n}\frac{\hat{e}_{I}\phi }{|d\phi |_{H_{n}}%
}\hat{e}_{I}\rho  \label{App12} \\
&=&\sum_{I=1}^{2n}\frac{(\hat{e}_{I}\phi )^{2}}{|d\phi |_{H_{n}}(4\rho ^{3})}
\notag \\
&=&\frac{|d\phi |_{H_{n}}}{4\rho ^{3}}=\frac{r}{\rho }  \notag
\end{eqnarray}

\noindent by (\ref{App10-1}). Substituting (\ref{App11}) and (\ref{App12})
into (\ref{App4}) gives%
\begin{equation*}
H_{\phi }=\rho (2n+1)\frac{r}{\rho ^{2}}-(2n+1)\frac{r}{\rho }=0.
\end{equation*}

\noindent So this means that Heisenberg spheres \{$\rho ^{4}$ $=$ $c\}$ are
horizontally ($p$-) minimal hypersurfaces in the Heisenberg cylinder $%
(H_{n}\backslash \{0\},$ $\rho ^{-2}\Theta ).$ We summarize what we obtain
so far as follows:

\bigskip

\textbf{Proposition 6.1.} \textit{Let }$\Sigma $\textit{\ be a hypersurface
in the Heisenberg cylinder }$(H_{n}\backslash \{0\},$\textit{\ }$\rho
^{-2}\Theta )$ \textit{with} $n\geq 1.$\textit{\ We have}

\textit{(a) Suppose }$\Sigma $\textit{\ is defined by }$z=cr^{2}$\textit{\
for a constant }$c.$\textit{\ Then }$\Sigma $\textit{\ is a hypersurface of
constant horizontal (}$p$\textit{-) mean curvature with constant }$\frac{%
2(2n-1)c}{(1+4c^{2})^{1/4}};$

\textit{(b) Suppose }$\Sigma $\textit{\ is defined by }$\rho ^{4}$\textit{\ }%
$=$\textit{\ }$c$\textit{\ for a constant }$c>0.$\textit{\ Then }$\Sigma $%
\textit{\ is a horizontally (}$p$\textit{-) minimal hypersurface.}

\bigskip

Observe that the dilation $\tau _{\lambda }$ $:$ $(x^{1},..,$ $x^{2n},$ $%
z)\rightarrow (\lambda x^{1},..,$ $\lambda x^{2n},$ $\lambda ^{2}z)$
preserves $\rho ^{-2}\Theta ,$ $i.e.,$ $\tau _{\lambda }^{\ast }(\rho
^{-2}\Theta )$ $=$ $\rho ^{-2}\Theta $ for any $\lambda $ $\in $ $%
R\backslash \{0\}.$ So $\tau _{\lambda }$ is a pseudohermitian isomorphism
of the Heisenberg cylinder\textit{\ }$(H_{n}\backslash \{0\},$\textit{\ }$%
\rho ^{-2}\Theta ).$ We can now prove a uniqueness result stated in Theorem
L in Section 1.

\bigskip

\proof
\textbf{(of Theorem L)} For the case (a), we take a Heisenberg sphere $%
S(c_{1})$ defined by $\rho ^{4}$\textit{\ }$=$\textit{\ }$c_{1}$ for $c_{1}$
large enough so that the interior region \{$\rho ^{4}$\textit{\ }$<$\textit{%
\ }$c_{1}\}$ of $S(c_{1})$ contains $\Sigma .$ Decrease (or take) $c_{1}$ to
reach a constant $c_{2}$ $>$ $0$ so that $S(c_{2})$ is tangent to $\Sigma $
at some point $p_{0}$ while $\Sigma $ lies in \{$\rho ^{4}$\textit{\ }$\leq $%
\textit{\ }$c_{2}\}.$ Observe that 
\begin{equation*}
H_{\Sigma }\leq 0=H_{S(c_{2})}
\end{equation*}

\noindent near $p_{0}$ by the assumption and Proposition 6.1 (b). It follows
from the SMP (Theorem C and Theorem J$^{\prime })$ that $\Sigma $ must
coincide with $S(c_{2}).$

Similarly for the case (b), we can first find a Heisenberg sphere $S(c_{3})$
defined by $\rho ^{4}$\textit{\ }$=$\textit{\ }$c_{3}$ for $c_{3}$ small
enough so that $S(c_{3})$ is contained in the interior region of $\Sigma .$
Increase \ (or take) $c_{3}$ to reach a constant $c_{4}$ $>$ $0$ so that $%
S(c_{4})$ is tangent to $\Sigma $ at some point $q$ while \{$\rho ^{4}$%
\textit{\ }$\leq $\textit{\ }$c_{4}\}$ is contained in the interior region
of $\Sigma .$ Now observe that%
\begin{equation*}
H_{\Sigma }\geq 0=H_{S(c_{4})}
\end{equation*}

\noindent near $q$ by the assumption and Proposition 6.1 (b). So it follows
that $\Sigma $ $=$ $S(c_{4})$ by the SMP.

\endproof%

\bigskip

We can also show a nonexistence result (pseudo-halfspace theorem).

\bigskip

\proof
\textbf{(of Theorem N) }For the first case $\Omega $ $=$ $\{z$\textit{\ }$>$%
\textit{\ }$\varphi (\sqrt{x_{1}^{2}+..+x_{2n}^{2}})\},$ we consider
comparison hypersurfaces -horizontal hyperplanes $\{z=c,$ a constant$\}.$
Starting from $c$ $=$ $c_{0}$ $<$ $\min_{\tau \in \lbrack 0,\infty )}\varphi
(\tau )$ (existence by the assumption: $\lim_{\tau \rightarrow \infty
}\varphi (\tau )$ = $\infty ),$ we increase $c$ to reach $c$ $=$ $c_{1}$
such that the hyperplane $\{z=c_{1}\}$ is tangent to $\Sigma $ at some point 
$p_{1}$ at the first time. The existence of such $p_{1}$ $\in $ $\{z=c_{1}\}$
$\cap $ $\Sigma $ is due to the immersion being proper. Note that the
hyperplane $\{z=c_{1}\}$ is horizontally ($p$-) minimal and its singular set
consists of one isolated singular point $(0,$ $..,$ $0,$ $c_{1})$. We can
then apply the SMP (Corollary I or Theorem J) at $p_{1}$ to conclude that $%
\Sigma $ $\subset $ the hyperplane $\{z=c_{1}\}$ which touches $\partial
\Omega $ $=$ $\{z$\textit{\ }$=$\textit{\ }$\varphi (\sqrt{%
x_{1}^{2}+..+x_{2n}^{2}})\}.$ But such $\Sigma $ is not properly immersed in 
$\Omega .$

For the second case $\Omega $ $=$ $\{x_{1}$\textit{\ }$>$\textit{\ }$\varphi
(\sqrt{x_{2}^{2}+..+x_{2n}^{2}+z^{2}})\},$ we consider comparison
hypersurfaces -vertical hyperplanes $\{x_{1}$ $=$ $c,$ a constant$\}.$ By a
similar reasoning as for the first case, we can find $c$ $=$ $c_{1}$ such
that the hyperplane $\{x_{1}$ $=$ $c_{1}\}$ is tangent to $\Sigma $ at some
point $q$. Note that the hyperplane $\{x_{1}$ $=$ $c_{1}\}$ is horizontally (%
$p$-) minimal and has no singular points. Apply the SMP (Corollary D) to
this situation to conclude that $\Sigma $ $\subset $ $\{x_{1}$ $=$ $c_{1}\}$
which touches $\partial \Omega $ $=$ $\{x_{1}$\textit{\ }$=$\textit{\ }$%
\varphi (\sqrt{x_{2}^{2}+..+x_{2n}^{2}+z^{2}})\}.$ But such $\Sigma $ is not
properly immersed in $\Omega .$

\endproof%

\bigskip

The simplest example is $\varphi (\tau )$ $=$ $a\tau $ with $a$ $>$ $0.$
Call associated domains wedge-shaped. Theorem N tells us nonexistence of $p$%
-minimal hypersurfaces in wedge-shaped domains.But Theorem N does not hold
for the case $a$ $=$ $0.$ That is, halfspace theorem does not hold since
there are catenoid type horizontal ($p$-) minimal hypersurfaces with finite
height (\cite{RR0}) in $H_{n}$ for $n$ $\geq $ $2.$ On the other hand, we do
have halfspace theorem for $H_{1}$ (see \cite{ch}).

\bigskip

\section{Appendix}

\textbf{A: Bony's strong maximum principle}

For completeness and reader's convenience, we review some material in Bony's
original paper \cite{Bony}.

Let $\Omega $ be a domain of $R^{n}.$ We consider a differential operator of
second order:%
\begin{equation*}
Lu(x)=\sum_{i,j=1}^{n}a_{ij}(x)\frac{\partial ^{2}u}{\partial x_{i}\partial
x_{j}}(x)+\sum_{i=1}^{n}a_{i}(x)\frac{\partial u}{\partial x_{i}}(x)+a(x)u(x)
\end{equation*}

\noindent with the following properties:

($\alpha )$ The quadratic form ($a_{ij}(x))$ is nonnegative for each $x$ $%
\in $ $\Omega ,$ that is,%
\begin{equation*}
\sum_{i,j=1}^{n}a_{ij}(x)\xi _{i}\xi _{j}\geq 0\text{ \ for any }x\in \Omega
,\text{ }\xi \in R^{n}.
\end{equation*}

($\beta )$ $a(x)\leq 0$ in $\Omega $ and $a(x)\in C^{\infty }.$

($\gamma )$ There exist vector fields $X_{1},...,X_{r}$ and $Y$ of class $%
C^{\infty }$ such that%
\begin{equation*}
Lu=\sum_{k=1}^{r}X_{k}^{2}u+Yu+au.
\end{equation*}

We remark that an operator $L$ with property $(\alpha )$ may not have
property $(\gamma ).$ Writing $X_{k}$ $=$ $b_{k}^{j}\frac{\partial }{%
\partial x_{j}}$ (summation convention), we compute%
\begin{equation*}
\sum_{k=1}^{r}X_{k}^{2}u=\sum_{k=1}^{r}b_{k}^{j}b_{k}^{i}\frac{\partial ^{2}u%
}{\partial x_{i}\partial x_{j}}+b_{k}^{j}\frac{\partial b_{k}^{i}}{\partial
x_{j}}\frac{\partial u}{\partial x_{i}}.
\end{equation*}

\noindent Therefore we have%
\begin{equation}
a_{ij}=\sum_{k=1}^{r}b_{k}^{j}b_{k}^{i}.  \label{A1}
\end{equation}

\textbf{Lemma A1}. $\sum_{i,j=1}^{n}a_{ij}(x)\xi _{i}\xi _{j}=0$\textit{\ if
and only if }$\xi =(\xi _{i})$\textit{\ is perpendicular to }$X_{k}$\textit{%
\ for any }$k,$\textit{\ }$1\leq k\leq r.$

\bigskip

\proof
By (\ref{A1}) we have 
\begin{eqnarray}
\sum_{i,j=1}^{n}a_{ij}(x)\xi _{i}\xi _{j}
&=&\sum_{k=1}^{r}\sum_{i,j=1}^{n}b_{k}^{j}b_{k}^{i}\xi _{i}\xi _{j}
\label{A2} \\
&=&\sum_{k=1}^{r}(\sum_{j=1}^{n}b_{k}^{j}\xi
_{j})(\sum_{i=1}^{n}b_{k}^{i}\xi _{i})  \notag \\
&=&\sum_{k=1}^{r}(X_{k}\cdot \xi )^{2}  \notag
\end{eqnarray}

\noindent where "$\cdot "$ means the standard inner product on $R^{n}.$ The
result follows from (\ref{A2}).

\endproof%

\bigskip

From (\ref{A2}) we also learn that $(\gamma )$ implies $(\alpha ).$

\bigskip

\textbf{Proposition A2 }(Proposition 1.1 in \cite{Bony}). \textit{Suppose }$%
u\in C^{2}(\Omega )$\textit{\ attains a nonnegative local maximum at }$%
x_{0}\in \Omega .$\textit{\ Then we have }$Lu(x_{0})$\textit{\ }$\leq $%
\textit{\ }$0.$\textit{\ If we further assume the maximum is positive and }$%
a(x_{0})$\textit{\ }$<$\textit{\ }$0.$\textit{\ Then we have }$Lu(x_{0})$%
\textit{\ }$<$\textit{\ }$0.$

\bigskip

\proof
Observe that the matrix ($u_{ij}(x_{0}))$ is nonpositive and the matrix $%
(a_{ij}(x_{0}))$ is nonnegative. It follows that%
\begin{equation*}
\sum_{i,j=1}^{n}a_{ij}(x_{0})u_{ij}(x_{0})=Trace[(a_{ij}(x_{0}))(u_{ij}(x_{0}))]\leq 0
\end{equation*}

\noindent by elementary linear algebra. We then have 
\begin{eqnarray*}
Lu(x_{0})
&=&\sum_{i,j=1}^{n}a_{ij}(x_{0})u_{ij}(x_{0})+%
\sum_{i=1}^{n}a_{i}(x_{0})u_{i}(x_{0})+a(x_{0})u(x_{0}) \\
&\leq &0
\end{eqnarray*}

\noindent in which we have used $u_{i}(x_{0})$ $=$ $0.$

\endproof%

\bigskip

Let $F$ be a closed subset of $\Omega .$ We say a vector $\nu $ is normal to 
$F$ at a point $x_{0}$ $\in $ $F$ if there exists a ball $B$ $\subset $ $%
\Omega \backslash F$ with center $x_{1},$ such that $x_{0}$ $\in $ $\partial
B$ and the vectors $x_{1}$ $-$ $x_{0}$ and $\nu $ are parallel.

We say a vector field $X$ is tangent to $F$ if for any $x_{0}$ $\in $ $F$
and any vector $\nu $ normal to $F$ at $x_{0},$ $X(x_{0})$ is perpendicular
to $\nu $.

\bigskip

\textbf{Theorem A3} (Theorem 2.1 in \cite{Bony}) \textit{Let }$\Omega $%
\textit{\ be a domain of }$R^{n}.$\textit{\ Let }$F$\textit{\ be a closed
subset of }$\Omega .$\textit{\ Suppose }$X$\textit{\ is a Lipschitz
continuous vector field on }$\Omega $\textit{, which is tangent to }$F.$%
\textit{\ Then any integral curve of }$X$\textit{\ meeting }$F$\textit{\ at
a point is entirely contained in }$F.$

\bigskip

\proof
(outline) Assume, on the contrary, there exists a curve $x(t)$ satisfying%
\begin{equation*}
\frac{dx(t)}{dt}=X(x(t)),
\end{equation*}

\noindent and meeting $F$ at a point, which is not contained in $F.$ Then we
can find an interval $[t_{0},t_{1}]$ such that $x(t_{0})$ $\in $ $F$ and $%
x(t)$ $\notin $ $F$ for all $t$ $\in $ $(t_{0},t_{1}].$ We then prove the
following two facts:

(a) (Lemma 2.1 in \cite{Bony}) Let $\delta (t)$ denote the distance between $%
x(t)$ and $F.$ Then there exists a constant $K_{0}$ $>$ $0$ such that for
all $t$ $\in $ $(t_{0},t_{1}]$, there holds%
\begin{equation*}
\lim \inf_{h\rightarrow 0}\frac{\delta (t+h)-\delta (t)}{|h|}\geq
-K_{0}\delta (t).
\end{equation*}

(b) (Lemma 2.2 in \cite{Bony}) Suppose $\varphi (t)$ is a continuous
function on $[t_{0},t_{1}]$ and satisfies%
\begin{equation*}
\lim \inf_{h\rightarrow 0}\frac{\varphi (t+h)-\varphi (t)}{|h|}\geq -M
\end{equation*}

\noindent for all $t$ $\in $ $(t_{0},t_{1}).$ Then $\varphi $ is
Lipschitzian of ratio $M$ on $[t_{0},t_{1}].$

Now let $\theta $ $=$ $\inf \{t_{1}-t_{0},$ $\frac{1}{2K_{0}}\}$ and $%
\varepsilon $ $=$ $\sup_{t\in \lbrack t_{0},t_{0}+\theta ]}\delta (t).$ From 
$(a),$ $(b)$ we obtain that the function $\delta $ is Lipschitzian of ratio $%
K_{0}\varepsilon $ on $[t_{0},t_{0}+\theta ].$ So we have%
\begin{eqnarray}
|\delta (t)-\delta (t_{0})| &\leq &\varepsilon K_{0}(t-t_{0})  \label{A3} \\
&\leq &\varepsilon K_{0}\theta \leq \frac{\varepsilon }{2}  \notag
\end{eqnarray}

\noindent for any $t$ $\in $ $[t_{0},t_{0}+\theta ].$ On the other hand, we
observe that $\delta (t_{0})$ $=$ $0$ in (\ref{A3}) and hence $\varepsilon $ 
$=$ $\sup_{t\in \lbrack t_{0},t_{0}+\theta ]}\delta (t)=\sup_{t\in \lbrack
t_{0},t_{0}+\theta ]}|\delta (t)-\delta (t_{0})|\leq \frac{\varepsilon }{2},$
a contradiction.

\endproof%

\bigskip

Let us denote by $\pounds (X_{1},...,X_{r})$ the smallest $C^{\infty }$%
-module which contains $X_{1},$ $...,$ $X_{r}$ and is closed under the Lie
bracket. That is, if $Z$ $\in $ $\pounds (X_{1},...,X_{r}),$ then $Z$ is a
sum of finite terms of the form%
\begin{equation}
\lambda \lbrack X_{i_{1}},[X_{i_{2}},...,[X_{i_{l-1}},X_{i_{l}}]]]
\label{A4}
\end{equation}

\noindent where $\lambda $ $\in $ $C^{\infty }(\Omega )$ and $i_{k}$ $\in $
\{$1,$ $...,r\}.$ The rank of $\pounds (X_{1},...,X_{r})$ at a point $x$ is
the dimension of the vector space spanned by the vectors $Z(x)$ for all $Z$ $%
\in $ $\pounds (X_{1},...,X_{r}).$

\bigskip

\textbf{Proposition A4 }(Proposition 2.1 in \cite{Bony}). \textit{Let }$%
X_{1},...,X_{r}$\textit{\ be vector fields of class }$C^{\infty }.$\textit{\
Take }$Z$\textit{\ }$\in $\textit{\ }$\pounds (X_{1},...,X_{r}).$\textit{\
Then any integral curve of }$Z$\textit{\ can be uniformly approximated by
piecewise differentiable curves each of which is an integral curve of one of
vector fields }$X_{1},$\textit{\ }$...,$\textit{\ }$X_{r}.$

\bigskip

\textbf{Theorem A5 }(Theorem 2.2 in \cite{Bony}). \textit{Let }$\Omega $%
\textit{\ be an open set in }$R^{n}$\textit{. Let }$F$\textit{\ be a closed
subset of }$\Omega .$\textit{\ Suppose }$X_{1},...,X_{r}$\textit{\ are
vector fields of class }$C^{\infty },$\textit{\ each of which is tangent to }%
$F.$\textit{\ Then for each }$Z$\textit{\ }$\in $\textit{\ }$\pounds %
(X_{1},...,X_{r}),$\textit{\ }$Z$\textit{\ is tangent to }$F$\textit{\ and
any integral curve of }$Z$\textit{\ meeting }$F$\textit{\ at a point is
entirely contained in }$F.$

\bigskip

We refer the reader to the original paper of Bony for the proof of the above
two results. We then discuss the propagation of maximums.

\bigskip

\textbf{Proposition A6} (Proposition 3.1 in \cite{Bony}) \textit{Let }$u$%
\textit{\ be a function of class }$C^{2}$\textit{\ on }$\Omega $\textit{\
such that }$Lu$\textit{\ }$\geq $\textit{\ }$0.$\textit{\ Suppose the
maximum of }$u$\textit{\ is nonnegative and attained at a point in }$\Omega
. $\textit{\ Let }$F$\textit{\ be the set of all points where }$u$\textit{\
attains the maximum. Then for each }$k$\textit{\ }$=$\textit{\ }$1,$\textit{%
\ }$...,$\textit{\ }$r,$\textit{\ the vector field }$X_{k}$\textit{\ is
tangent to }$F.$

\bigskip

\proof
Let $B(x_{0},\rho )$ $\subset $ $\Omega \backslash F$ be a ball of radius $%
\rho $ and center $x_{0},$ such that $\partial B(x_{0},\rho )$ $\cap $ $F$ $%
= $ $\{x_{1}\}.$ We are going to show 
\begin{equation*}
\alpha :=\sum_{i,j}a_{ij}(x_{1})(x_{0}^{i}-x_{1}^{i})(x_{0}^{j}-x_{1}^{j})=0.
\end{equation*}

\noindent The conclusion follows from Lemma A1 and the definition of tangent
vector. Suppose, on the contrary, $\alpha $ $>$ $0.$ Consider an auxiliary
function $v$ defined by%
\begin{equation*}
v(x)=e^{-k|x-x_{0}|^{2}}-e^{-k\rho ^{2}}
\end{equation*}

\noindent where $k$ is a positive constant. A direct computation gives%
\begin{eqnarray*}
Lv(x_{1}) &=&e^{-k\rho
^{2}}[4k^{2}\sum_{i,j}a_{ij}(x_{1})(x_{0}^{i}-x_{1}^{i})(x_{0}^{j}-x_{1}^{j})
\\
&&-2k(a_{ii}(x_{1})+a_{i}(x_{1})(x_{1}^{i}-x_{0}^{i}))] \\
&=&e^{-k\rho ^{2}}[4k^{2}\alpha
-2k(a_{ii}(x_{1})+a_{i}(x_{1})(x_{1}^{i}-x_{0}^{i}))]
\end{eqnarray*}

\noindent in which we have used $v(x_{1})$ $=$ $0.$ Therefore if $k$ is
large enough, we have%
\begin{equation*}
Lv(x_{1})>0,
\end{equation*}

\noindent and hence $Lv$ $>$ $0$ in a neighborhood $V$ of $x_{1}.$ We now
consider the function%
\begin{equation*}
w(x)=u(x)+\lambda v(x),\text{ }\lambda >0.
\end{equation*}

\noindent It follows that $Lw$ $=$ $Lu+\lambda Lv$ $>$ $0$ in $V$ since $Lu$ 
$\geq $ $0$ by assumption. From Proposition A2 we conclude that $w$ cannot
achieve a nonnegative local maximum in $V.$ On the other hand, we are going
to show $w|_{\partial V}$ $\leq $ $w(x_{1}),$ a contradiction.

Let $B(x_{0},\rho )^{c}$ denote the complement of $B(x_{0},\rho ).$ Observe
that for $x$ $\in $ $\partial V$ $\cap $ $B(x_{0},\rho )^{c},$ $v(x)$ $\leq $
$0$ and hence%
\begin{equation*}
w(x)\leq u(x)\leq u(x_{1})=w(x_{1})
\end{equation*}

\noindent since $u$ attains the maximum at $x_{1}$ and $v(x_{1})$ $=$ $0.$
Now observe that $\sup_{\partial V\cap B(x_{0},\rho )}u$ $<$ $u(x_{1})$ by
compactness$.$ So for $x$ $\in $ $\partial V$ $\cap $ $B(x_{0},\rho ),$we
have%
\begin{equation*}
w(x)\leq u(x_{1})=w(x_{1})
\end{equation*}

\noindent for small $\lambda .$

\endproof%

\bigskip

\textbf{Theorem A7} (Theorem 3.1 in \cite{Bony}). \textit{Let }$u$\textit{\
be a function of class }$C^{2}$\textit{\ in }$\Omega $\textit{\ and }$Z$%
\textit{\ }$\in $\textit{\ }$\pounds (X_{1},...,X_{r}).$\textit{\ Suppose }$%
Lu$ $\geq $ $0$\textit{\ and }$u$\textit{\ attains a nonnegative maximum at
a point of an integral curve }$\Gamma $\textit{\ of }$Z.$\textit{\ Then the
maximum is attained at all points of }$\Gamma .$

\bigskip

\textbf{Corollary A8} (Corollary 3.1 in \cite{Bony}) \textit{Suppose the
rank of }$\pounds (X_{1},...,X_{r})$\textit{\ is }$n$\textit{\ for all
points. Then a function of class }$C^{2}$\textit{\ in }$\Omega $\textit{\
satisfying }$Lu$\textit{\ }$\geq $\textit{\ }$0$\textit{\ cannot achieve a
nonnegative maximum in }$\Omega $\textit{\ unless it is constant.}

\bigskip

\textbf{B: Subriemannian geometry from the viewpoint of differential forms}

A subriemannian manifold is a ($C^{\infty })$ smooth manifold $M$ equipped
with a nonnegative definite inner product $<\cdot ,\cdot >^{\ast }$ on $%
T^{\ast }M,$ its cotangent bundle. Clearly if $<\cdot ,\cdot >^{\ast }$ is
positive definite, $(M,<\cdot ,\cdot >^{\ast })$ is a Riemannian manifold.
For $M$ being the Heisenberg group $H_{n}$ of dimension $m$ $=$ $2n+1,$ we
recall that the multiplication $\circ $ of $H_{n}$ reads%
\begin{eqnarray*}
&&(a_{1},..,a_{n},b_{1},..,b_{n},c)\circ (x_{1},..,x_{n},y_{1},..,y_{n},z) \\
&=&(a_{1}+x_{1},...,b_{n}+y_{n},c+z+\sum_{j=1}^{n}(b_{j}x_{j}-a_{j}y_{j})).
\end{eqnarray*}%
\noindent Let 
\begin{equation*}
\hat{e}_{j}=\frac{\partial }{\partial x_{j}}+y_{j}\frac{\partial }{\partial z%
},\hat{e}_{j^{\prime }}=\frac{\partial }{\partial y_{j}}-x_{j}\frac{\partial 
}{\partial z},
\end{equation*}

\noindent $1\leq j\leq n$ be the left-invariant vector fields on $H_{n},$ in
which $x_{1},$ $..,$ $x_{n},$ $y_{1},$ $..,$ $y_{n},$ $z$ denote the
coordinates of $H_{n}$ (instead of $x^{1},$ $..,$ $x^{n},$ $x^{n+1},..,$ $%
x^{2n},$ $z$ used previously)$.$ The (contact) 1-form $\Theta $ $\equiv $ $%
dz+\sum_{j=1}^{n}(x_{j}dy_{j}-y_{j}dx_{j})$ annihilates $\hat{e}_{j}^{\prime
}s$ and $\hat{e}_{j^{\prime }}^{\prime }s.$ We observe that $dx_{1},$ $%
dy_{1},$ $dx_{2},$ $dy_{2},$ $...,$ $dx_{n},$ $dy_{n},$ $\Theta $ are dual
to $\hat{e}_{1},$ $\hat{e}_{1^{\prime }},$ $\hat{e}_{2},$ $\hat{e}%
_{2^{\prime }},$ ...,$\hat{e}_{n},$ $\hat{e}_{n^{\prime }},$ $\frac{\partial 
}{\partial z}.$ Define a nonnegative inner product $<\cdot ,\cdot >_{H_{n}}$
or $<\cdot ,\cdot >^{\ast }$ by%
\begin{eqnarray}
&<&dx_{j},dx_{k}>^{\ast }=\delta _{jk},\text{ }<dy_{j},dy_{k}>^{\ast
}=\delta _{jk},<dx_{j},dy_{k}>^{\ast }=0,  \label{H1} \\
&<&\Theta ,dx_{j}>^{\ast }=<\Theta ,dy_{k}>^{\ast }=<\Theta ,\Theta >^{\ast
}=0.  \notag
\end{eqnarray}

We can extend the definition of the above nonnegative inner product to the
situation of a general pseudohermitian manifold. Take $e_{j},$ $e_{j^{\prime
}}=Je_{j}$, $j=1,$ $2,$ $...,n$ to be an orthonormal basis in the kernel of
the contact form $\Theta $ with respect to the Levi metric $\frac{1}{2}%
d\Theta (\cdot ,J\cdot ).$ Let $T$ be the Reeb vector field of $\Theta $
(such that $\Theta (T)$ $=$ $1$ and $d\Theta (T,\cdot )$ $=$ $0).$ Denote
the dual coframe of $e_{j},$ $e_{j^{\prime }},$ $T$ by $\theta ^{j},$ $%
\theta ^{j^{\prime }}$ (and $\Theta ).$ Now we can replace $dx_{j}$, $dy_{j}$
by $\theta ^{j},$ $\theta ^{j^{\prime }}$ in (\ref{H1}) to define a
nonnegative inner product on a general pseudohermitian manifold:

\begin{eqnarray}
&<&\theta ^{j},\theta ^{k}>^{\ast }=\delta _{jk},\text{ }<\theta ^{j^{\prime
}},\theta ^{k^{\prime }}>^{\ast }=\delta _{jk},<\theta ^{j},\theta
^{k^{\prime }}>^{\ast }=0,  \label{H2} \\
&<&\Theta ,\theta ^{j}>^{\ast }=<\Theta ,\theta ^{k^{\prime }}>^{\ast
}=<\Theta ,\Theta >^{\ast }=0.  \notag
\end{eqnarray}

Define the bundle morphism $G:T^{\ast }M\rightarrow TM$ by 
\begin{equation}
\omega (G(\eta ))=<\omega ,\eta >^{\ast }  \label{H3}
\end{equation}

\noindent for $\omega ,\eta $ $\in $ $T^{\ast }M.$ In the Riemannian case, $%
G $ is in fact an isometry. In the pseudohermitian case, $G(T^{\ast }M)$ is
the contact subbundle $\xi $ of $TM,$ the kernel of $\Theta .$ By letting $%
\eta $ $=$ $\Theta $ in (\ref{H3})$,$ we get $G(T^{\ast }M)$ $\subset $ $\xi
.$ On the other hand, it is easy to see that $G(\theta ^{j})$ $=$ $e_{j},$ $%
G(\theta ^{j^{\prime }})$ $=$ $e_{j^{\prime }}$ (and $G(\Theta )$ $=$ $0)$.
Since $e_{j},$ $e_{j^{\prime }}$, $j=1,$ $2,$ $...,n$ span $\xi ,$ we have $%
\xi $ $\subset $ $G(T^{\ast }M).$ For a smooth function $\varphi $ on $M,$
we define the gradient $\nabla \varphi $ $:=$ $G(d\varphi ).$ In the
pseudohermitian case, this $\nabla \varphi $ is nothing but the subgradient $%
\nabla _{b}\varphi $ $:=$ $\sum_{j=1}^{n}\{e_{j}(\varphi )e_{j}$ $+$ $%
e_{j^{\prime }}(\varphi )e_{j^{\prime }}\}.$

Let $\Sigma $ $\subset $ $M$ be a (smooth) hypersurface in $M$ with a
defining function $\phi $ such that $\Sigma $ $=$ $\{\phi =0\}.$ Fix a
background volume form $dv_{M},$ i.e., a nonvanishing $m+1$ form on $M,$
where $m+1$ $=$ $\dim M$ (hence $\dim \Sigma =m).$ For a point $\zeta $
where \TEXTsymbol{\vert}$d\phi |_{\ast }^{2}:=<d\phi ,d\phi >$ $\neq $ $0,$
for any $p$ $\geq $ $0,$ we define subriemannian area (or volume) element $%
dv_{\phi ,p}$ and (generalized) mean curvature $H_{\phi ,p}(\zeta )$ for the
hypersurface $\{\phi $ $=$ $\phi (\zeta )\}$ by%
\begin{equation*}
dv_{\phi ,p}:=\frac{d\phi }{|d\phi |_{\ast }^{1-p}}\lrcorner \text{ }dv_{M},
\end{equation*}

\begin{equation}
d(dv_{\phi ,p})=d(\frac{d\phi }{|d\phi |_{\ast }^{1-p}}\lrcorner \text{ }%
dv_{M}):=H_{\phi ,p}dv_{M},  \label{3.1}
\end{equation}

\noindent respectively. Here the interior product for forms is defined so
that%
\begin{equation*}
\eta \wedge (\omega \lrcorner \text{ }dv_{M})=<\eta ,\omega >^{\ast }dv_{M}.
\end{equation*}

\noindent The above notion of subriemannian area unifies those in Riemannian
and pseudohermitian geometries. For more details (on the case $p$ $=$ $0)$,
please see the Appendix: Generalized Heisenberg Geometry in \cite{ch1}.
Define unit normal $\nu _{p}$ to a hypersurface $\Sigma $ $:=$ $\{\phi $ $=$ 
$c\}$ by the formula%
\begin{equation}
\nu _{p}\lrcorner \text{ }dv_{M}=dv_{\phi ,p}.  \label{3.1.1}
\end{equation}%
\noindent Given area element $dv_{\phi ,p}$ and unit normal $\nu _{p},$ we
can also define (generalized) mean curvature $H_{\Sigma ,p}$ through a
variational formula:%
\begin{equation}
\delta _{f\nu _{p}}\int_{\Sigma }dv_{\phi ,p}=\int_{\Sigma }fH_{\Sigma
,p}dv_{\phi ,p}  \label{3.2}
\end{equation}

\noindent for $f$ $\in $ $C_{0}^{\infty }(\Sigma ).$

\bigskip

\textbf{Proposition B.1}. \textit{On }$\Sigma ,$ \textit{we have }$H_{\phi
,p}$ $=$ $H_{\Sigma ,p}.$

\bigskip

\proof
Let $\iota _{\mu }$ denote the interior product with vector $\mu .$ From (%
\ref{3.2}) we compute%
\begin{eqnarray*}
\int_{\Sigma }fH_{\Sigma ,p}dv_{\phi ,p} &=&\delta _{f\nu _{p}}\int_{\Sigma
}dv_{\phi ,p}=\int_{\Sigma }L_{f\nu _{p}}(dv_{\phi ,p}) \\
&=&\int_{\Sigma }(d\circ \iota _{f\nu _{p}}+\iota _{f\nu _{p}}\circ
d)(dv_{\phi ,p}) \\
&=&\int_{\Sigma }d(f\nu _{p}\lrcorner \text{ }dv_{\phi ,p})+\iota _{f\nu
_{p}}(d(dv_{\phi ,p})) \\
&=&\doint\limits_{\partial \Sigma }f\nu _{p}\lrcorner \text{ }dv_{\phi
,p}+\int_{\Sigma }\iota _{f\nu _{p}}(H_{\phi ,p}dv_{M})\text{ (by (\ref{3.1}%
))} \\
&=&0+\int_{\Sigma }fH_{\phi ,p}dv_{\phi ,p}
\end{eqnarray*}

\noindent since $f$ $\in $ $C_{0}^{\infty }(\Sigma )$ and by (\ref{3.1.1})$.$
It follows that $H_{\Sigma ,p}$ $=$ $H_{\phi ,p}.$

\endproof%

\bigskip

Let $(M,J,\Theta )$ be a pseudohermitian manifold of dimension $2n+1,$
considered as a subriemannian manifold$.$ We take the background volume form 
$dv_{M}$ $:=$ $\Theta \wedge (d\Theta )^{n}.$ Let $\Sigma $ be a (smooth)
hypersurface of $M.$ At nonsingular points (where $T\Sigma $ is transversal
to the contact bundle $\xi ),$ we choose orthonormal (with respect to the
Levi metric $<\cdot ,\cdot >_{Levi}$ $:=$ $\frac{1}{2}d\Theta (\cdot ,J\cdot
))$ basis $e_{1},$ $e_{1^{\prime }},$ $...$, $e_{n-1},$ $e_{(n-1)^{\prime
}}, $ $e_{n}$ in $T\Sigma \cap \xi ,$ where $e_{j^{\prime }}$ $=$ $Je_{j}.$
We choose the horizontal (or Legendrian) normal $\nu $ $=$ $e_{n^{\prime }}$
and a defining function $\phi $ satisfying $d\phi (\nu )$ $>$ $0$ such that
the $p$-area element $dv_{\phi }:=\frac{d\phi }{|d\phi |_{\ast }}\lrcorner $ 
$dv_{M}$ has the expression $\Theta \wedge $ $e^{1}\wedge e^{1^{\prime }}$ $%
\wedge ...$ $e^{n-1}$ $\wedge $ $e^{(n-1)^{\prime }}$ $\wedge $ $e^{n}$
where $\Theta ,$ $e^{1},$ $e^{1^{\prime }},$ $...,$ $e^{n-1}$, $%
e^{(n-1)^{\prime }}$, $e^{n},$ $e^{n^{\prime }}$ are dual to $T,$ $e_{1},$ $%
e_{1^{\prime }},$ $...$, $e_{n-1},$ $e_{(n-1)^{\prime }},$ $e_{n},$ $%
e_{n^{\prime }}.$ It follows that the associated mean curvature $H_{\Sigma }$%
, called horizontal mean curvature, is the trace of the second fundamental
form:%
\begin{equation*}
H_{\Sigma }=\sum_{K=1,1^{\prime },..,n}<\nabla _{e_{K}}^{p.h.}\nu
,e_{K}>_{Levi}
\end{equation*}

\noindent (see the Appendix in \cite{C} for more details; note sign
difference in \cite{C}). Let $|X|_{Levi}$ $:=$ $<X,X>_{Levi}^{1/2}$ for $X$ $%
\in $ $\xi .$ We can also express $H_{\Sigma }$ in terms of $\phi $ as
follows:%
\begin{equation*}
H_{\Sigma }=\func{div}_{b}\frac{\nabla _{b}\phi }{|\nabla _{b}\phi |_{levi}}
\end{equation*}

\noindent where $\nabla _{b}$ and $\func{div}_{b}$ denote subgradient and
subdivergence in $(M,J,\Theta ).$ For a graph $\Sigma $ $:=$ $\{(x_{1},$ $%
.., $ $x_{n},$ $y_{1},$ $..,$ $y_{n},$ $u(x_{1},$ $..,$ $x_{n},$ $y_{1},$ $%
..,$ $y_{n}))\}$ in the Heisenberg group $H_{n},$ we take the defining
function $\phi (x_{1},$ $..,$ $x_{n},$ $y_{1},$ $..,$ $y_{n},$ $z)$ $=$ $%
u(x_{1},$ $.., $ $x_{n},$ $y_{1},$ $..,$ $y_{n})$ $-$ $z.$ Then the
resulting $H_{\Sigma }$ coincides with the definition given by (\ref{0.1})
with $\vec{F}$ $:=$ $(-y_{1},x_{1},$ $..,$ $-y_{n},$ $x_{n})$.

\end{document}